\documentclass[opre,nonblindrev]{informs_bl}
\usepackage{amsmath,amssymb}
\usepackage{enumerate}
\usepackage{framed}

\def\Neil#1{{{\color{red}\texttt{ #1}}}}

\usepackage{subfigure}
\usepackage{tabularx}
\usepackage{multirow}
\usepackage{comment}
\usepackage{bigints}
\usepackage{algorithm}
\usepackage{algpseudocode}
\usepackage{lineno}
\usepackage[pagebackref=false, bookmarks]{hyperref}
\hypersetup{
backref,
pageanchor=true,
bookmarks,
bookmarksnumbered,
colorlinks=true,
menubordercolor= [rgb]{0, 0, 1},
linkbordercolor= [rgb]{0, 1, 1},
citecolor = red,
linkcolor= blue}
\usepackage{url}

\OneAndAHalfSpacedXI 

\usepackage{endnotes}
\let\footnote=\endnote
\let\enotesize=\normalsize
\def\notesname{Endnotes}%
\def\makeenmark{\hbox to1.275em{\theenmark.\enskip\hss}}
\def\enoteformat{\rightskip0pt\leftskip0pt\parindent=1.275em
  \leavevmode\llap{\makeenmark}}

\newcommand{\expct}[1]{\ensuremath{\text{{\bf E}$\left[#1\right]$}}}
\newcommand{\ord}[1]{\ensuremath{\mathcal{O}\left(#1\right)}}

\newcommand{\ceil}[1]{\ensuremath{\left\lceil#1\right\rceil}}
\newcommand{\floor}[1]{\ensuremath{\left\lfloor#1\right\rfloor}}

\newcommand{\junk}[1]{}

\newcommand{\ignore}[1]{}

\newcommand{\E}{\mathbb{E}} 

\usepackage[normalem]{ulem}

\usepackage{natbib}
 \bibpunct[, ]{(}{)}{,}{a}{}{,}%
 \def\bibfont{\small}%
 \def\bibsep{\smallskipamount}%
 \def\bibhang{24pt}%
\TheoremsNumberedThrough     
\ECRepeatTheorems

\EquationsNumberedThrough    

\newcommand{\bP}{\mathbb{P}}
\newcommand{\bE}{\mathbb{E}}

\newcommand{\bZ}{\mathbb{Z}}

\newcommand{\bI}{{\mathbf I}}

\newcommand{\Ik}{I^{(k)}}
\newcommand{\Nk}{N^{(k)}}
\newcommand{\Mk}{M^{(k)}}

\newcommand{\Vk}{V^{(k)}}

\graphicspath{{./figures/}}

\begin{document}

\RUNAUTHOR{Gupta and Walton}
\RUNTITLE{JSQ in NDS}

\TITLE{Load Balancing in the Non-Degenerate Slowdown Regime}
\ARTICLEAUTHORS{%
\AUTHOR{Varun Gupta}
\AFF{Booth School of Business, University of Chicago, Chicago, IL 60637, \EMAIL{varun.gupta@chicagobooth.edu}} 
\AUTHOR{Neil Walton}
\AFF{ Alan Turing Building, University of Manchester, Manchester
M13 9PL, UK,
 \EMAIL{neil.walton@manchester.ac.uk}}
} 


\ABSTRACT{%
We analyse Join-the-Shortest-Queue in a contemporary scaling regime known as the Non-Degenerate Slowdown regime. 
Join-the-Shortest-Queue (JSQ) is a classical load balancing policy for queueing systems with multiple parallel servers.
Parallel server queueing systems are regularly analysed and dimensioned by diffusion approximations achieved in the Halfin-Whitt scaling regime. 
However, when jobs must be dispatched to a server upon arrival,
we advocate the Non-Degenerate Slowdown regime (NDS) to compare different load-balancing rules. 

In this paper we identify novel diffusion approximation and timescale separation that provides insights into the performance of JSQ.
We calculate the price of irrevocably dispatching jobs to servers and prove this to be within 15\% (in the NDS regime) of the rules that may manoeuvre jobs between servers. 
We also compare ours results for the JSQ policy with the NDS approximations of many modern load balancing policies such as Idle-Queue-First and Power-of-$d$-choices policies which act as low information proxies for the JSQ policy.
Our analysis leads us to construct new rules that have identical performance to JSQ but require less communication overhead than power-of-2-choices.
}

\KEYWORDS{Non-Degenerate Slowdown, Load balancing, Join-the-Shortest-Queue, $M/M/k$}

\maketitle


\section{Introduction} \label{sec:introduction}

We investigate a parallel server queueing system using the Join-the-Shortest-Queue (JSQ) dispatch rule.
Join-the-Shortest-Queue is a simple, classical and widely deployed load balancing policy in Operations Research. 
It remains of central importance, in part, due to its known optimality properties.
Given the obstacles in implementing the JSQ rule, there is currently substantial interest in load balancing policies that approximate JSQ rule without requiring global state information. 
This is relevant for data center applications where communication between dispatcher and servers is expensive.
Two of the most popular policies in the class of low-overhead load balancing policies are the Power-of-$d$ choices policies (Po$d$) (see \cite{mitzenmacher2001power} and \cite{vvedenskaya1996queueing}), and the Idle-Queue-First rule (IQF), see \cite{lu2011join}.
In this paper we quantify the optimality of JSQ and compare it with other well known dispatch policies.
From this, we design simple dispatch rules which have asymptotically optimal performance with respect to JSQ.
%
%

To achieve this goal, we are obliged to investigate a non-standard diffusion limit. Parallel queueing models are commonly studied with fluid limits, heavy traffic limits, mean field limits and in the Halfin-Whitt regime. 
However as we discuss (in Section \ref{sec:OtherAsymp}), these asymptotic frameworks are unable to distinguish between several well-known dispatch policies. 
We show that the recently introduced \emph{Non-Degenerate Slowdown} regime (NDS) (see \cite{atar2012diffusion}) emerges naturally from our goal to delineate the performance of load balancing policies.

In this paper, for a parallel queueing system with arrival rate $\lambda$ and $k$ servers each with service rate $\mu$, the \emph{Non-Degenerate~Slowdown} regime is defined as a many server asymptotic framework where, for some constant $\alpha$, 
\begin{equation}\label{Intro:NDS}
\lambda -k \mu \rightarrow -\alpha \mu \qquad \text{as}\qquad k\rightarrow\infty.
\end{equation}
In other words, as the number of servers $k$ grows, the number of ``spare servers,'' $\alpha$, remains fixed. (The above definition may appear somewhat different than the definition by \cite{atar2012diffusion}, but is equivalent up to scaling of time by some function of $k$.)
To explain the moniker, recall that the \emph{slowdown} of a job is the ratio between its sojourn time and its service requirement, and is greater than or equal to $1$ by definition. (A slowdown equal to $1$ corresponds to uninterrupted service). As described in \cite{atar2012diffusion},  an asymptotic regime is considered to be degenerate if the slowdown distribution of a typical customer either converges to $1$ or diverges, else the regime has \emph{non-degenerate slowdown}. 
We show that, for JSQ, Non-Degenerate Slowdown is the asymptotic regime given by \eqref{Intro:NDS}.

We analyze JSQ and related dispatch policies when arrivals are Poisson and job sizes are exponentially distributed. For this simple Markov chain model a rich set of probabilistic techniques including fluid and diffusion approximations, couplings, time-scale separation, reversibility and stochastic averaging principles are required. We characterize the limiting diffusion of JSQ and related policies in the NDS regime. 

Our main results provide new insights into Join-the-Shortest-Queue. These contributions are briefly summarized as follows:

\begin{enumerate}
		\item \emph{We provide a novel characterization of JSQ:}
 We provide the first analysis of JSQ under the NDS regime, characterizing its diffusion approximation and stationary distribution. 
	\item \emph{We quantify the difference between optimal dispatch and optimal pooling:} As we will review shortly, Join-the-shortest queue is the optimal size agnostic dispatch rule, while maintaining a single centralized queue provides optimal workload based dispatch. In NDS, we can precisely quantify the optimality gap between these two rules.  We show that in the NDS regime, JSQ has a mean response time that is at most 15\% larger than Centralized Queueing scheme (CQ) used by an $M/M/k$ queue. In essence, this quantifies the impact of not being able to pool jobs or jockey jobs between queues.
	\item \emph{We provide a new low-information dispatch rule achieving JSQ optimality:} The Idle-Queue-First policy, which is used as a low-overhead proxy for JSQ, can have mean response time up to a 100\% larger than Central Queue $M/M/k$. However, we show that a minor modification that prioritizes idle servers first and servers with one job will lead to the same asymptotic performance as JSQ for the distribution of total number of jobs in system, and hence for mean response time. We call this policy Idle-One-First (I1F).
\end{enumerate}

More broadly, we find that, unlike conventional limit regimes, NDS is able to distinguish between different load-balancing rules. For instance, we remark that the 2nd bullet point cannot be observed in Heavy Traffic, the two processes of interest being indistinguishable. Similarly for bullet 3 under a Halfin Whitt regime the two processes there are also indistinguishable. 
We should note that we do not advocate NDS as a regime by which one should provision capacity under a given load. We recommend \cite{borst2004dimensioning} for such an analysis.
However, when a system experiences load that induces fluctuations on the order of magnitude of a system's size or when we wish to compare performance of load balancing policies, then NDS is a meaningful regime to gain qualitative insight into system performance.

\subsection{Outline}
The remainder of the paper proceeds as follows. 
In Section \ref{sec:model}, we define our parallel queueing model; we introduce the load balancing policies of interest; we define and motivate the NDS Regime; and, we compare NDS with several well known asymptotics.
In Section \ref{sec:literature}, we review relevant literature on parallel server models, policies, NDS and other asymptotic regimes.
In Section \ref{sec:NDS1}, we state our main result on JSQ. 
Proofs and technical results are deferred to Appendix. 
In Section \ref{sec:NDS2}, we present our analysis of the other polices introduced in Section~\ref{sec:model}:  Central Queue, Idle-Queue-First and Idle-One-First. 
%
We conclude in Section \ref{sec:conclusions} by discussing various avenues of potential future investigation.

\section{Parallel Server Model} \label{sec:model}

We consider a queueing system that consists of $k$ servers each serving at unit rate. 
Jobs arrive to this queueing system as a Poisson process of rate $\lambda$.
Each job has a service requirement that is independent exponentially distributed with rate parameter $\mu$. 
Upon arrival a job is dispatched to one of $k$ servers where it is queued and receives service. The assignment of jobs to servers is irrevocable, and the jobs can not jockey between queues after being dispatched.
Thus the load on each queue is 
\[
\rho =\frac{\lambda}{k \mu}.
\]
We use $N$ to denote the random variable for the total number of jobs in the system (we will subscript $N$ by the load balancing policy when it is not clear from the context). Let $N/k$ give the (empirical) mean number of jobs per server, and let $\bar{\pi}(n)$ be the stationary probability that the mean number of jobs per server is more than $n$.
Since we will consider a sequence of systems parametrized by the number of servers, we will superscript these quantities as $\lambda^{(k)}, N^{(k)}, \bar{\pi}^{(k)}(n)$.
 
\subsection{Policies} \label{subsec:policies}

We now give a brief definition of the policies that we will address in the present paper:
\begin{itemize}
\item[\bf Join-the-Shortest-Queue (JSQ)] dispatches an arriving job to a server with the least number of jobs. Ties are broken randomly. The JSQ policy is the central object of study of this paper. 
\item[\bf Idle-Queue-First (IQF)] is a low overhead proxy for JSQ that prioritizes idle queues. If there are any idle servers then an arriving job is dispatched to an idle server, else a server is chosen at random. This policy was introduced in \cite{lu2011join}, where the authors call this Join-Idle-Queue.
\item[\bf Idle-One-First (I1F)] policy prioritizes idle queues and then queues of length one. That is, if there is an idle server then an arriving job is dispatched to an idle server 
else if there is a queue with one job 
then an arriving job is dispatched to a server whose queue has one job, 
otherwise a server is chosen at random. 
\item[\bf Power-of-$d$-Choices (Po$d$)] selects $d$ queues at random when a job arrives and the job is dispatched to the shortest of these $d$ queues with ties broken randomly.
\item[\bf Central Queue (CQ)] is an idealized policy where jobs are not immediately dispatched on arrival, but instead are kept in a single central buffer.
%
On arrival, jobs join a first-in first-out queue. From this queue jobs are routed to the next available server. 
In our setting, this policy results in an $M/M/k$ queueing model. 
\end{itemize}

Each of IQF, I1F and Po$d$ attempts to emulate the decisions of JSQ, but do not require knowing the state of all servers to make their decision. We study CQ as a comparison benchmark for the above policies.

\subsection{Non-Degenerate Slowdown Regime} \label{sec:motivation}
%
The \emph{Non-Degenerate Slowdown} regime (NDS) is a many server limit where the number of servers $k$ approaches infinity, and for a positive constant $\alpha$ the sequence of arrival rates $\lambda^{(k)}$ satisfy
\begin{equation}\label{motivation:NDS}
\lambda^{(k)} - k\mu \xrightarrow[k\rightarrow\infty]{} -\alpha \mu \quad \text{as}\quad k\rightarrow \infty.
\end{equation}
Here the service rate $\mu$ is fixed.

We now explain the nomenclature of the regime, and motivate it for the load balancing application we are studying. Recall that the \emph{slowdown} of a job is defined as the ratio between its sojourn time and its service requirement. 
The slowdown is at least $1$ by definition. An asymptotic regime has degenerate slowdown if the slowdown of a typical customer either converges to $1$ or diverges; otherwise the asymptotic regime exhibits non-degenerate slowdown. 
Throughout this paper, we consider the non-degenerate slowdown regime where
\begin{equation}\label{ourNDS}
\lambda^{(k)} = (k - \alpha)\mu.
\end{equation}

To explain why the limit regime \eqref{motivation:NDS} exhibits a non-degenerate slowdown, consider a Central Queue ($M/M/k$ queue) with service rate $\mu$ and arrival rate $\lambda^{(k)}= (k -\alpha) \mu $. From the Detailed Balance equations, we can calculate the stationary mean queue size $\bar{\pi}_{CQ}^{(k)}(n)$ from $\bar{\pi}_{CQ}^{(k)}(1)$ as follows
\begin{equation}\label{motivation:calc}
\bar{\pi}_{CQ}^{(k)}(n) = \left(\frac{\lambda^{(k)}}{k\mu}\right) \bar{\pi}_{CQ}^{(k)}\big(n-\tfrac{1}{k}\big)=...=\left(\frac{\lambda^{(k)}}{k\mu}\right)^{(n-1)k} \bar{\pi}_{CQ}^{(k)}\big(1\big).
\end{equation}
For the slowdown to be non-degenerate, we require the typical customer to visit a system where the mean number of jobs per-server, $n$, is such that $1<n<\infty$ (i.e. the typical customer does not see an idle server on arrival nor is there an infinite waiting time). For a limit distribution to have this support, $\bar{\pi}_{CQ}^{(k)}(n)$ must be some positive fraction smaller than $\bar{\pi}_{CQ}^{(k)}(1)$.
Given \eqref{motivation:calc}, this is equivalent to 
\[
\lim_{k\rightarrow \infty}\left(\frac{\lambda^{(k)}}{k\mu}\right)^{k} = e^{-\alpha}
\]
for some positive constant $\alpha>0$. This in turn implies that condition \eqref{motivation:NDS} must hold. So we see, for this example, that non-degenerate slowdown corresponds to the asymptotic \eqref{motivation:NDS}.

In addition, for stability to hold, the equilibrium fraction of non-idle servers must satisfy $\bar{\pi}_{CQ}^{(k)}(1)= \lambda^{(k)}/(k\mu)$. Thus for the Central Queue policy (CQ) we find that
\begin{equation}\label{CQDist}
\bar{\pi}_{CQ}^{(k)}(n) \rightarrow e^{-\alpha(n-1)} \quad\text{ as }\quad   k\rightarrow \infty.
\end{equation}
This result which is known to Atar, Gurvich and Whitt, is summarized in Section \ref{Comparison}.


\begin{remark}\label{Po2Remark}
	Although we will see that the asymptotic regime \eqref{motivation:NDS} corresponds to a non-degenerate slowdown limit for Central Queue and JSQ, less efficient policies can have a divergent slowdown in this asymptotic. For instance, the key observation on Power-of-two choices (Po2) is that, for $k$ large, \emph{each queue} has stationary distribution
	\[
	\bar{\pi}^{(k)}_{Po2}(n) = \Big(\frac{\lambda}{k\mu}\Big)^{2^{n-1}-1}.
	\]
	
	This stationary distribution is substantially more light-tailed than the geometric distribution achieved by randomized load-balancing. However, Po2 has infinite queue sizes under the NDS scaling \eqref{motivation:NDS}: for all $n$,
	\[
	\bar{\pi}^{(k)}_{Po2}(n)\xrightarrow[]{} 1,\quad\text{ as }\quad  \lambda - k\mu \xrightarrow[]{} -\alpha \mu.
	\]
	This suggests that Power-of-2 choices has substantially worse performance than that found for the Central Queue in Equation \eqref{CQDist}. To get non-degenerate slowdown for Power-of-2, one should set $\frac{\lambda^{(k)}}{k\mu} \to \rho < 1 $, which is the mean-field regime. 
\end{remark}

	Because of the observation remarked above, we do not -- and, indeed, can not -- analyze Power-of-two choices in the NDS regime. 
	We instead use the NDS regime to  discern between the efficient policies, such as CQ, JSQ, IQF, and I1F. 

\subsection{Comparison with other asymptotic regimes}\label{sec:OtherAsymp}

Parallel queueing systems are commonly analyzed in the following asymptotic regimes: Fluid Limits, Heavy Traffic Limits, Mean Field Limits and the Halfin-Whitt (Quality-Efficiency-Driven) Regime. 
We can view each of these regimes as the quantity $(\lambda - k \mu)$ taking on different limiting values after appropriately rescaling by some function of $k$. See the 2nd column of Table~\ref{Table:summary_regimes} for a list containing these scalings. 

For fluid and heavy traffic limits the number of servers is assumed fixed. Asymptotically, the number of jobs per-server diverges to infinity and, after appropriate rescaling, the limit for the total queue size process will have identical behaviour under policies such as JSQ and CQ. For $k=2$ and with servers of potentially different service rates, this fact was proved by \cite{foschini1978basic}. For general $k$ with homogeneous servers (equal service rates), this follows from the bounds obtained in \cite{adan1994upper}. 

Under a mean field limit, the number of servers grows to infinity but the load per-server approaches a limit bounded strictly away from 0 and 1.
The observation that a large under-loaded parallel server model will behave as an $M/M/\infty$ queue under JSQ is understood in \cite{badonnel2008dynamic}. 
 Therefore in mean field regime, there must always be a proportion of idle servers. Thus, this observation can then be applied to any policy that gives higher priority to idle servers (see \cite{lu2011join} for further discussion and  \cite{stolyar2015pull} for a rigorous analysis). So policies that prioritize idle servers such as CQ, JSQ, IQF and I1F always send jobs to an empty queue in a mean field limit and thus have an identical zero-one queue size distribution.  

A similar observation holds in the Halfin-Whitt regime. Here the number of servers, $k$, grows to infinity but the load per server is less that one by a $\Theta\left(\frac{1}{\sqrt{k}}\right)$ term (unlike $\Theta\left( \frac{1}{k}\right)$ for NDS). In this limiting regime, all but $\ord{\sqrt{k}}$ servers have queue size one, with the $\ord{\sqrt{k}}$ error characterized by a limiting diffusion. Thus mean queue sizes are $1$ in the Halfin-Whitt regime. See \cite{halfin1981heavy} for the classical analysis of CQ in this regime, \cite{eschenfeldt2015join} for a recent article analyzing the process level convergence to a limiting diffusion for JSQ, and \cite{braverman2018steady} which proves tightness of the JSQ system and hence establishes limiting steady-state distribution. 

Each scaling regime provides different insights into performance. However, we see that none of the standard limit regimes are able to discern even between policies such as CQ and JSQ. (See the last column of Table~\ref{Table:summary_regimes} for a summary). A contribution of this paper is to show that the NDS regime will differentiate between JSQ and CQ policies and thus provides a mechanism for the design and performance evaluation of different parallel server polices. In particular, we use our analysis to design a new load balancing policy, Idle-One-First (I1F) which is easier to implement than JSQ, has lower information overhead than Po$d$ and preserves asymptotic optimality of JSQ.



\begin{table}[] 
	\centering
	\begin{tabular}{l|c|c}
					& Asymptotic Framework	 & Equilibrium jobs per server   \\ \hline
 & &  \\				
 Fluid Limit	& \(\displaystyle \lambda -k\mu \xrightarrow[q \rightarrow \infty]{} -\alpha \mu \)$\qquad$  $\mu,k$ fixed & $\bar{\pi}$ is undefined  \\
 & &  \\ \hline
 & & \\				
	Heavy Traffic	& \(\displaystyle q (\lambda -k\mu) \xrightarrow[q \rightarrow \infty]{} -\alpha \mu \)$\qquad$  $\mu,k$ fixed & $\bar{\pi}$ is undefined   \\
 & &  \\
 \hline
  & &  \\			
	Mean Field		& \(\displaystyle\frac{\lambda- k\mu}{k}\xrightarrow[k \rightarrow \infty]{} -\alpha \mu \), $\qquad$  $\mu, \frac{q}{k}$ fixed 	&  $\bar{\pi}_{CQ} \equiv \bar{\pi}_{JSQ} \equiv \bar{\pi}_{I1F} \equiv \bar{\pi}_{IQF} $   \\ 
 & & \\			
	\hline
 & &  \\			
	Halfin-Whitt	& \(\displaystyle \frac{\lambda-k \mu}{\sqrt{k}} \xrightarrow[k \rightarrow \infty]{} -\alpha \mu \), $\qquad$$\mu, \frac{q - k}{\sqrt{k}}$ fixed   	& $\bar{\pi}_{CQ} \equiv \bar{\pi}_{JSQ} \equiv \bar{\pi}_{I1F} \equiv \bar{\pi}_{IQF} $   \\ 
 & &  \\		
	\hline
 & &  \\	
	Non-Degenerate Slowdown  & \(\displaystyle \lambda-k \mu \xrightarrow[k \rightarrow \infty]{} -\alpha \mu \),$\qquad$ $\mu, \frac{q}{k}$ fixed & $\bar{\pi}_{CQ} < \bar{\pi}_{JSQ} = \bar{\pi}_{I1F} < \bar{\pi}_{IQF} $ \\
	 & &  \\
\end{tabular}
	\caption{We detail limiting regimes given parameters: $\lambda$ arrival rate, $\mu$ service rate, $k$ number of servers, $q$ the initial number of jobs, $\bar{\pi}(n)$ stationary probability of greater than or equal to $n$ jobs per server. We indicate that, in contrast to NDS regime conventional asymptotic frameworks do not give meaningful stationary description of jobs per server. This holds for policies Central Queue (CQ), Joint-the-Shortest-Queue (JSQ), Idle-One-First (I1F)), Idle-Queue-First (IQF). }
	\label{Table:summary_regimes}
\end{table}


\section{Literature Review} \label{sec:literature}
We review policies used for parallel server models, paying particular attention to results on JSQ and its variants; then we review asymptotic methods used to analyze these systems; finally, we review recent work on the non-degenerate slowdown regime.

\subsubsection*{Parallel Server Models.}  

The earliest mathematical analysis' of Join-the-Shortest Queue are \cite{haight1958two, kingman1961two, flatto1977two}. In each case, a system with two servers is analyzed. 
It was subsequently proved by \cite{winston1977optimality} that, for exponentially distributed job sizes, JSQ minimizes the number of jobs among policies with past queue size and arrival information. See \cite{foss1980approximation,foss1984queues} for an elegant coupling proof of this result.

The Central Queue (CQ) policy goes back to Erlang and Kendall. When the arrival and service distribution is specified, the model is called the $M/M/k$ or $G/G/k$ queue. With FIFO scheduling rule at servers, Central Queue is equivalent to the policy where jobs are dispatched to the queue with least work, see \cite{foss1980approximation} and \cite[Chapter XII]{asmussen2008applied}. Thus, CQ is optimal among a broad class of dispatch policies. (In particular those having knowledge of each jobs work load, its distribution and past arrivals and assuming FIFO service at each queue). Again, see \cite[Chapter XII]{asmussen2008applied}.

Although CQ and JSQ are practical policies in some applications, such as call centers where the number of servers is moderate, in other applications, such as data centers, dispatch decision are made on the arrival of each job and with limited state information. For this reason, there has been considerable interest in finding simple policies that emulate the optimal behavior of JSQ and CQ.

\cite{mitzenmacher2001power} and \cite{vvedenskaya1996queueing} analyzed the Power-of-$d$ choices policies. One motivation for Power-of-$d$ rule is as a policy that tries to approximately emulate the performance of JSQ by randomly sampling a subset of servers. A second motivation comes from data center applications where compute jobs must be co-located with the data. By replicating each data on $d$ randomly chosen servers, the load balancing of compute tasks behaves as if the the $d$ servers were randomly sampled on the arrival of the task.   
 An alternative mechanism to emulate JSQ is to prioritize idle servers. This approach is emphasized in \cite{Laws00}. In the parallel server models, dispatching tasks as is called the Idle-Queue-First Policy, see \cite{lu2011join}. As discussed in \cite{badonnel2008dynamic},  IQF emulates the behavior of an $M/M/\infty$ queue in a mean-field asymptotic and thus can be argued to outperform Po$d$. See \cite{stolyar2015pull} for a mathematically rigorous analysis and discussion, as well as results with more general job sizes in \cite{foss2017large}. Accordingly the analysis of these policies and variants have attracted growing research interest with articles such as \cite{tsitsiklis2012power}, \cite{mukherjee2015universality}, \cite{gamarnik2016delay} and \cite{jonckheere2016asymptotics} each analyzing different trade-offs present in this modeling framework.

As indicated above, although analysis of classical policies such as CQ and JSQ goes back more than fifty years, the development of more modern data center applications have necessitated new policy variants and new methods to compare and analyze them.

\subsubsection*{Methods for Asymptotic Analysis.}
We now review scaling regimes applied to CQ and JSQ. 
 Heavy traffic limits were first proposed by Kingman in order to provide simplified descriptions of queueing system that extract the key characteristics and differences in different queueing processes (a point which motivates for our NDS analysis). 
\cite{IgWh70} first characterized the limiting diffusion of the Central Queue policy under a heavy traffic scaling.
A heavy traffic analysis of the Join-the-Shortest-Queue for 2 servers is provided by \cite{foschini1978basic}. Following this work \cite{halfin1981heavy} developed a scaling regime for the central queue policy where the number of servers was allowed to grow with the load parameters (often called the Halfin-Whitt, or Quality-Efficiency-Driven regime). \cite{borst2004dimensioning} motivate, reformulate and refine the Halfin-Whitt regime as the cost minimizing regime for call center dimensioning. More recently \cite{eschenfeldt2015join} provide an analysis of JSQ in the Halfin-Whitt regime, and \cite{braverman2018steady} uses Stein's method to prove the tightness of JSQ thus establishing validity of limit interchange needed to extend the result of \cite{eschenfeldt2015join} to steady-state distribution under Halfin-Whitt regime. In addition to further literature review, \cite{pang2007martingale} provides a good review of methods used to characterize the limiting diffusions in the Halfin-Whitt regime. 

Mean field limits have commonly been applied to the modeling of large queueing systems, particularly in settings where queues behave approximately independently (as in the case for the Power-of-$d$ load balancing systems -- see \cite{Bramson2012}). A mean field analysis of Po$d$ was given by \cite{vvedenskaya1996queueing} and further by \cite{bramson2013decay}. Introducing dependence between queues can substantially improve performance. A mean field analysis of this effect is given by \cite{tsitsiklis2012power} and for the IQF policies by \cite{stolyar2015pull}. \cite{ying2016rate} employs Stein's method to establish the rate of convergence of the steady-state measure to the mean-field equlibrium for Power-of-2 system.

\subsubsection*{Time-Scale Separation.}
{
From a mathematical perspective, there are a rich set of phenomena that occur in our analysis. For instance, we find that there is separation of time-scales between several important components of our process. This is somewhat analogous to the snapshot principle, as applied in heavy traffic scaling. 
Separation of timescale in Markov processes is first analyzed by \cite{has1966stochastic}. This approach is generalized by \cite{Veretennikov91}. A broad text book treatment is given in \cite{freidlin1998random}.
The analysis of \cite{hunt1994large} provides a Markov generator approach to the analysis of queueing systems. A more recent paper \cite{luczak2013} provides a time-scale separation analysis for Power-of-$d$ policies. There are difficulties in directly applying these work in our context, particularly, since our diffusion limit has non-Lipschitz coefficients.
Our analysis benefited much from the methods of \cite{robert2015stochastic}. Coupling approaches of \cite{walton2012flow} are also applied.  
Understanding the separation of time-scales of the idleness process leads us to understand optimal queue size behavior, and to the development of I1F policy.
}

\subsubsection*{Non-Degenerate Slowdown.}
The Non-Degenerate Slowdown Regime was first introduced by \cite{atar2012diffusion}. The NDS regime has system size fluctuations which are larger than those considered in the Halfin-Whitt regime. Atar motivates the use of NDS in call centers from empirical work, where customers' slowdown can be large, see \cite{brown2005statistical}.  
Atar attributes first use of this scaling to \cite{whitt2003multiserver} and \cite{gurvich2004design} (and to private communications with A. Mandelbaum and G. Shaikhet), and 
Subsequent work in the non-degenerate slowdown regime developed notions of optimal control, see \cite{atar2011asymptotically} and \cite{atar2014scheduling}.
Independently of Atar, the work of \cite{maglaras2013optimal} consider the non-degenerate slowdown regime. Here the scaling regime is motivated by a cost minimization objective. {A recent work by \cite{He15}  views prior work on NDS as a form of space-time scaling and provides diffusion approximations with general service times for call-centre models.}

The NDS regime is a relatively new scaling asymptotic with most key papers appearing in the last 6 years and originating from motivations on call center applications, where customer are served in order of arrival. Our work advocates the use of NDS in order to distinguish between load balancing policies and is primarily motivated by recent data center applications.

\section{Join-the-Shortest-Queue in Non-Degenerate Slowdown Regime} \label{sec:NDS1}
In this section, we formally state the limiting behavior of the parallel server system with JSQ dispatch rule as number of servers $k\to \infty$ under NDS regime (Theorem~\ref{JSQ:Thrm}). While a rigorous proof of Theorem~\ref{JSQ:Thrm} is quite intricate and is presented in the appendix, the intuition behind the result is quite simple and we present a heuristic derivation of the result in Section~\ref{sec:Heuristic}.

Let $N^{(k)}(t)$ be the total number of jobs at time $t$ in the $k$ server system in the NDS scaling $\lambda^{(k)} =  k\mu - \alpha\mu$. 
We assume that, at time $t=0$, queues are balanced, that is all queue sizes are within one of each other. 
(As discussed shortly, this state of balance is naturally reached by the JSQ policy.) In an NDS limit we re-scale the process as follows
\begin{equation}\label{NNDS}
\hat{N}^{(k)}(t) : = \frac{{N^{(k)}}(kt)}{k}. 
\end{equation}

Note that $\hat{N}^{(k)}(t)$ represent the mean number of jobs per-server. 
We assume that $\hat{N}^{(k)}(0)$ converges in distribution to a random variable $\hat{N}(0)$ where ${\mathbb P} (\hat{N}(0) \leq 1)=0$. (Again, we see such states such that $\hat{N}(t) \leq 1$ are not attained in the NDS limit.) Here and throughout the rest of the paper we let $\Rightarrow$ indicate weak convergence with respect to the uniform topology on compact time intervals, see \cite[Chapter 2 and Section 15]{billingsley2013convergence} . 

The following is the main result proven in this paper.

\begin{theorem}[Join-the-Shortest Queue in NDS]\label{JSQ:Thrm}
	For  a parallel server model operating under Join-the-Shortest-Queue as stated above,
	 the diffusion-scaled process for number of jobs in the system, converges in distribution:
	\begin{equation}
	\hat{N}^{(k)}\Rightarrow \hat{N}
	\end{equation} 
	uniformly on compact time intervals where $\hat{N}$ is the path-wise unique solution to the stochastic differential equation
	\begin{equation}\label{JSQ:Diffusion}
	d \hat{N} (t) = \mu  \Big[ \frac{(2-\hat{N})_+}{\hat{N}-1}-\alpha\Big] dt + \sqrt{2\mu} \;dB (t), \qquad t\geq 0
	\end{equation}
	where $B=(B(t):t\geq 0)$ is a standard Brownian motion. Moreover, $\hat{N}$ has stationary distribution $\pi_{JSQ}$ with density 
	\begin{equation}\label{JSQ:Equilibrium}
	\frac{d\pi_{JSQ}}{dn} = 
	\begin{cases}
	C \cdot e^{ -\alpha (n-2) }, & n \geq 2, \\
	C \cdot (n-1)e^{ -(\alpha+1)(n-2) }, & 1\leq n\leq 2;
	\end{cases}
	\end{equation}
	 mean
	\begin{align}
	\label{JSQ:EN_Equilibrium}
	{\mathbb E}[\hat{N}] &= 1 + C\left[  \frac{\alpha^2+4\alpha+1}{\alpha^2 (1+\alpha)^3} + \frac{2 e^{\alpha+1}}{(1+\alpha)^3} \right],
	\end{align}
	and normalizing constant 
\[
C = \left[ \frac{1}{\alpha (1+\alpha)^2} + \frac{e^{1+\alpha}}{(1+\alpha)^2} \right]^{-1}.
\]
\end{theorem}

\begin{figure}[t]
	\begin{center}
		\includegraphics[width=0.6\textwidth]{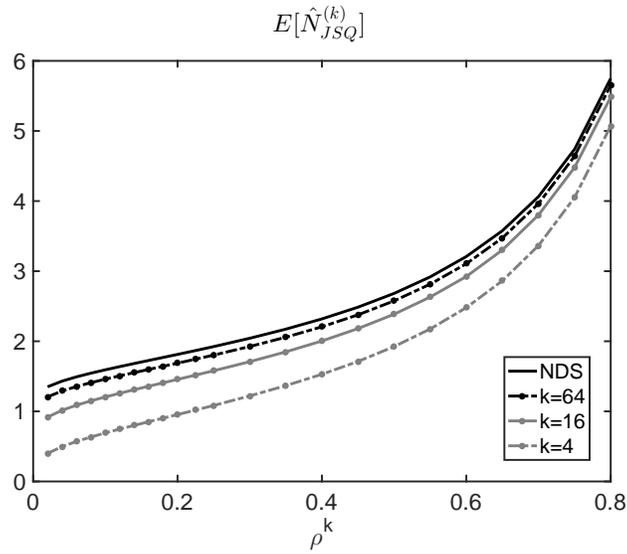}
	\end{center}
	\caption{Simulation results comparing  ${\mathbb E}{\hat{N}^{(k)}_{JSQ}}$ for $k=4,16,64$ servers with the limiting expression under an NDS regime. (The values on the $x$-axis are expressed in terms of $\rho^k$ since $\rho^k \to e^{-\alpha}$ in the NDS regime.)}
\label{fig:Sim_EN_JSQ}
\end{figure}

Thus we see that, like the Central Queue policy, Join-the-Shortest Queue has a non-trivial limit in the NDS regime (and as we show later, distinct from that of Central Queue). Simulation results demonstrating convergence to the NDS limit are presented in Figure \ref{fig:Sim_EN_JSQ}. To understand the implications of this result, we must compare the diffusion process found above with the diffusion processes corresponding to the other dispatch rules described in Section~\ref{subsec:policies}. This is done in Section \ref{sec:NDS2}. We note that, by Little's law, the expected response time divided by the expected job size is precisely ${\mathbb E}[\hat{N}]$, as given above. In this sense we can characterize the slowdown asymptotic in NDS.

To understand why the above NDS diffusion arises, we need to analyze the interplay between the total queue size, the service at individual queues and the quantity of system idleness, all of which interact on different timescales. We discuss this in the following subsection.

\subsection{Heuristic Derivation of Theorem \ref{JSQ:Thrm}}\label{sec:Heuristic}
To better position Theorem \ref{JSQ:Thrm}, we present a heuristic derivation of the result. The formal argument is substantially more intricate and is presented in the Appendix.  

Consider the $k$ parallel server model under JSQ, as described in Theorem \ref{JSQ:Thrm}. Further let $I^{(k)}$ be the total number of idle servers in this queueing system.
The total number of jobs $N^{(k)}$ and the number of idle servers $I^{(k)}$ are the main quantities of interest. 

First observe that, given the number of idle servers $I^{(k)}(t)$, the transitions in $N^{(k)}$ at time $t$ occur at the following transition rates
\begin{align*}
N^{(k)} &\mapsto  N^{(k)}+1,\qquad \text{at rate} \qquad (k -\alpha) \mu ,\\
N^{(k)} &\mapsto  N^{(k)}-1, \qquad \text{at rate } \qquad (k-I^{(k)}(t)) \mu.
\end{align*}
So to understand the behaviour of $N^{(k)}$, we must also analyze $I^{(k)}$.

Second, we focus on a given value of $N^{(k)}$ and ask what the typical joint state of the $k$ servers looks like under JSQ. Note that, under JSQ, servers that have an above average number of jobs per server will receive no arrivals, and thus will decrease until less than or equal to the average. This suggests that queues balance in a way that the difference between the shortest queue and the longest is close to, and frequently must be equal to, one. 

\begin{figure}[t]%
	\centering
	\includegraphics[width=0.8\textwidth]{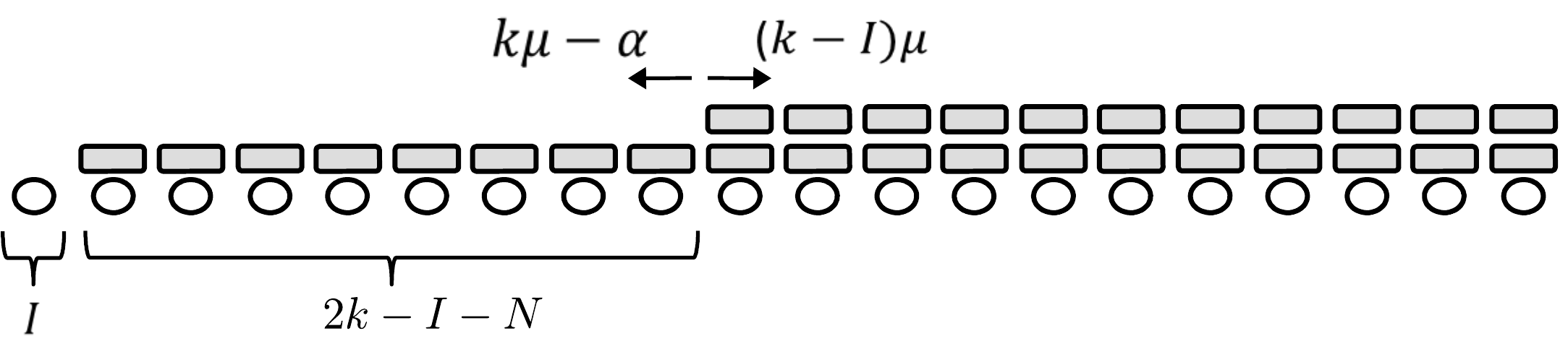}
	\caption{Representation of JSQ with $k$ servers and less than $2$ jobs per-server. 
		\label{TimeScalePic}}
\end{figure}%

Next let us consider the fluctuations in this difference and what this means for the number of idle servers. As an example, let $N^{(k)}=3k/2$ where each server has $3/2$ jobs on average. See Figure~\ref{TimeScalePic} for reference.
Approximately, half the queues are of length 1, the other half of length 2 and the number of idle queues is $\ord{1}$. 
The critical observation is that 
only departures from queues of length one can create an idle server and this idle server can not be made busy until a subsequent external arrival occurs. (This is in contrast to call center applications, where we could have moved a job from another queue in order to maintain busy servers). More concretely, idle servers evolve according to transition rates:
\begin{subequations}\label{ITrans}
\begin{align}
I^{(k)} &\mapsto  I^{(k)}+1,\qquad\;\; \text{at rate} \qquad\; k\mu \cdot \hat{M}^{(k)}_1  ,\\
I^{(k)} &\mapsto  (I^{(k)}-1)_+, \quad \text{at rate } \qquad k\mu-\alpha \mu,
\end{align}	
\end{subequations}
where $\hat{M}^{(k)}_1(t)$ denotes the proportion of queues with one job, which in our example is approximately $1/2$. 
Note that the transitions in $I^{(k)}$ occur on a timescale of order $\ord{k^{-1}}$. Further, 
given $\hat{M}^{(k)}_1$ is \emph{approximately} constant, the idle server process is \emph{approximately} an $M/M/1$ queue with expected value:
\begin{align}
{\mathbb E} [I^{(k)}|N^{(k)}] 
&= 
\frac{\hat{M}_1^{(k)}}{1-\hat{M}_1^{(k)}} + o(1)\notag \\
&=
\frac{
	(2-\hat{N}^{(k)})_+
}{
	\hat{N}^{(k)} - 1
}
+ 
o(1) \label{Iprox}
\end{align}
where we observe that, given queues are balanced, we have $\hat{M}_1^{(k)}=(2-\hat{N}^{(k)})_+$.
In our example, this yields $\expct{I^{(k)}|N^{(k)}=3k/2} = 1 + o(1)$, i.e. when average queue size is 3/2, on average we lose the capacity of one server due to JSQ dispatch.

Of course $\hat{M}^{(k)}_1$ is not constant. However, note that the $M/M/1$ defined by \eqref{ITrans} makes transitions at rate $\ord{k}$ and thus has an $\ord{k^{-1}}$ mixing time; while, on the other hand, ${N}^{(k)}$ requires an $\ord{k}$ time to make an order $1$ change in the value of  $\hat{M}^{(k)}$. (To see this note that arrival and departure rates are $\ord{k}$ but their difference between arrival and departure rates of ${N}^{(k)}$ are $\ord{1}$.) In other words there is a timescale separation between a fast process, $I^{(k)}$, and a slow processes, ${N}^{(k)}$. 

 We prove this time-scale separation, yielding a \emph{stochastic averaging principle} where the idleness process $I^{(k)}$ evolves on a fast time scale with its stationary distribution a function of $N^{(k)}$ and where the stationary value of $I^{(k)}$ drives the slow timescale evolution of $N^{(k)}$. This makes rigorous the following  approximation for the evolution of $N$
\begin{align*}
N^{(k)} &\mapsto  N^{(k)}+1,\qquad \text{at rate} \qquad (k -\alpha) \mu ,\\
N^{(k)} &\mapsto  N^{(k)}-1, \qquad \text{at rate } \qquad (k-{\mathbb E}[I^{(k)} | N^{(k)}(t)]) \mu.
\end{align*}
Under the scaling $\hat{N}^{(k)}(t)=N^{(k)}(kt)/k$ and given \eqref{Iprox}, this leads to the diffusion approximation:
\begin{align*}
	d \hat{N} (t) &= \mu  \Big[ {\mathbb E}[{I^{(k)} | N^{(k)}}]-\alpha\Big] dt + \sqrt{2\mu} \;dB (t)\\
& =\mu  \Big[ \frac{(2-\hat{N})_+}{\hat{N}-1}-\alpha\Big] dt + \sqrt{2\mu} \;dB (t)
\end{align*}
which is the form made precise in Theorem \ref{JSQ:Thrm} and its proof.

\section{Informal discussion using NDS to compare JSQ with other policies} \label{sec:NDS2}

As we mentioned earlier, Join-the Shortest-Queue is a classical and widely studied dispatch rule due to its simplicity and optimality properties. 
In this section, we want to use our NDS analysis to better understand how JSQ compares with other policies which have greater or lesser control over a parallel queueing system. 
Specifically, we use our NDS analysis to answer the following questions:
\begin{itemize}
	\item [\textbf{Q1.}] What is the price of immediately dispatching jobs, rather than storing jobs in a central buffer or jockeying jobs from another queue  when a server becomes idle? 
\item [\textbf{Q2.}] Given that JSQ requires all queue size information, what is the impact of using more limited queue size information (such as just idleness information or just a finite subsample of queue sizes)?
 
\item [\textbf{Q3.}] Can we use our NDS analysis of JSQ to improve upon previous well studied policies such as IQF (and Po$d$)?
\end{itemize}

 Below we specify the NDS approximation and stationary distribution for the Central Queue, Idle-Queue-First and Idle-One-First policies.  We, then, compare the performance of these policies with Join-the-Shortest-Queue providing simulations to support our theoretical findings. 

\subsection{CQ, IQF and I1F in NDS}\label{Comparison}
The NDS approximations for Central Queue, Idle-Queue-First and Idle-One-First are as follows.

\textbf{CQ in NDS:} The NDS limit of the Central Queue policy (or $M/M/k$ queue) is a reflected Brownian motion with drift and a reflection at one. Here is $\hat{N}_{CQ}$ is the stochastic process with range $[1,\infty)$ satisfying the stochastic differential equation
\begin{equation*}
 d\hat{N}_{CQ}(t) = - \alpha \mu dt + \sqrt{2\mu}d B(t) + dL(t),\qquad t\geq 0,
\end{equation*}
where the increasing processes $L=(L_t: t\geq 0)$ is the reflection term at $1$. 
This diffusion has stationary distribution $({\pi}_{CQ}(n): n \geq 1)$ with density:
\begin{equation*}
	\frac{d\pi_{CQ}}{dn}
	=
	\alpha 
	\exp \{ -\alpha(n-1)\},
	\qquad n\geq 1.
\end{equation*}
Here the stationary expected queue size is
\begin{equation}\label{CQED}
\mathbb{E}[\hat{N}_{CQ}] = 1+ \frac{1}{\alpha} 
\end{equation}

This NDS approximation is proven by \cite{atar2012diffusion}. 
This form can be anticipated: when all servers are busy the change in queue size is the difference between two Poisson processes. This gives the Brownian motion with drift in the NDS limit. Further since the difference between the arrival and departure rates is $\alpha$ under an NDS scaling, the number of idle servers is $\ord{1}$ for any $k$. This provides the reflection at $\hat{N}_{CQ}=1$.


\textbf{IQF in NDS:} The Idle-Queue-First policy has an NDS approximation given by
\begin{equation*}
	d\hat{N}_{IQF} (t)
= 
\mu \left[
		\frac{1}{\hat{N}_{IQF}-1}
		-
		\alpha 
	\right]  dt
+
\sqrt{2\mu}dB(t),\qquad t\geq 0.
\end{equation*}
Note that $\hat{N}_{IQF} $ is a Bessel process with additional drift $-\alpha \mu t$. 
This process has stationary distribution
$({\pi}_{IQF}(n): n \geq 1)$ with density:
\begin{equation*}
	\frac{d\pi_{IQF}}{dn}
	\propto 	
\begin{cases}
\exp \{-\alpha(n-2) \}	 &	n\geq 2, \\
(n-1) \exp \{ -(\alpha +1)(n-2)\}	& 1\leq n \leq 2.
\end{cases}
\end{equation*}
Here the stationary expected queue size is
\begin{equation*}
	\mathbb{E} [ \hat{N}_{IQF} ] 
	=
	1+ \frac{2}{\alpha}.
\end{equation*}

\begin{figure}[t]
	\begin{center}
		\subfigure[Idle-Queue-First]{
		\includegraphics[width=0.47\textwidth]{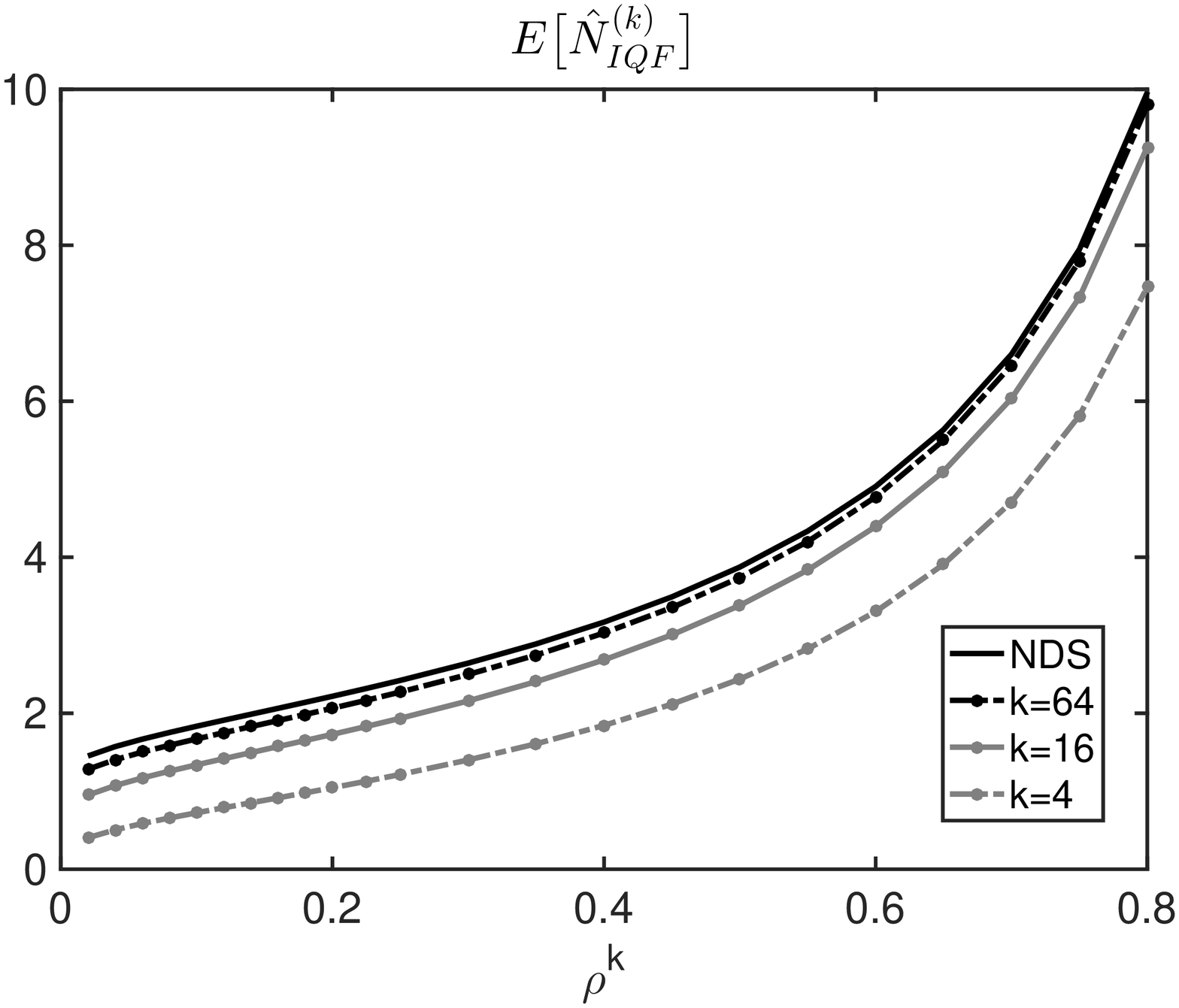}
		\label{fig:Sim_EN_IQF}
		}
		\subfigure[Idle-One-First]{
		\includegraphics[width=0.47\textwidth]{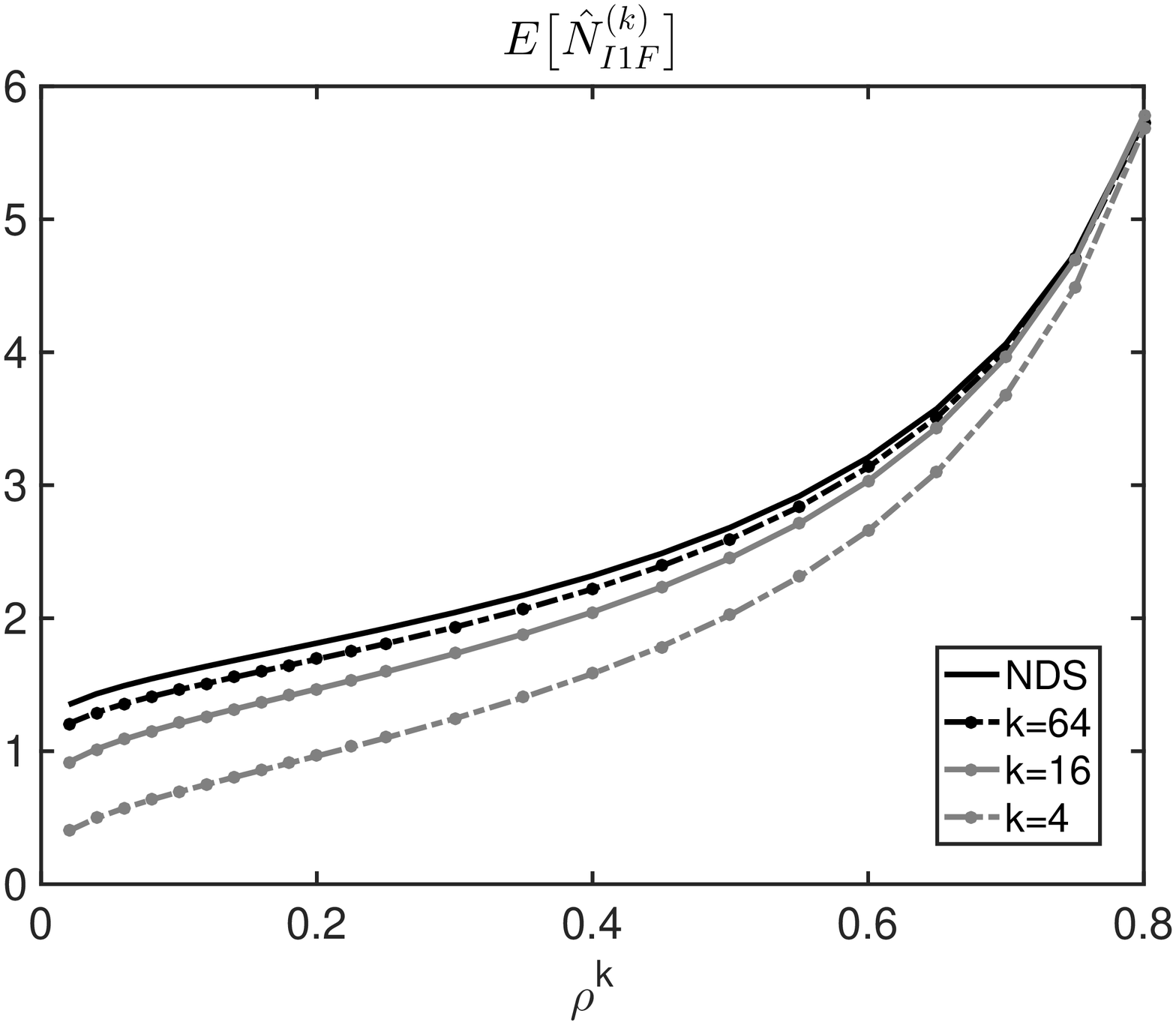}
		\label{fig:Sim_EN_I1F}
		}
	\end{center}
	\caption{Simulation results comparing  $\expct{\hat{N}^{(k)}_{IQF}}$ and $\expct{\hat{N}^{(k)}_{I1F}}$ for $k=4,16,64$ servers and the NDS diffusion approximation.}
\label{fig:Sim_EN_IQF_I1F}
\end{figure}

We do not provide a formal proof of this NDS approximation. Simulations in Figure \ref{fig:Sim_EN_IQF} help justify this diffusion approximation.
 An argument along similar lines to our JSQ proof is possible. 
Similar to the discussion in Section \ref{sec:Heuristic}, the intuition for this limit is as follows. When the mean number of jobs per server is $\hat{N}^{(k)}_{IQF}$, the number of idle servers is $\ord{1}$ and is governed by the dynamics of an $M/M/1$ queue with departure rate $\lambda^{(k)}$ and arrival rate given by the number of queues of length $1$. 
Again, there is a time-scale separation between idle-server process and the proportion of queues of length $1$. However, for IQF, jobs are now routed to a random server when all servers are non-idle and, due to time-scale separation, this randomized routing occurs as a Poisson process for each server. The behaviour of each non-idle server is effectively that of a single server queue -- which is geometrically distributed. The total number of jobs in non-idle servers is $\hat{N}^{(k)}_{IQF}$. This specifies the mean of this geometric distribution to be $\hat{N}^{(k)}_{IQF}$ and, consequently, the fraction of queues of length one is $ \hat{M}^{(k)}_{1,IQF} \approx \frac{1}{\hat{N}_{IQF}^{(k)}}$, giving
$$ \mathbb{E}[I^{(k)}_{IQF} | \hat{N}^{(k)}_{IQF} ] = \frac{1}{\hat{N}^{(k)}_{IQF} -1 } + o(1). $$ 

This account for the loss of capacity in this policy which is expressed in the term in square brackets for our $IQF$ approximation. 

\textbf{I1F in NDS:} 
The Idle-One-First policy has an NDS approximation that is identical to the NDS approximation found for JSQ. That is
$\hat{N}_{I1F}$ has stationary distribution $\pi_{JSQ}$ with density 
	\begin{equation*}
	\frac{d\pi_{I1F}}{dn} \propto  
	\begin{cases}
	 e^{ -\alpha (n-2) }, & n \geq 2, \\
	 (n-1)e^{ -(\alpha+1)(n-2) }, & 1\leq n\leq 2;
	\end{cases}
	\end{equation*}
	and mean
	\begin{align*}
	{\mathbb E}[\hat{N}_{I1F}] &= 1 + \left[ \frac{1}{\alpha (1+\alpha)^2} + \frac{e^{1+\alpha}}{(1+\alpha)^2} \right]^{-1}\left[  \frac{\alpha^2+4\alpha+1}{\alpha^2 (1+\alpha)^3} + \frac{2 e^{\alpha+1}}{(1+\alpha)^3} \right]
	.\end{align*}
Simulations in Figure \ref{fig:Sim_EN_I1F} show that the stationary distribution under I1F converges to the NDS approximation for JSQ. The intuition for this is as follows. From the analysis of JSQ and IQF, we see that, to reduce queue lengths in the NDS regime, we must minimize the number of queues of length $1$. This is how JSQ retains its optimality properties in the NDS regime.
When $1< \hat{N}^{(k)}_{I1F} < 2$, the behaviour of I1F is identical to JSQ and the heuristic followed in Section \ref{sec:Heuristic} is identical for this policy.
When $\hat{N}^{(k)}_{I1F} \geq 2$, (similar to the argument for Idle-Queue-First) all but $\ord{1}$ queues must have length less than $2$, and consequently, $\ord{1/k}$ have length $1$. Thus there is no loss of capacity when $\hat{N}^{(k)}_{I1F} \geq 2$.
So, in summary, I1F behaves identically to JSQ for  $1< \hat{N}^{(k)}_{I1F} < 2$ and, like JSQ, has no idle servers for $\hat{N}^{(k)}_{I1F} > 2$ in the NDS regime. This yields an identical NDS limit to JSQ, while retaining substantially less communication overhead.

\begin{remark}[Po$d$ in NDS]
	We do not include the Power-of-$d$-choices policies in our discussions above, because as discussed in Remark \ref{Po2Remark}, the NDS limit of the Power-of-d-choices gives an infinite queue length, which again suggests that inefficiencies exist for this class of policies. 
\end{remark} 

In summary, we see that along the lines discussed in Section~\ref{sec:Heuristic}, it is possible to give diffusion approximations for the NDS regime for various policies and these are validated by simulation. Moreover, the rationale behind these diffusion approximations can be used to construct new policies such as I1F. The systematic study of the NDS regime for these policies, along the lines presented for JSQ in the Appendix, remains an area of active research interest. 

\subsection{NDS Policy Comparison}\label{Comparison}
In this section, we compare the NDS diffusions presented above with the diffusion found for JSQ. We use these approximations to answer the three questions posed at the beginning of Section~\ref{sec:NDS2}. 

First, we observe stationary distributions found satisfy the following inequalities
\[
	\pi_{CQ} \leq_{st} \pi_{JSQ} =_{st} \pi_{I1F} \leq_{st} \pi_{IQF}.
\]
Here $\pi_X \leq_{st} \pi_Y$ when $\mathbb{P}(X\leq x) \geq \mathbb{P}(Y \leq x)$.
The proof of the above follows because the drift term in square brackets for CQ, JSQ, IQF and IQF each dominate each other. (We provide a proof in Lemma \ref{St_Dom} in the Appendix.) This results in stationary distributions which each dominate the other. We now investigate the magnitude of these differences.

\begin{figure}[t]
\begin{center}
		\includegraphics[width=0.45\textwidth]{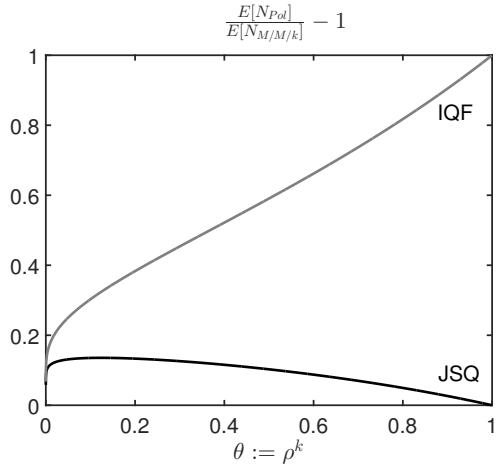}
\end{center}
\caption{Relative performance loss of JSQ and IQF compared to Central Queue  in the NDS asymptotics as a function of $\theta := e^{-\alpha}$.\label{fig:NDS_JSQ_IQF}}
\label{fig:sim_policy_comparison}
\end{figure}

\textbf{A1.} Let us answer \textbf{Q1}, as given at the start of Section \ref{sec:NDS2}. 
 Under JSQ, we must dispatch a job irrevocably to a server. Without this restriction, the optimal blind policy with exponential service distribution is the \emph{Central Queue} policy where jobs queue up at a central buffer, and whenever a server completes service, it picks the next job from this central buffer.
From our NDS analysis, we find that the price of immediate dispatch is at most $14\%$:
\begin{equation}\label{JSQ/CQ}
 \sup_{\alpha > 0} \frac{{\mathbb E}[\hat{N}_{JSQ}]}{{\mathbb E}[\hat{N}_{CQ}]} \leq 1.14 .
\end{equation}
See Figure \ref{fig:NDS_JSQ_IQF}. Thus, 
the impact of dispatching jobs (compared to maintaining a central pooled resource) does not have as severe an impact on total queue size as one might first expect. 

\textbf{A2.} Let us answer \textbf{Q2}, what is the price of ignoring queue length information?
The {Idle-Queue-First} is a cheaper alternative to JSQ which prioritizes dispatching to idle servers over busy servers, but ignores the queue lengths of the busy servers. For the given diffusion approximations for JSQ and IQF, we see that the mean sojourn time of IQF in NDS regime can be $100\%$ larger than Central Queue:  
\begin{equation}\label{IQF/CQ}
\sup_{\alpha>0} \frac{\mathbb{E}\big[{\hat{N}_{IQF}}\big]}{\mathbb{E}{\big[\hat{N}_{CQ}\big]}} = 2 .
\end{equation} 
See Figure \ref{fig:NDS_JSQ_IQF}. Again, IQF is within a constant factor of CQ. 

\begin{figure}[t]
	\begin{center}
		\subfigure[Idle-Queue-First]{
		\includegraphics[width=0.47\textwidth]{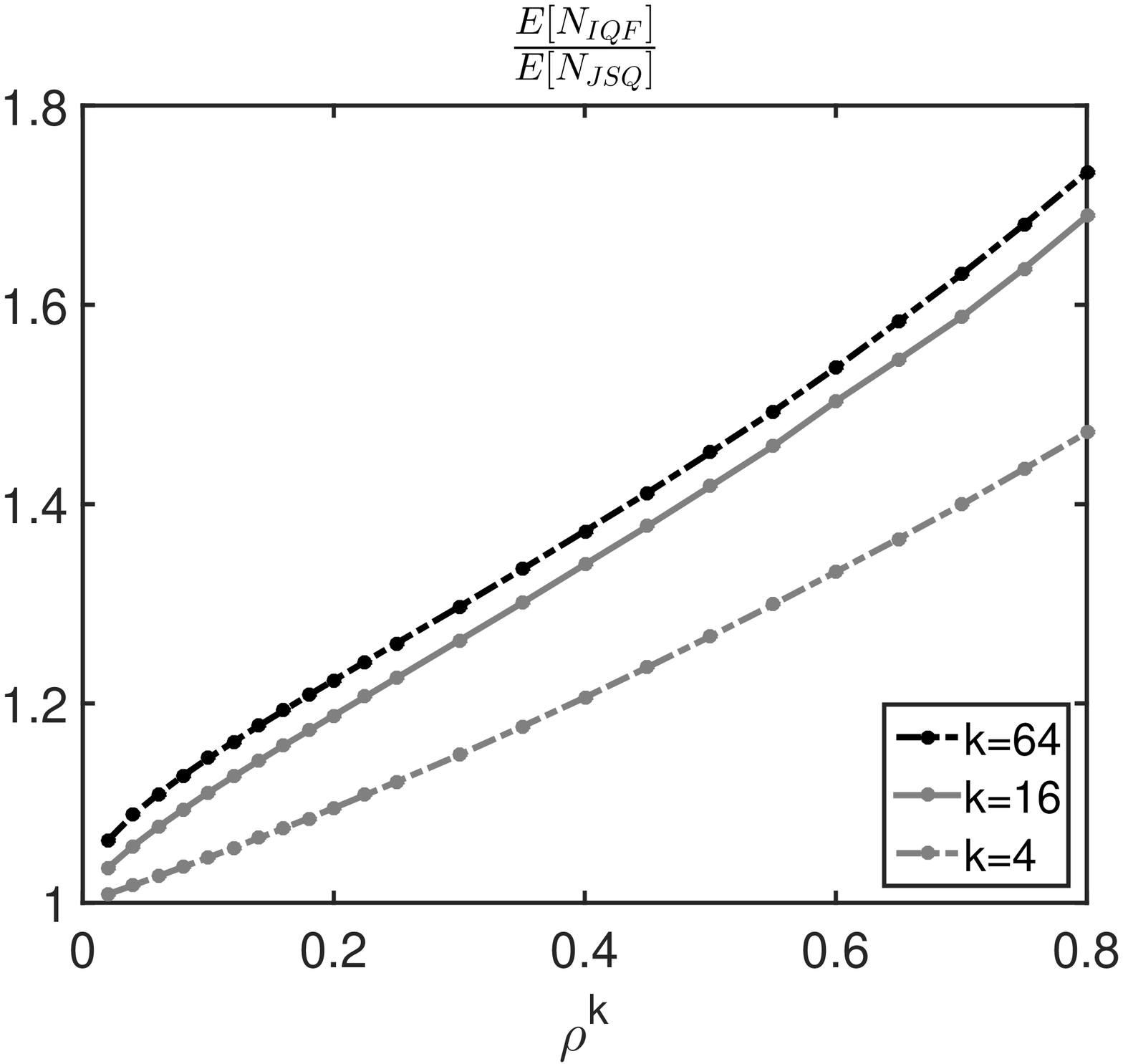}
		}
		\subfigure[Idle-One-First]{
		\includegraphics[width=0.47\textwidth]{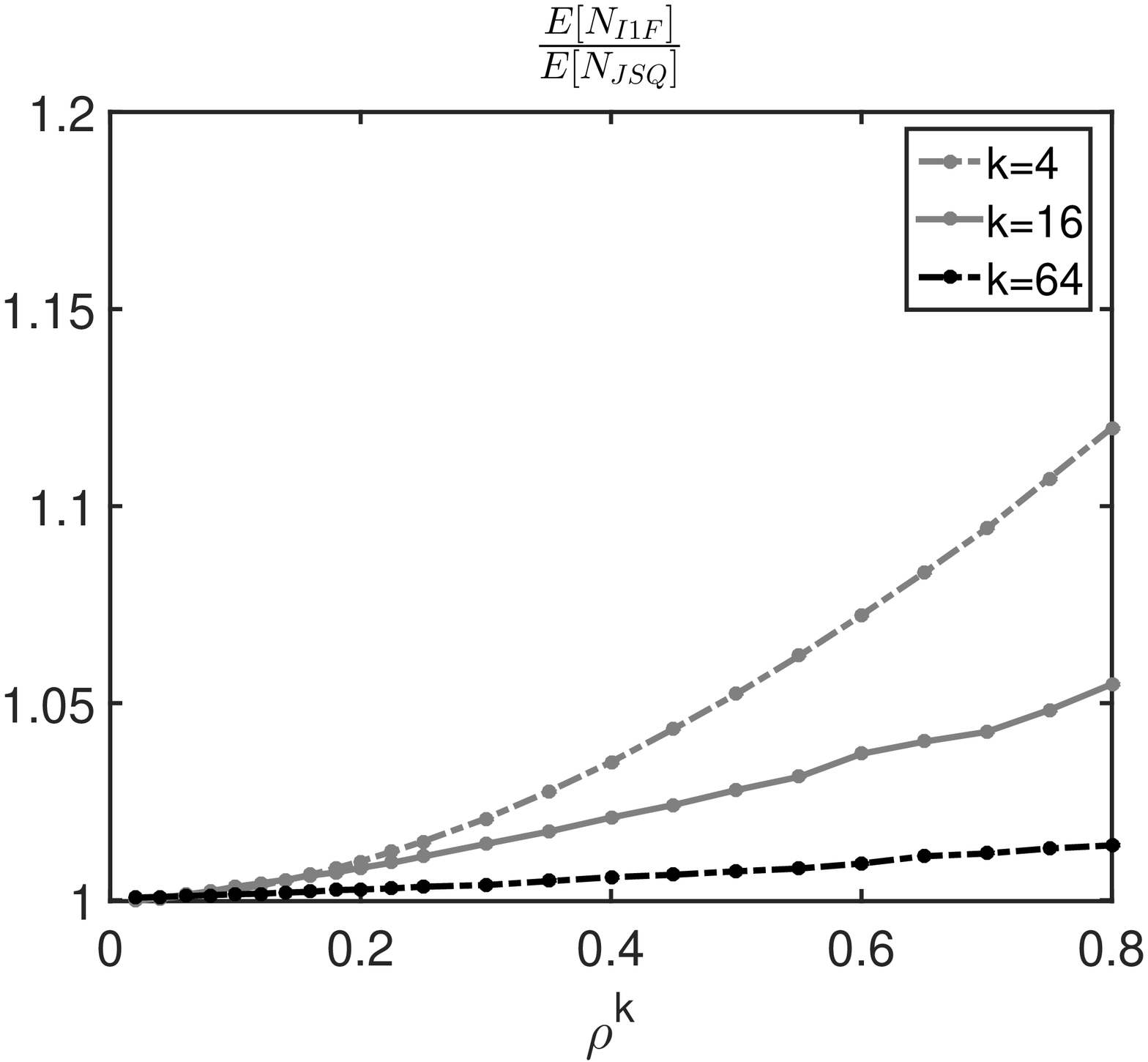}
		}
	\end{center}
	\caption{Simulation results comparing  $\expct{N^{(k)}_{IQF}}/\expct{N^{(k)}_{JSQ}}$ and $\expct{N^{(k)}_{I1F}}/\expct{N^{(k)}_{JSQ}}$ for $k=4,16,64$ servers and various values of $\rho$.}
\label{fig:Sim_EN_IQF_I1F_vs_Sim_JSQ}
\end{figure}

\textbf{A3.}  
As we saw, the Idle-One-First policy has identical performance to JSQ in the NDS regime.
Thus it is possible to achieve the same asymptotic performance as JSQ with significantly less queue size information.
A priori when introducing the Idle-One-First policy, one might expect a family of policies where we use JSQ dispatch rule among queues of length at most $q$, and use random routing if all queues are longer than $q$. 
However, through our analysis of JSQ, we see that little improvement will be achieve through this. 
Surprisingly, while there is benefit in looking at queues of length 1, there is asymptotically no benefit of going beyond this for the metric of mean response time (higher moments of response time would decrease as we move from I1F to JSQ).
This is because the NDS approximations for I$q$F are identical to JSQ for $q$ large than $1$. 
(This is analogous to the mean field analysis of power-of-$d$-choice, where there is an exponential factor improvement for power of $d=2$ choices but only constant factor improvement thereafter.) See Figure \ref{fig:Sim_EN_IQF_I1F_vs_Sim_JSQ} for a comparison between IQF to JSQ and I1F to JSQ for different numbers of servers.

\strut\vspace{0.0cm}

We remark that these insights comparing IQF, I1F and JSQ are only apparent in the NDS regime. 
In the Halfin-Whitt (QED) regime the number of queues with 2 or more jobs is asymptotically negligible. 
As discussed previously, the Halfin-Whitt asymptotic is not appropriate approximation for fluctuations that occur relative to system size, but does provide good engineering rules when stochastic fluctuations are $\ord{1}$ relative to system size. 
Thus we see that the NDS regime forms an appropriate mechanism to find new design rules and to analyse properties for policy selection.

\subsection{Join-the-Shortest-Queue with general processing time distributions} \label{sec:JSQ_PS}

Thus far we have been comparing centralized load balancing policies with exponential processing time distribution. In addition to such comparisons, NDS regime can also be used to uncover the effect of general processing time distributions on the performance of load balancing policies. 
This task would certainly require a substantial amount of technical work. (Work which we do not undertake in this paper.) As a demonstration, in this Section we show experimentally that an observation about JSQ load balancing with Processor Sharing servers made by \cite{bonomi1990job} and \cite{gupta2007analysis} can be established in the NDS regime, but is unlikely to be established in other conventional asymptotic regimes such as heavy-traffic and Halfin-Whitt. 

\begin{figure}
\centering
\subfigure[Simulation results of a finite system]{
\includegraphics[width=0.5\textwidth]{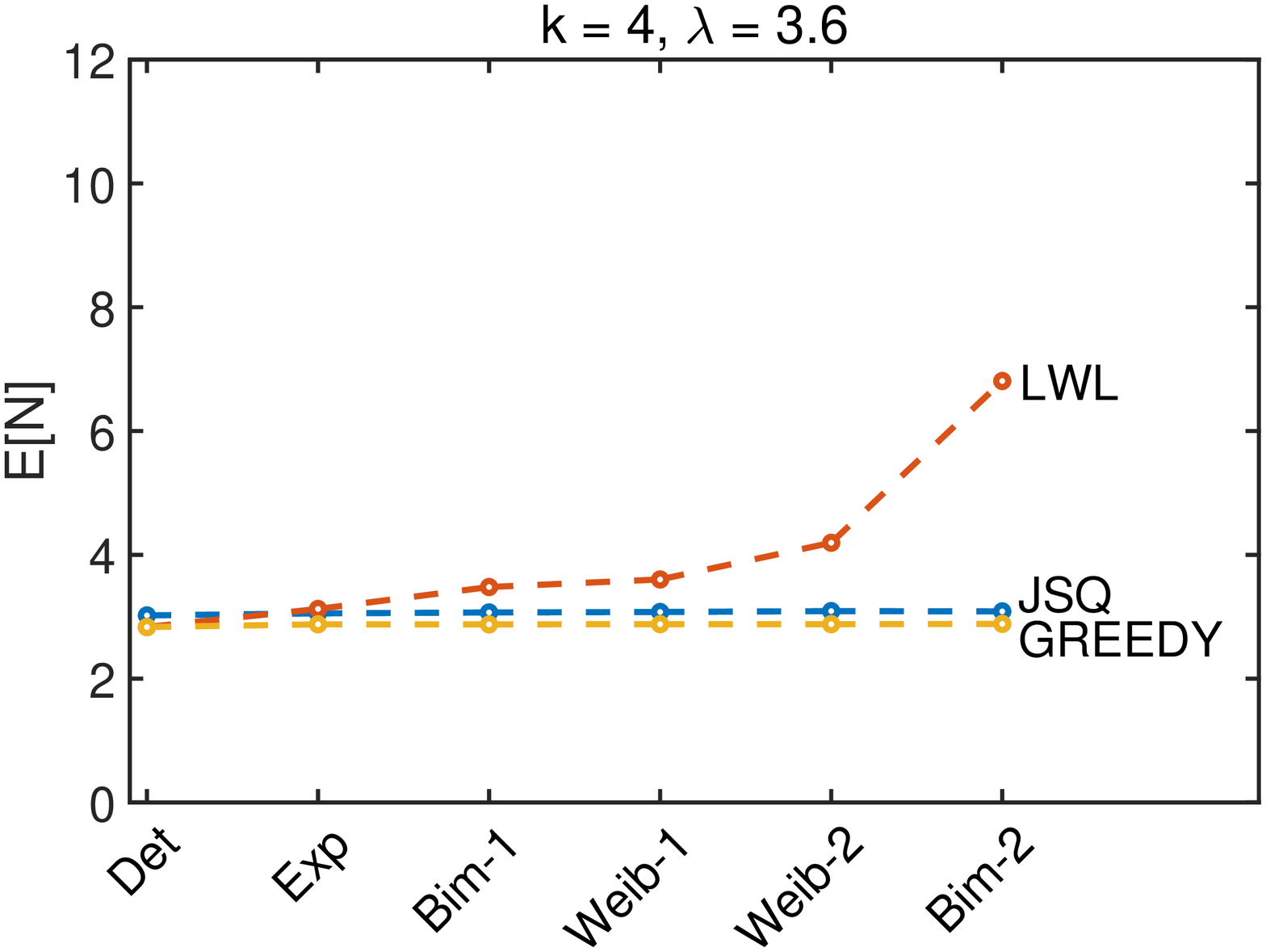}
\label{fig:JSQPS_finite}
}
\begin{minipage}{\textwidth}
\hspace{-0.3in}
\subfigure[NDS]{
\fbox{
\begin{minipage}{0.33\textwidth}
\includegraphics[width=1.12\textwidth]{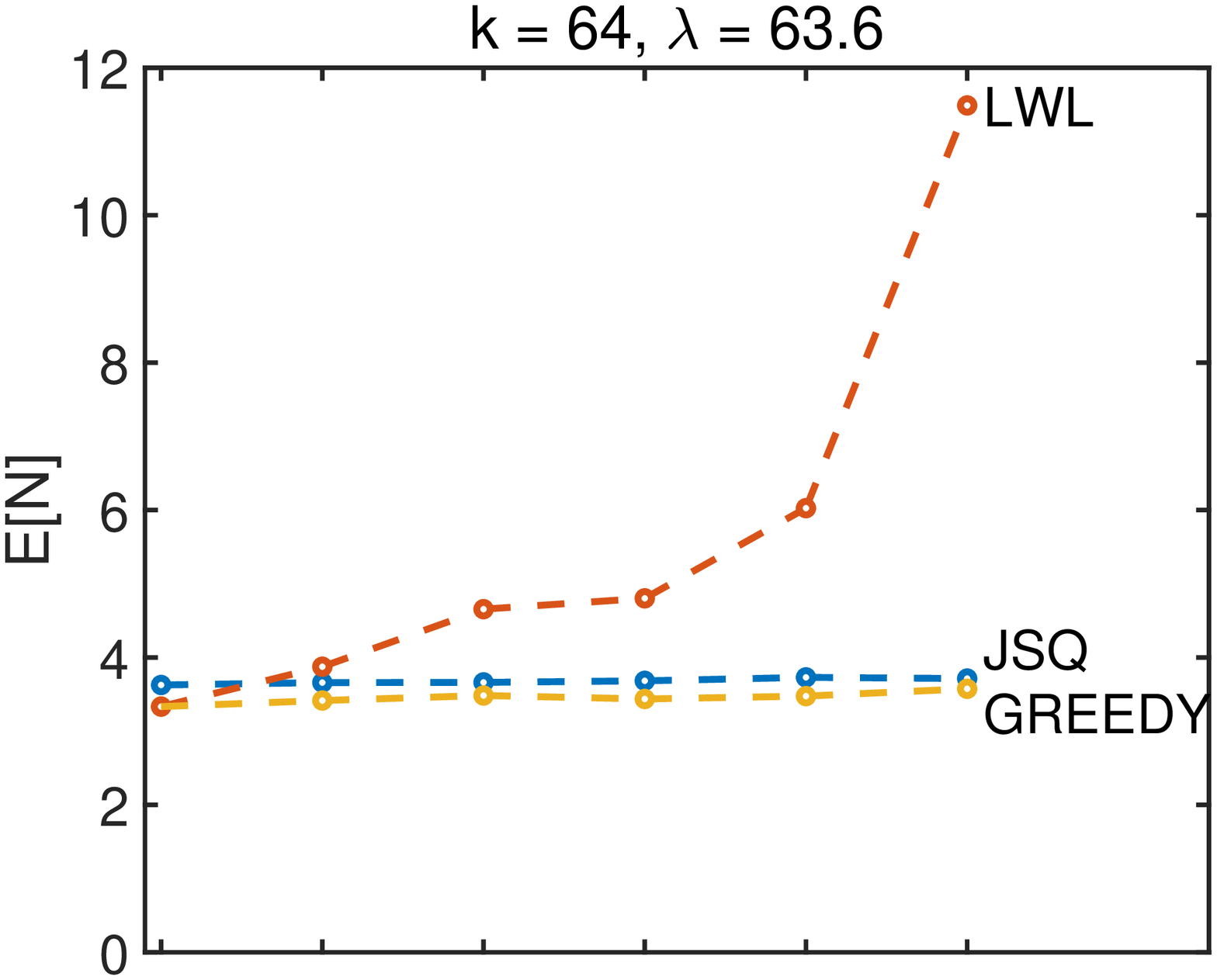}
\includegraphics[width=1.12\textwidth]{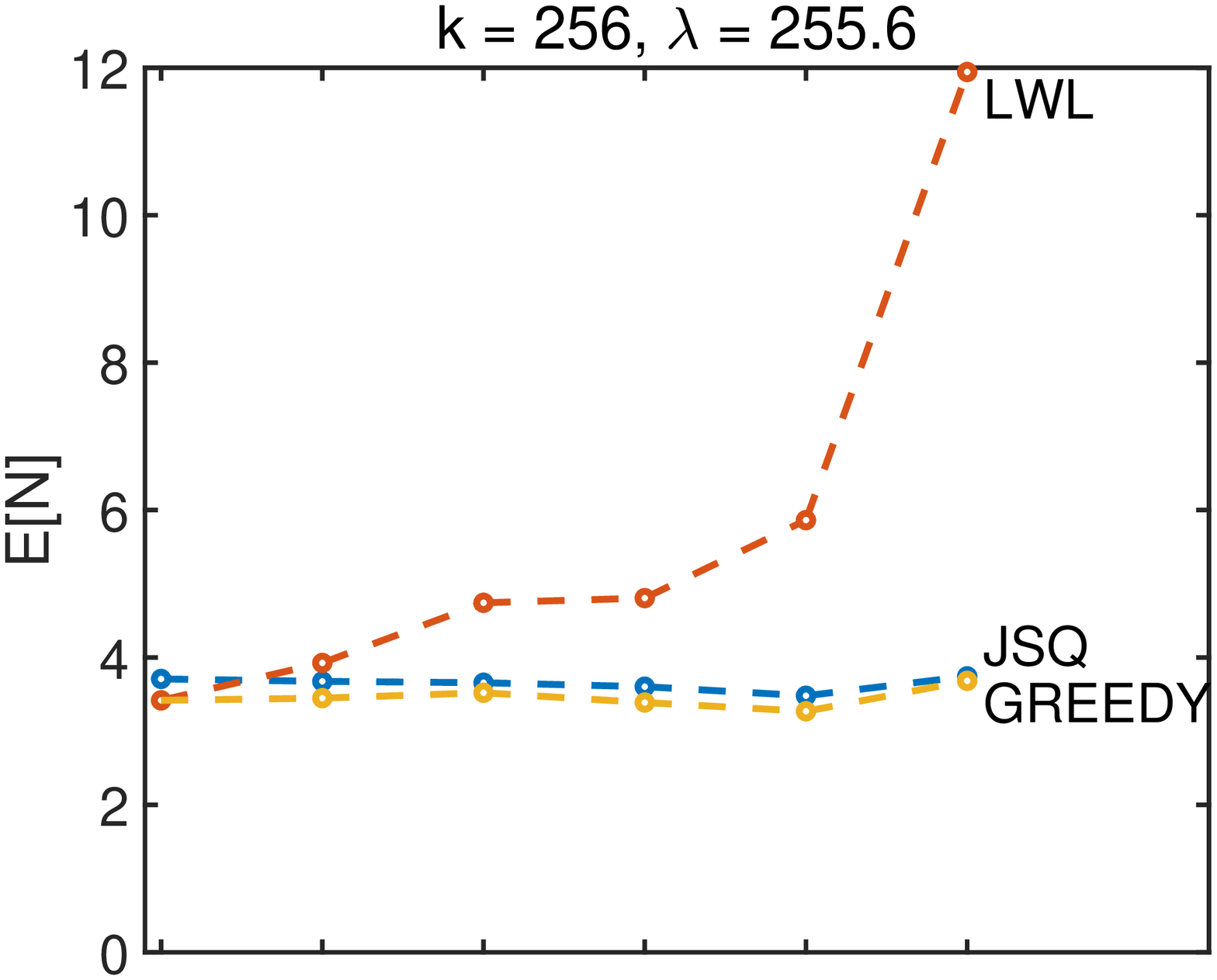}
\end{minipage}
}
\label{fig:JSQPS_NDS}
}%
\subfigure[Heavy-traffic]{
\fbox{
\begin{minipage}{0.33\textwidth}
\includegraphics[width=1.12\textwidth]{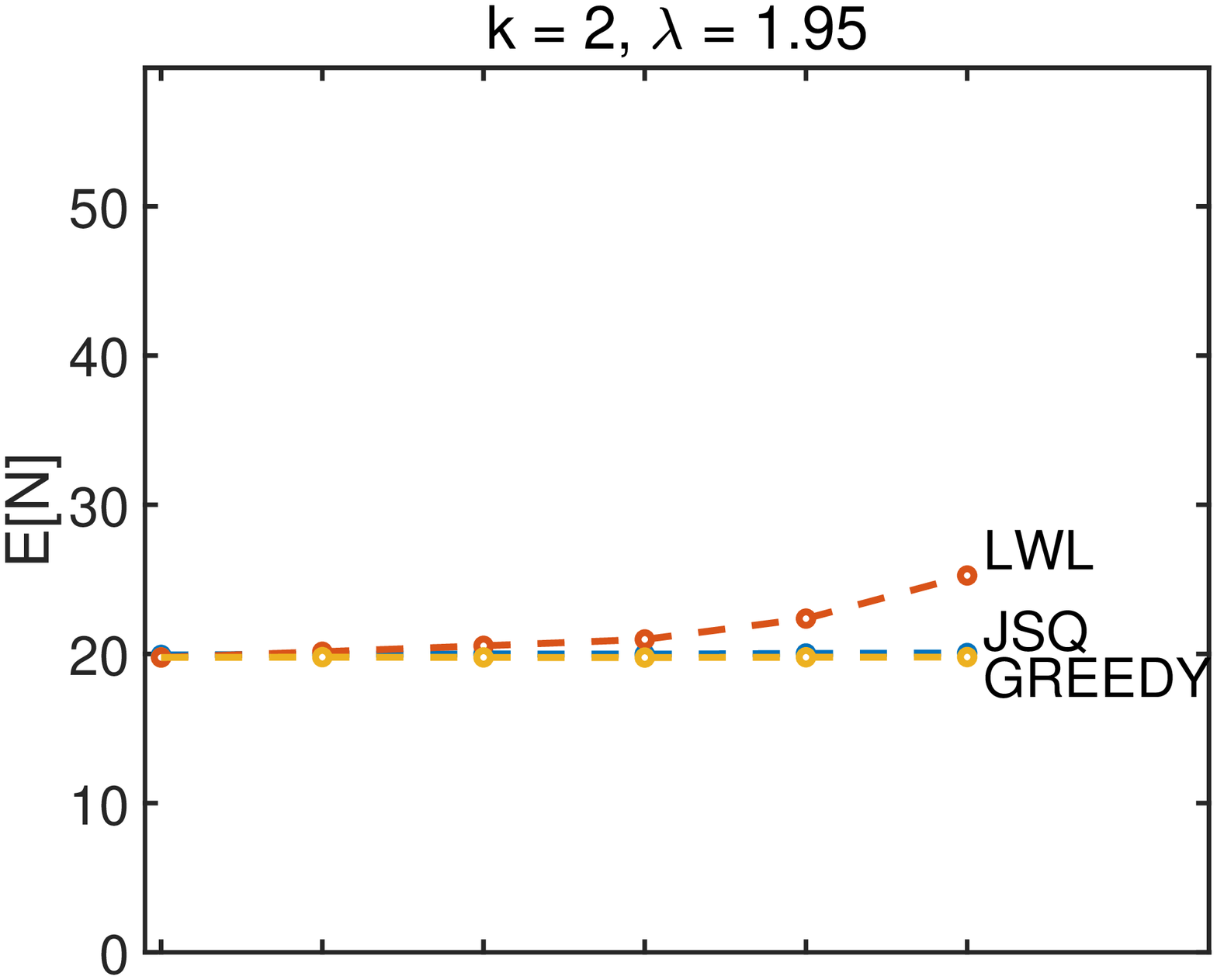}
\includegraphics[width=1.12\textwidth]{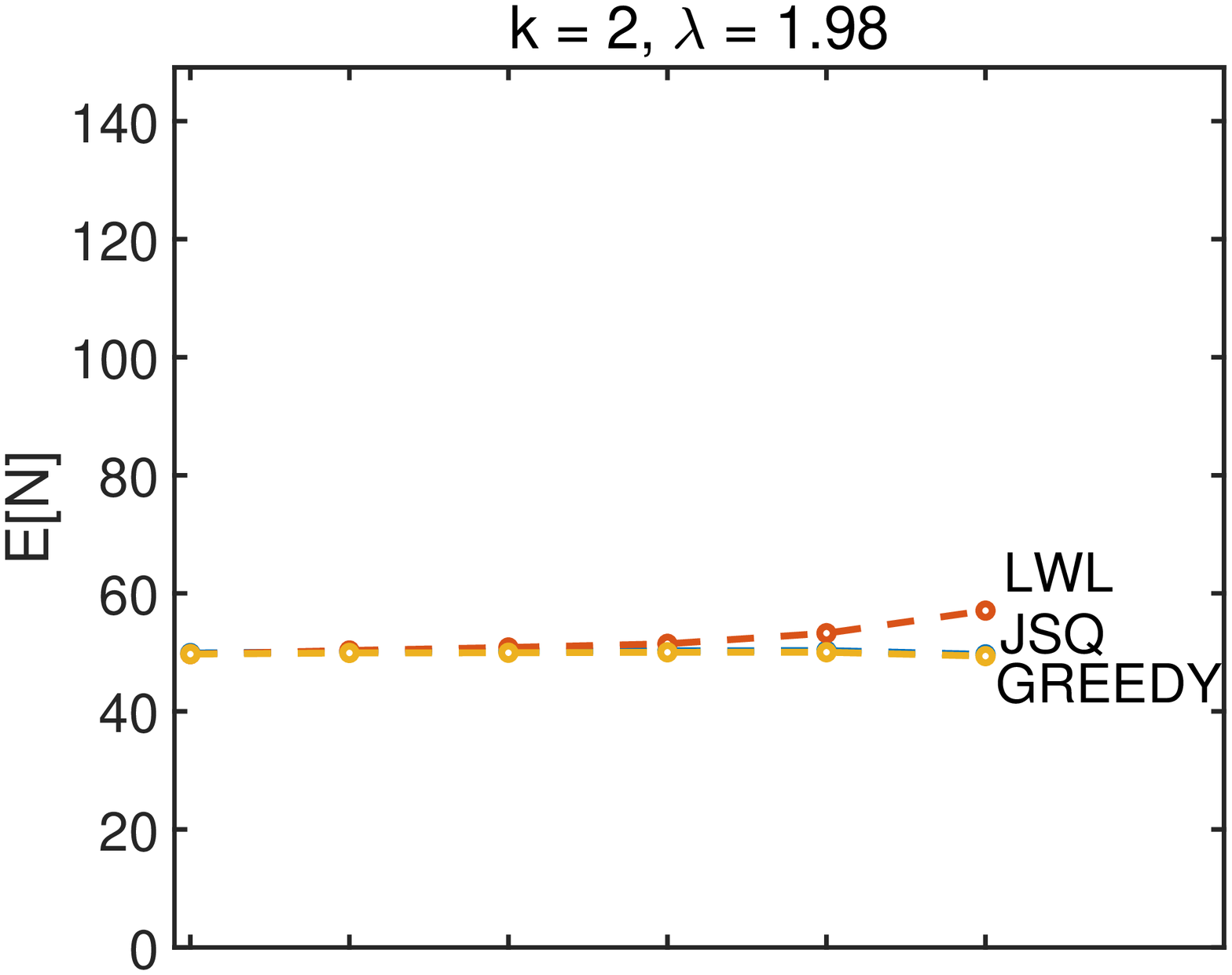}
\end{minipage}
}
\label{fig:JSQPS_HT}
}%
\subfigure[Halfin-Whitt]{
\fbox{
\begin{minipage}{0.33\textwidth}
\includegraphics[width=1.12\textwidth]{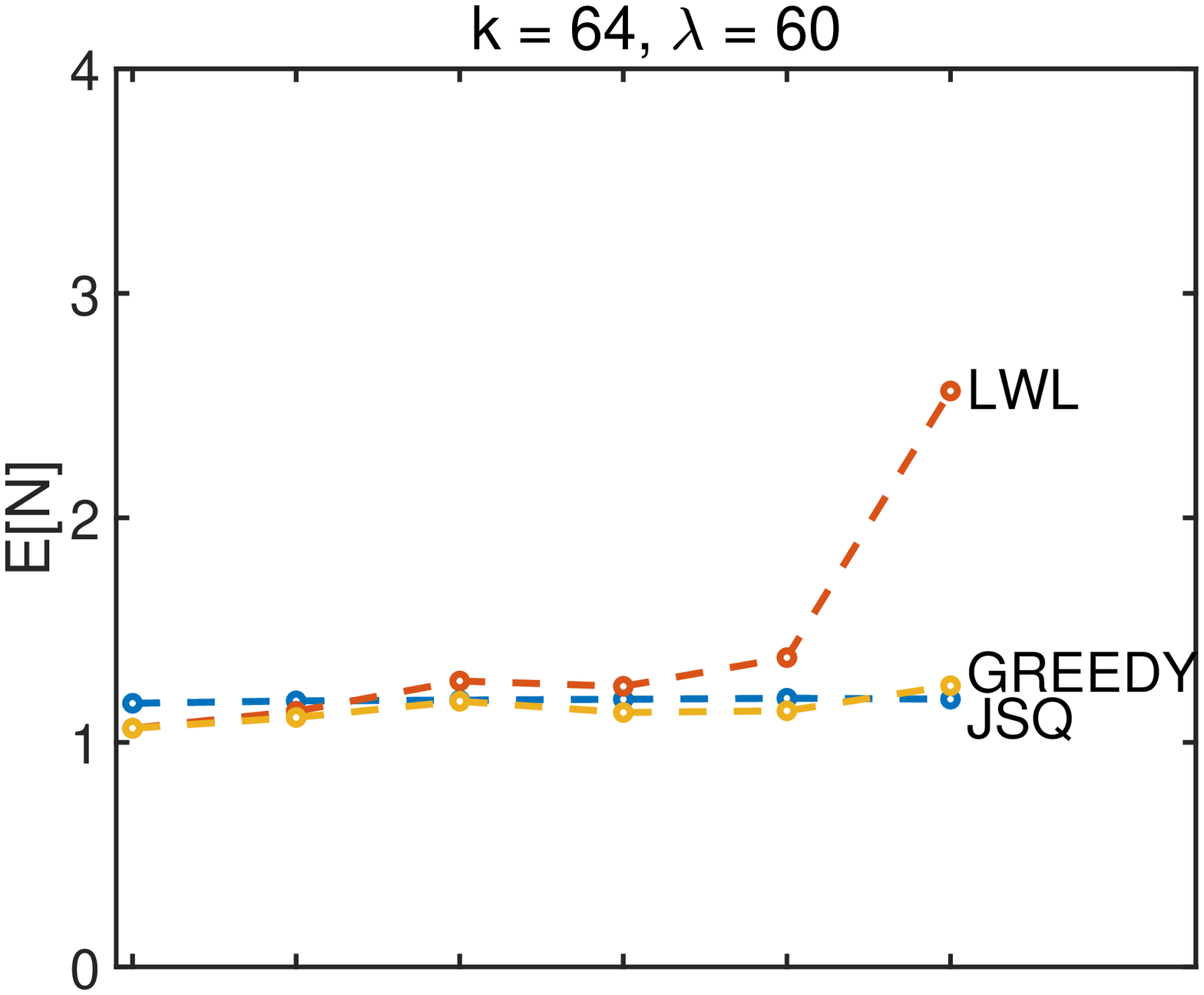}
\includegraphics[width=1.12\textwidth]{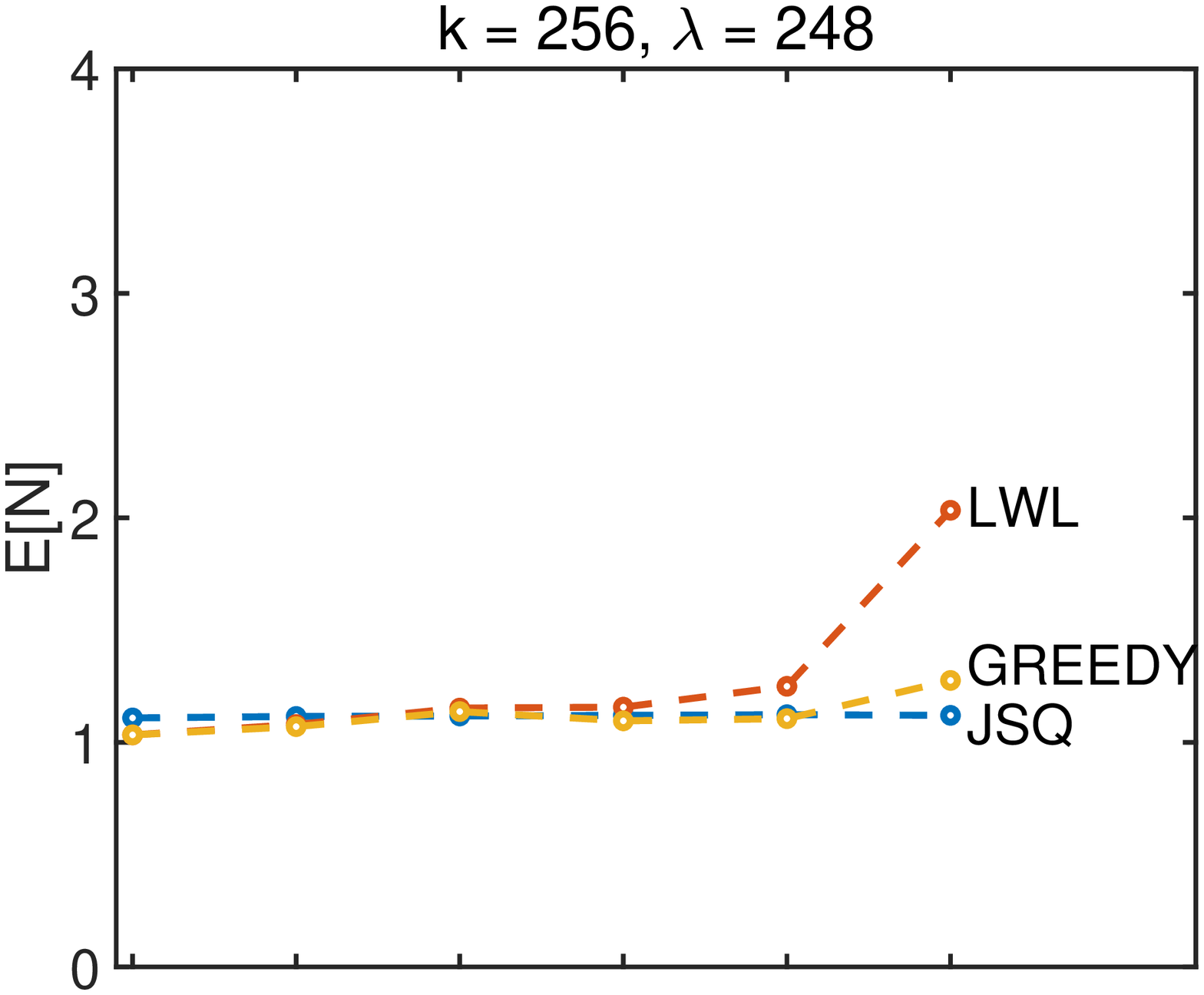}
\end{minipage}
}
\label{fig:JSQPS_HW}
}
\end{minipage}
\caption{{\bf Comparison of load balancing policies for general processing time distributions with Processor Sharing servers.} The $Y$-axis show the mean number of jobs per server, and the $X$-axis denotes the processing time distribution, all with mean 1 in increasing order of variance.
(a) Simulation for a finite system with $k=4$ servers, $\rho=0.9$. The remaining graphs are Pre-limit results for 
(b) NDS: number of servers increases from $k=64$ (top) to $k=256$ (bottom).
(c) Conventional heavy-traffic: $k=2$ while load is increased from $\rho=0.975$ (top) to $\rho=0.99$ (bottom). 
(d) Halfin-Whitt: Number of servers is increased from $k=64$ (top) to $k=256$ (bottom).
\label{figy}}
\end{figure}

Figure~\ref{fig:JSQPS_finite} shows simulation results for JSQ load balancing for a $k=4$ server system and six different processing time distributions ordered in increasing coefficient of variation on the $X$-axis. These distributions are listed in Table \ref{table}. The $Y$-axis shows the mean queue length per server. In addition to JSQ, simulations results are also shown for LWL (Least-Work-Left) load balancing which routes jobs to the server with the least total unfinished work, and a GREEDY policy which routes jobs to minimize the sum of sojourn times of all jobs currently in the system assuming there will be no further arrivals. Note that LWL and GREEDY require information about the exact processing times while JSQ does not. The values shown in Figure \ref{figy} are the means over 10 sample paths of $10^9$ arrivals each.

\begin{table}
\begin{center}
\begin{tabular}{l|ccc}
 Name & Distribution & Mean & Variance\\
  \hline
{\sf Det}:  & point mass at $1$   & $1$ & $0$ \\
{\sf Exp}: & exponential distribution & $1$ & $1$ \\
{\sf Bim-1}: & $X=0.5$  w.p. $0.9$ and  $X=5.5$ w.p. $0.1$ & $1$ & $2.25$\\
{\sf Weib-1}: & Weibull with shape~$=0.5$ and scale~$=0.5$ & $1$ & $5$ \\
{\sf Weib-2}: & Weibull with shape~$=\frac{1}{3}$ and scale~$=\frac{1}{6}$ & $1$ & $19$\\
{\sf Bim-2}: & $X=0.5$  w.p. $0.99$ and  $X=50.5$ w.p. $0.01$ & $1$ & $24.75$\\
\end{tabular}
\end{center}
\caption{\label{table}}
\end{table}

{\bf Insensitivity of JSQ load balancing: }  As was observed in \cite{bonomi1990job} and \cite{gupta2007analysis}, JSQ seems to exhibit some degree of insensitivity of the mean sojourn time to the processing time distribution, while LWL does not. If one were to conduct a formal study of the effect of processing time distribution on the performance of JSQ, what asymptotic regime will meaningfully establish this insensitivity? As we argue below, the answer is the NDS regime.

{\bf NDS regime:} Figure~\ref{fig:JSQPS_NDS} shows the NDS asymptotics by considering $k=64$ and $k=256$ server systems. As can be seen, as we increase $k$ in NDS regime, JSQ is insensitive (within simulation errors) while LWL is not. As a heuristic explanation: this insensitivity of JSQ again holds because of the separation of timescales between the queue length at a server (fast time scale) and the total number jobs in the system (slow timescale). On the faster timescale that jobs arrive and depart each server, $\hat{N}$ remains approximately fixed, and each processor sharing server will have a queue size that will fluctuate stochastically between a length $\lceil \hat{N} \rceil$ and $\lfloor \hat{N} \rfloor$. With $\hat{N}$ fixed, this stochastic process is reversible and insensitive. This follows from standard arguments, see \cite[Section 3.3]{kelly2011reversibility}. Consequently we anticipate that JSQ is insensitive in the NDS regime. 

{\bf Heavy-traffic regime:} Figure~\ref{fig:JSQPS_HT} shows the conventional Heavy-traffic asymptotics by considering $k=2$ servers with increasing values of $\lambda$. While JSQ is again insensitive in this asymptotic regime, so is LWL. This happens because the queue length for each server become large while the difference in queue lengths is of a smaller order. Therefore the entire system collapses approximately to a single PS server fed by a Poisson process -- an $M/G/1$ Processor Sharing queue which is known to be insensitive. Thus heavy traffic analysis will not distinguish between the insensitivity of JSQ and the sensitivity of LWL.

{\bf Halfin-Whitt regime:} Figure~\ref{fig:JSQPS_HW} shows the Halfin-Whitt asymptotics by considering $k=64$ and $k=256$. As the scale of the system increases, we see that $\expct{N}$ for all policies degenerates to 1 as we argued in the introduction. Thus again, this regime is not adequate in explaining the observations of Figure~\ref{fig:JSQPS_finite}.

\section{Conclusion}
\label{sec:conclusions}
This paper set out three primary goals: first, to argue that the many-servers NDS regime is a meaningful regime to study the performance of load balancing heuristics for finite systems; second, to employ the NDS regime to provide the first concrete analysis of the classical Join-the-Shortest-Queue (JSQ) dispatch rule; third and finally, to use this analysis to provide insight into more modern proxies for the JSQ policy, specifically, Power-of-$d$-choices, Idle-Queue-First and Idle-One-First.

Towards these goals, we emphasize the rationale behind the NDS regime: that the scaling regime is able to delineate properties that may not be present in alternative scalings such as heavy traffic, mean field and the Halfin-Whitt regime. 
We presented a first rigorous analysis of JSQ in the NDS regime. 
In doing so we find that a timescale separation between the number of idles servers and the mean number of jobs per servers informs much of the queue size behaviour of this policy. 
This in turn requires a number of novel mathematical approaches to provide a rigorous proof of convergence. On the other hand, the analysis reveals that the number of idle-servers \emph{and} the number of servers of length $1$ are the important factors which determine queue sizes for this policy. This informs our analysis of alternative proxies for JSQ.
We find that power-of-$d$ choices policies do not scale with load and system size. Idle-queue-first policy remains within a factor of $2$ of optimal policies (without jockeying) while JSQ lies within a factor of 1.15. We introduce the Idle-One-First policy that has identical behaviour to JSQ in NDS while having a lower communication overhead than Power-of-2-choices.

We advocate the use of NDS as a new regime for the analysis of multi-server queueing systems. 
Our analysis lends further credence to our claim that NDS regime faithfully replicates the performance of finite queueing systems where delay is of the same order of magnitude as processing time or where fluctuations in load are of the same magnitude as system size. 

{There are many aspects of this analysis that warrant further investigation.
A systematic theory that deals jointly with the time-scale separation and mean-field terms in the NDS diffusion asymptotic is required. Limit interchange results would be valuable;
Generalizations of the results taken here for general job sizes and complementary insensitivity proof would be very beneficial. Rigorous derivation of the stated NDS limits for Idle-Queue-First, Idle-One-First should be possible.
Study of more recent topics such as job replication, workload aware dispatchers, multiple dispatchers, resource pooling, and communication delay in NDS also require examination.}




\bibliographystyle{ormsv080} 
\bibliography{reference} 

\newpage
\begin{APPENDICES}
\section{Proof of Theorem~\ref{JSQ:Thrm}}

In this Appendix we prove our main result. A number of auxiliary results which are not specific to the line of investigation are contained in Appendix~\ref{sec:auxiliary_lemmas}.

\subsection{Additional Notation}\label{sec:additional_notation}

Our proofs will require the following additional notation.
%
%
For the $k$ server JSQ system, we let $M^{(k)}_l(t)$ denote the number of servers that have $l$ jobs at time $t$. We, also, define
\[
M^{(k)}_{\geq l}(t) = \sum_{l' \geq l} M^{(k)}_{l'}(t)\qquad \text{and} \qquad M^{(k)}_{\leq l}(t) = \sum_{l'\leq l} M^{(k)}_{l'}(t).
\]
Note that $I^{(k)}(t)=M^{(k)}_0(t)$. Further, we let
\[
N^{(k)}_{\geq l}(t) = \sum_{l' \geq l} (l'-l + 1) M^{(k)}_{l'}(t)
\]
Think of $N^{(k)}_{\geq l}$ the number of jobs above queue level $l$. Further, given the NDS scaling \eqref{NNDS}, we define
\[
\hat{N}^{(k)}_{\geq l}(t) = \frac{N^{(k)}_{\geq l}(kt) }{k},\qquad \hat{M}^{(k)}_{\geq l}(t) =\frac{M^{(k)}_{\geq l}(kt)}{k},\quad  \text{and}\quad \hat{M}^{(k)}_{\leq l}(t) = \frac{M^{(k)}_{\leq l}(kt)}{k}.
\]

\subsection{Outline of Theorem \ref{JSQ:Thrm}'s Proof.}\label{sec:outline}

The process $N^{(k)}$ can be represented in the following form
\begin{equation}\label{NExpression}
N^{(k)}(t) - N^{(k)}(0) = {\mathcal N}_a\Big( (k -\alpha)\mu t\Big) - {\mathcal N}_d \left( \int_0^t \mu k - \mu I^{(k)}(s)  ds \right), \qquad t \geq 0 
\end{equation}
where ${\mathcal N}_a$ and ${\mathcal N_d}$ are independent unit rate Poisson processes corresponding to arrivals and departures. We can subtract the mean from these Poisson processes.
Under the NDS scaling, this gives the following expression
\begin{subequations}\label{NTerms}
	\begin{align}
	\hat{N}^{(k)}(t)- \hat{N}^{(k)}(0)  =& {\mathcal M}_a\left( \mu t - \frac{\alpha}{k}\mu t \right) - {\mathcal M}_d\left( \mu t - \frac{\mu}{k^2} \int_0^{kt} I^{(k)}(s)ds  \right) \label{NTermsa} \\
	& + \frac{\mu}{k^2}\int_0^{k^2t} \hat{I}^{(k)}(s)ds - \mu \int_0^{t}  \frac{(2-\hat{N}^{(k)}(s))_+}{(\hat{N}^{(k)}_1(s)-1 ) } ds \label{NTermsb}\\
	&- \alpha \mu t + \mu \int_0^{t}  \frac{(2-\hat{N}^{(k)}(s))_+}{(\hat{N}^{(k)}_1(s)-1 ) } ds.\label{NTermsc}
	\end{align}
\end{subequations}
where
\[
\hat{I}^{(k)}(t) ={I}^{(k)}\Big(\frac{t}{k}\Big) ,\qquad {\mathcal M}^{(k)}_a (s) = \frac{{\mathcal N}_a (k^2 s)}{k} - k s\qquad \text{and}\qquad {\mathcal M}^{(k)}_d (s) = \frac{{\mathcal N}_d (k^2 s)}{k} - ks.
\]
Both ${\mathcal M}_a^{(k)}$ and ${\mathcal M}_d^{(k)}$ are Martingales.

We consider the convergence of each term in \eqref{NTerms} separately. We will argue that the first term \eqref{NTermsa} converges to a Brownian motion; the second term \eqref{NTermsb} converges to zero; and third \eqref{NTermsc} has sufficient continuity of solutions with respect to the arrival/departure process given in \eqref{NTermsa} to converge to the required SDE. We now give a little more detail in reference to the following subsections. 

First consider the term \eqref{NTermsa}. In Proposition \ref{lem:IdleBound}, it is shown that 
\[
\frac{1}{k^2} \int_0^{kt} I^{(k)}(s) ds \Rightarrow 0.
\]
This implies that arguments applied to $\mathcal{M}_a(\cdot)$ and $\mathcal{M}_d(\cdot)$ both weakly converge to the map $(t:t\geq 0) \mapsto (\mu t : t\geq 0)$. Thus, the Martingale Functional Central Limit Theorem gives that  \eqref{NTermsa} converges to a Brownian motion:

\[
 \left({\mathcal M}_a\left( \mu t - \frac{\alpha}{k}\mu t \right) - {\mathcal M}_a\left( \mu t - \frac{\mu}{k^2} \int_0^{kt} I^{(k)}(s)ds  \right) : t\geq 0\right)  \Rightarrow (\sqrt{2\mu}B(t) : t\geq 0) 
\]
where $B(t)$ is a standard Brownian motion. (A proof of the Martingale FCLT can be found in \cite{whitt2007} Theorem 2.1 and is stated in Theorem \ref{MgCLT}, below.) This deals with the term \eqref{NTermsa}.

Second we prove convergence to zero of \eqref{NTermsb}. Much of this article's technical novelty lies in proving the time-scale separation that 
\begin{align}
\label{eqTSS}
\left(\frac{1}{k^2}\int_0^{k^2t} \hat{I}^{(k)}(s) ds -\int_0^t  \frac{(2-\hat{N}^{(k)}(s))_+}{(\hat{N}^{(k)}(s) -1) } ds\; :\; t\geq 0\right) \Rightarrow 0.
\end{align}
\noindent As discussed in Section \ref{sec:Heuristic} the idleness process, $\hat{I}^{(k)}$, is informed by the proportion of queues of size one, $\hat{M}^{(k)}_1$, which in turn is closely related to the mean queue size $\hat{N}^{(k)}$. We must obtain the limiting relationship between these quantities. 
This is divided into the following parts: 
\begin{itemize}
\item We show in Proposition \ref{Thrm:StateSpaceCollapse} that, with high probability,
\[
\sup_{0\leq t \leq T} \Big| \big(2-\hat{N}^{(k)}(t) \big)_+ - \hat{M}^{(k)}_1(t) \Big| = k^{-\frac{1}{2} + o(1)} .
\]
This implies that, with high probability, $\hat{I}^{(k)}$ undergoes transition rates of the form
\begin{subequations}\label{IMM1}
	\begin{align}
	&\hat{I}^{(k)}(t) \rightarrow \;\; \hat{I}^{(k)}(t) +1\;   &&\text{at rate }\quad \mu \hat{M}^{(k)}_1\Big(\frac{t}{k^2}\Big) =\mu\Big(2-\hat{N}^{(k)}\Big(\frac{t}{k^2}\Big) \Big)_+ + \mu k^{-\frac{1}{2} + o(1)}, &\\
	&\hat{I}^{(k)}(t) \rightarrow \Big(\hat{I}^{(k)}(t) -1\Big)_+\;   &&\text{at rate }\quad \mu \left( 1 -\frac{\alpha}{k}\right). &
	\end{align}
\end{subequations}
Section \ref{sec:PropEC2} is devoted to proving Proposition \ref{Thrm:StateSpaceCollapse}. 

\item We then prove that \eqref{eqTSS} holds in Proposition \ref{TimescaleIdleness}. To prove Proposition \ref{TimescaleIdleness}, we must show that the transitions \eqref{IMM1} are well approximated by the dynamics of an M/M/1 queue. To do this we first refine Proposition \ref{Thrm:StateSpaceCollapse} to show that the result holds over each excursion of the process $\hat{I}^{(k)}$. This is Proposition \ref{TightProp} proven in subsubsection \ref{subsec:Tight}. We then bound excursions and their transition rates above and below in Section \ref{BoundingI}. Finally we combine these in Section \ref{sec:Prop3Proof} where Proposition \ref{TimescaleIdleness} is proven. 
\end{itemize}
This verifies the convergence of \eqref{NTermsb}.

Third and final, we show convergence to the required SDE and analyze its properties.  
The processes \eqref{NTerms} and their limit diffusion are sequence of jump/differential equations of the form
\[
n(t) = n(0) + \sqrt{2\mu} b(t) + \mu \int_0^t \left[  \frac{(2-n(s))_+}{(n(s)-1)} - \alpha \right] ds.
\]
In this sequence, the terms corresponding to $b(t)$ converges to a Brownian motion. Similar to \cite{pang2007martingale}, we can argue that the solution $n$ to the above equation is a continuous function of $b$. This is proven in Lemma \ref{CTS_SDE} below. Thus the sequence of processes \eqref{NTerms} converges to an SDE of the required form
\[
d\hat{N}(t) =  \sqrt{2\mu} dB(t) + \mu\left[  \frac{(2-\hat{N}(t))_+}{(\hat{N}(t)-1)} - \alpha \right] dt.
\]
The sample-path properties of this SDE are analyzed in Lemma \ref{DiffZero}.
The stationary density of the above SDE is found in Proposition \ref{JSQ:EqPro}. The stationary expected queue size is a routine calculation.

These arguments, as outlined above, then give the proof of Theorem \ref{JSQ:Thrm}.
\begin{remark}
The above argument is made somewhat complicated by the discontinuity in coefficients of the above SDE at $\hat{N}=1$. However, it can be argued that $\hat{N}$ never hits $1$. This holds since at $\hat{N}=1$ the diffusion coefficient of the above SDE is similar to that of a Bessel Process of dimension $2$ -- a process which is known to never hit zero, cf. \cite[Theorem V.40.1.]{rogers2000diffusions} and Proposition \ref{DiffZero}. This allows us to prove weak convergence under an appropriate localizing sequence of stopping times, $(\tau^{(k)}(\eta) : k\geq 0)$, which then implies the required weak convergence. Further, it implies path-wise uniqueness of our SDE.
\end{remark}

\subsection{A First Idleness Bound}

The lemma below shows that the number of idle servers is always less than $k^{1/2+o(1)}$. This is then sufficient to prove convergence of the Brownian term found in the NDS limit. 

\begin{proposition}\label{lem:IdleBound}
	For $\epsilon>0$,
	\[
	\bP \left( \sup_{0\leq t\leq kT} \Ik(t) \geq k^{\frac{1}{2}+\epsilon} \right) \leq a k^2 e^{-b k^{\epsilon}} 
	\]
	where $a$ and $b$ are positive constants depending on $\mu$, $\alpha$ and $T$ only. Consequently,
	\[
	\frac{1}{k^2}\int_0^{kt} I^{(k)}(s)ds \Rightarrow 0
	\]
where convergence is uniform on compact time intervals.
\end{proposition}

Notice this Lemma implies that, when placed in the same probability space, a sequence of JSQ networks will eventually never see a number of idle servers above $k^{\frac{1}{2}+\epsilon}$. 
\proof{Proof}
Given the number of queues of length $1$, the number of idle servers, $\Ik$, has transition rates
\begin{align*}
\Ik \nearrow \Ik +1 \quad &\text{at rate }\quad \mu \Mk_1\\
\Ik \searrow \Ik -1 \quad &\text{at rate }\quad \mu(k-\alpha) \quad \text{for}\quad \Ik > 0.
\end{align*}
Since $M^{(k)}_1 \leq k - I^{(k)}$, we can stochastically dominate the process $I^{(k)}$ by the Markov chain $\tilde{I}$ given by the transition rates
\begin{align*}
\tilde{I} \nearrow \tilde{I} +1 \quad &\text{at rate }\quad \mu( k- \tilde{I})\\
\tilde{I} \searrow \tilde{I} -1 \quad &\text{at rate }\quad \mu(k-\alpha) \quad (\text{when}\quad \tilde{I} > 0).
\end{align*}
By Lemma \ref{HittingLemma}, the probability of $\tilde{I}$ hitting $x$ before return to zero from state one is given by
\begin{align*}
h^x_1 = \frac{1}{\sum_{n=1}^x \prod_{m \leq n-1} \left( \frac{ k-\alpha}{ k - m}\right)} \label{Tagh}
&\leq 
\prod_{m \leq x-1} \left( \frac{  k - m}{ k-\alpha}\right) \\
&\leq \exp\left\{ -\sum_{m\leq x-1}\left( \frac{m}{ k} - 2 \frac{\alpha}{ k} \right) \right\}
= \exp \left\{ -\frac{x(x-1)}{2 k} + \frac{2\alpha (x-1)}{ k}   \right\} \notag ,
\end{align*}
where, in the first inequality, we take only the lead term from the denominator of $h_1^x$; and, in the second inequality, we apply the bound $e^{-2z}\leq 1-z \leq e^{-z}$ for $0\leq z\leq 1/2$. (Here we assume that $k$ is sufficiently large so that $\alpha/ k< 1/2$).  
Thus if we take $x= k^{\frac{1}{2}+\epsilon}$, for $\epsilon >0$, we have that 
\begin{equation}\label{excurisionbound}
h^{x}_1 \leq e^{-\frac{k^{2\epsilon}}{2}+ o(1)} \leq e^{-\frac{k^{\epsilon}}{2}+ o(1)}.
\end{equation}

This bounds each individual excursion. Now we must bound the possible number of such excursions.
Since the rate of creating idle servers is given by departures which is bounded above by a Poisson process of rate $\mu k$. Therefore the number of excursions in $\Ik$ from zero in time interval $[0,kT]$ is bounded above by, $Po(\mu T k^2)$, a Poisson random variable of mean $\mu T k^2$. By Lemma \ref{PoissonBound}
\begin{equation}\label{PoissonBoundy}
\bP( Po(\mu T k^2) \geq 8\mu T k^2 ) \leq e^{-9 \mu T k^2} .
\end{equation}
Thus applying a union bound, the bounds \eqref{excurisionbound} and \eqref{PoissonBoundy} give that 
\[
	\bP \left( \sup_{0\leq t\leq kT} \Ik(t) \geq k^{\frac{1}{2}+\epsilon} \right)
=\bP\Big( \exists t\leq kT\; s.t.\; I^{(k)}(t) \geq k^{\frac{1}{2}+\epsilon} \Big) \leq 8\mu Tk^2 e^{-\frac{k^{\epsilon}}{2\mu}+ o(1)} +   e^{-9 \mu T k^2} \leq a e^{-bk^{\epsilon}} 
\]
where the final inequality holds for suitable choices of $a$ and $b$. This is the bound required in Proposition \ref{lem:IdleBound}.
The convergence 
\[
\frac{1}{k^2}\int_0^{kt} I^{(k)}(s)ds \Rightarrow 0
\]
 is now a straightforward consequence. For $\epsilon \in (0,\frac{1}{2})$, 
\begin{align*}
{\mathbb P} \left( \frac{1}{k^2}\int_0^{kT} I^{(k)}(s) ds \geq \epsilon \right) & \leq {\mathbb P} \left( \frac{T}{k} \cdot \!\!\sup_{0\leq t \leq kT} I^{(k)} (t) \geq \epsilon \right)\\
& \leq   {\mathbb P} \left(  \sup_{0\leq t \leq kT} I^{(k)}(t) \geq k^{\frac{1}{2}+\epsilon}\right)\xrightarrow[k\rightarrow \infty]{} 0.
\end{align*} 
In the 2nd inequality above, we use the fact that for all $k$ sufficiently large, $\epsilon > T k^{-1/2+ \epsilon}$.
\hfill \Halmos
\endproof

A consequence of the above proposition is the convergence to the Brownian motion term in our limiting SDE.

\begin{corollary}\label{Cor:MtoB}
\[
\left({\mathcal M}_a\left( \mu t - \frac{\alpha}{k} \mu t \right) - {\mathcal M}_a\left( \mu t - \frac{\mu}{k^2} \int_0^{kt} I^{(k)}(s)ds  \right) : t\geq 0\right)  \Rightarrow (\sqrt{2\mu}B(t) : t\geq 0) 
\]	
	\end{corollary}
\proof{Proof}
By Proposition \ref{lem:IdleBound}ii) 
\[
\left(\mu t - \frac{\mu}{k^2} \int_0^{kt} I^{(k)}(s)ds : t\geq 0\right) \Rightarrow  \mu \mathbf{e}
\]
Thus the Martingale Functional Central Limit Theorem applies (see Theorem \ref{MgCLT} for a statement and see \cite[Theorem 2.1 and Lemma 2.1]{whitt2007} for a proof) and so the result holds.\hfill \Halmos
\endproof

\subsection{State Space Collapse}\label{sec:PropEC2}
In this subsection, we provide a sequence of results that prove the following proposition.  

\begin{proposition}[State Space Collapse]\label{Thrm:StateSpaceCollapse} For $\epsilon >0$
\[
\bP\Big( \sup_{0\leq t\leq T} \Big|\big(2-\hat{N}^{(k)}(t)\big)_+ - \hat{M}^{(k)}_1(t) \Big| \geq k^{-1/2+\epsilon}\Big) \leq a \exp\{-b k^{\frac{\epsilon}{2}}\}
\]
for positive constants $a$ and $b$.
\end{proposition}

The proposition ensures that the number of jobs $\hat{N}^{(k)}$ can provide good proxy for the number of servers of length one, which in turn informs the transitions of the number of idle servers.

The proof of Proposition \ref{Thrm:StateSpaceCollapse}. Relies on the following lemma. After proving this lemma, we outline the proof of Proposition \ref{Thrm:StateSpaceCollapse}.

\begin{lemma}\label{Lem:AMIdentity}
	\begin{equation}\label{eq:AMIdentity}
	\big(2-\hat{N}^{(k)}(t)\big)_+ - \hat{M}^{(k)}_1(t)= \left(\hat{M}_1^{(k)}(t) + 2\frac{I^{(k)}(kt)}{k} - {\hat{N}^{(k)}_{\geq 3}(t)}\right)_+\!\! - \hat{M}_1^{(k)}(t).
	\end{equation}
\end{lemma}
\proof{Proof}
The following identity counts the number of jobs for queue sizes above and below $2$ and will be useful for our analysis of JSQ: 
\begin{equation*}
2k =  N^{(k)}  + 2I^{(k)} + M^{(k)}_1 - N^{(k)}_{\geq 3}
\end{equation*}
(See Figure \ref{Fig:Identity} for pictorial representation of this identity).
\begin{figure}[t]%
	\centering
	\includegraphics[width=0.9\textwidth]{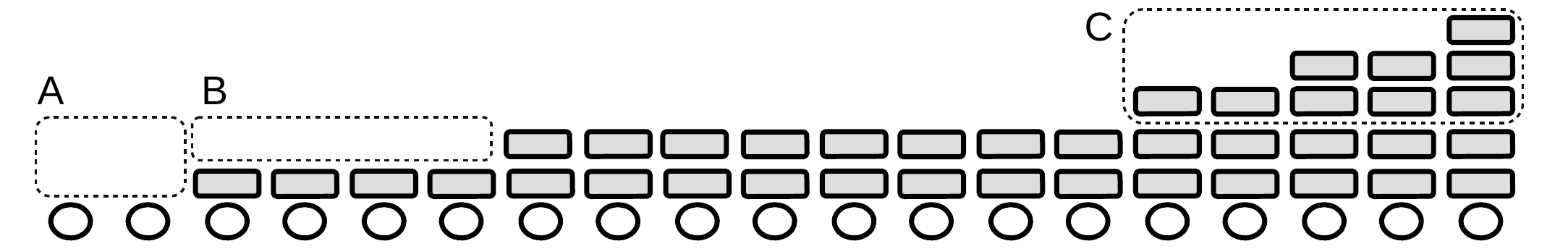}
	\caption{Representation of the identity $2k =  N^{(k)}  + 2I^{(k)} + M^{(k)}_1 - N^{(k)}_{\geq 3}$: $2k$ equals $N^{(k)}$ plus the number of jobs missing from box A which is $2 I^{(k)}$, plus the number missing from box B which is $M^{(k)}_1$, minus the number of jobs in box C which is $N^{(k)}_{\geq 3}$. 
		\label{Fig:Identity}}
\end{figure}%

The previous equality implies that
\begin{equation}\label{eq:MIdentity}
(2k - \Nk)_+ =  ( \Mk_1  +2\Ik - \Nk_{\geq 3})_+
\end{equation}
which after dividing by $k$ and subtracting $\hat{M}^{(k)}_1$ gives the required result. \hfill\halmos
\endproof

The proof of Proposition \ref{Thrm:StateSpaceCollapse} requires us to appropriately bound the right-hand side of \eqref{eq:AMIdentity} in Lemma \ref{Lem:AMIdentity}. First in Proposition \ref{lem:IdleBound}, we showed that $I^{(k)} = k^{\frac{1}{2}+o(1)}$ with high probability. This implies that the number of idle servers does not have a significant impact on Proposition \ref{Thrm:StateSpaceCollapse}. Consequently,  $(2-\hat{N}^{(k)})_+$ can only differ significantly from $\hat{M}_1^{(k)}$, when both $\hat{M}_1^{(k)}$ and $N^{(k)}_{\geq 3}$ are large. However, this is not possible under a JSQ rule. This is proven in Lemma \ref{lem:VBound} and Lemma \ref{lem:Nbound}, and this point was discussed heuristically in Section \ref{sec:Heuristic}.



\subsubsection{Further Idleness Bounds}\label{idlesness}
In this subsection we prove two somewhat more refined versions of Proposition \ref{lem:IdleBound} and each proof follows in a similar manner. Each lemma provides better bound on the terms in the identity \eqref{eq:AMIdentity}.
\begin{lemma}\label{lem:VBound}
	For $\epsilon>0$,
	\[
	\bP \left( \exists t \leq kT\text{ s.t. } (2\Ik (t) + \Mk_1(t) )   \geq k^{\frac{1}{2}}\text{ and }M^{(k)}_{\geq 3}(t) \geq k^{\frac{1}{2}+\epsilon} \right) \leq a_1  e^{-b_1 k^{\epsilon}} 
	\]
	where $a_1$ and $b_1$ are positive constants depending on $\mu$, $\alpha$ and $T$ only.	
\end{lemma}
\proof{Proof}
Define
\[
\Vk(t)= 2\Ik(t) + \Mk_1(t).
\]
(Notice this is the number of jobs that could be contained in boxes A and B in Figure \ref{Fig:Identity})

Notice that the transitions in $V^{(k)}(t)$ occur according to the following rates
\begin{align*}
\Vk \nearrow \Vk +1 \quad &\text{at rate }\quad \Mk_1 + \Mk_2 = k\mu - \mu \Mk_{\geq 3} - \mu \Ik \\
\Vk \searrow \Vk -1 \quad &\text{at rate }\quad \mu (k - \alpha)  \quad (\text{when}\quad \Vk > 0).
\end{align*}
Suppose that $M_{\geq 3} \geq y$. Notice that under JSQ $M_{\geq 3}$ increases only when $\Vk=0$. So we analyze the excursions of $\Vk$ from zero when $M_{\geq 3} \geq y$. In particular, we want to bound the probability that $\Vk\geq x$. Let ${\tilde{V}}$ be defined by the transition rates
\begin{align*}
\tilde{V} \nearrow \tilde{V} +1 \quad &\text{at rate} \quad k\mu -\mu y \\
\tilde{V} \searrow \tilde{V} -1 \quad &\text{at rate} \quad \mu (k -\alpha).
\end{align*}
By assumption excursions in $V^{(k)}$ are stochastically dominated by $\tilde{V}$. From Lemma \ref{HittingLemma}, we have that, when started from $\tilde{V}=1$, the probability of $\tilde{V}$ hitting $x$ before hitting zero is bounded as follows
\begin{align*}
h^x_1 = \frac{1}{\sum_{n=1}^x \prod_{m \leq n-1} \left( \frac{ k-\alpha}{  k - y}\right)} 
&\leq 
 \left( \frac{ 1 - \frac{y}{k}}{1-\frac{\alpha}{k}}\right)^{x-1} \leq \exp \left\{ -\frac{(x-1)y}{ k} + \frac{2\alpha (x-1)}{ k}   \right\} .
\end{align*}
where, in the first inequality, we take only the lead term from the denominator of $h_1^x$; and, in the second inequality, we apply the bound $e^{-2z}\leq 1-z \leq e^{-z}$ for $0\leq z\leq 1/2$.
Observe that if we let $x=k^{1/2}$ and $y=k^{1/2+\epsilon}$ then we have that 
\begin{equation}\label{hBound}
h_1^x \leq e^{-k^{\epsilon}+o(1)}.
\end{equation}

Again, this bounds each individual excursion. Now we must bound the possible number of such excursions.
Since the rate of increase in $\Vk$ is bounded above by a Poisson process of rate $\mu k$. Therefore the number of excursions in $\Vk$ from zero in time interval $[0,kT]$ is bounded above by, $Po(\mu T k^2)$, a Poisson random variable of mean $\mu T k^2$. Lemma \ref{PoissonBound} then gives that 
\begin{equation}\label{PoissonBound2}
\bP( Po(\mu T k^2) \geq 9 \mu T k^2 ) \leq e^{-10\mu T k^2} 
\end{equation}

Thus applying a union bound and the bounds \eqref{hBound}and \eqref{PoissonBound2}, we have that 
\[
\bP\Big( \exists t\leq kT\;\text{ s.t. }\; V^{(k)}(t) \geq k^{\frac{1}{2}}\text{ and } M^{(k)}_{\geq 3}(t) \geq k^{1/2+\epsilon}  \Big) \leq 9\mu Tk^2 e^{-\frac{k^{\epsilon}}{2\mu}+ o(1)} +   e^{-10 \mu T k^2} 
\]
It is not hard to see that for suitable choices of $a$ and $b$, we can achieve a bound of form required in Lemma \ref{lem:VBound}. \hfill \Halmos
\endproof

We now apply Lemma \ref{lem:VBound} to show that idleness occurs infrequently when $M_{\geq 3}$ is large. This will help us bound the maximum queue size in the next subsection. Again, the proof of Lemma \ref{2ndIdleBound} is a somewhat more sophisticated version of Proposition \ref{lem:IdleBound}.

\begin{lemma}\label{2ndIdleBound} For $\epsilon >0$, 
	\[
{\mathbb P}\left( \int_0^{kT} I^{(k)}(t) {\mathbb I} [ M_{\geq 3}(t) \geq k^{1/2 + \epsilon}]dt \geq k^{1/2+2\epsilon}\right)\leq a_2 {k^{2}} e^{-b_2 k^{\epsilon}}
	\]
	where $a_2$ and $b_2$ are positive constants depending on $\mu$, $\alpha$ and $T$ only.
\end{lemma}
\proof{Proof}
Given Lemma \ref{lem:VBound}, we consider the event 	
\[
A:=\left\{ \exists t \leq kT\text{ s.t. } (2\Ik (t) + \Mk_1(t) )   \geq k^{\frac{1}{2}}\text{ and }M^{(k)}_{\geq 3}(t) \geq k^{\frac{1}{2}+\epsilon}\right\}
\]
We know from Lemma \ref{lem:VBound} that ${\mathbb P}(A) \leq a_1k^2 e^{-b_1 k^{\epsilon}}$. 

Now consider $A^c$ the complement of $A$. On $A^c$, when $M^{(k)}_{\geq 3}(t)\geq k^{\frac{1}{2}+\epsilon}$ holds, it must be that
\begin{equation}\label{Mbound}
M_{1}^{(k)} \leq k^{1/2}.
\end{equation}
Recall the idleness process has transition rates  
 \begin{align*}
 &I^{(k)}\mapsto I^{(k)}+1 \qquad\text{at rate}\qquad \mu M_1^{(k)},\\
 &I^{(k)}\mapsto I^{(k)}-1 \qquad\text{at rate}\qquad \mu (k -\alpha), \quad\text{ when }\quad I^{(k)}>0.
 \end{align*}
So, given \eqref{Mbound}, we consider the CTMC 
 \begin{align*}
 &\tilde{I}\mapsto \tilde{I}+1 \qquad\text{at rate}\qquad \mu k^{1/2}\\
 &\tilde{I}\mapsto \tilde{I}-1 \qquad\text{at rate}\qquad \mu (k -\alpha), \quad\text{ when }\quad \tilde{I}^{(k)}>0.
 \end{align*}
Thus on event $A^c$, $I^{(k)}$ and $\tilde{I}$ are coupled so that $I^{(k)}(t)\leq \tilde{I}(t)$ for all $t$ such that $0\leq t \leq kT$.

We now bound the area under $(\tilde{I}(t) : 0\leq t \leq kT)$. We know that there are at most Poisson with mean $\mu k^{3/2}T$ excursions of $\tilde{I}(t)$ in time $0$ to $kT$ and, by Lemma \ref{PoissonBound}, we know that
\begin{equation}\label{thisPoissonBound}
{\mathbb P}\left( Po(\mu k^{3/2}T) \geq 9\mu k^{3/2} T \right) \leq e^{-10\mu T k^{3/2}}.
\end{equation}
Applying a union bound to these excursions, we have that 
\begin{equation}\label{yexpress}
{\mathbb P}  \left( \int_0^{kT}  \tilde{I} (t) dt \geq  y \right)
\leq
9\mu T k^{3/2 } {\mathbb P}  \left(  \int_0^{\tilde{T}} \tilde{I} (t)dt \geq \frac{y}{9\mu k^{3/2} T}\right) + e^{-10\mu T k^{3/2}}
\end{equation}
where $\tilde{T}$ is the length of an excursion of $\tilde{I}$.
(I.e. in the above union bound, either there are less that $9\mu k^{3/2}T$ excursions and then for the total area to be over $y$ there must be at least one excursion whose area is above $y/[9\mu k^{3/2} T]$, or, there are more than $9\mu k^{3/2}T$ excursions in which case \eqref{thisPoissonBound} is applied.) 

Now we apply Lemma \ref{ExcusionLemma} with parameter choices $\theta_0=1$, $c=\frac{k}{2}$, $x=k^{2\epsilon -1}$. (Notice $\phi(\theta_0))+ c \theta_0<0$ for all sufficiently large $k$). Thus we have that 
\[
{\mathbb P} \left( \int_0^{\tilde{T}} \tilde{I}(t) dt \geq k^{ 2\epsilon-1}  \right)
\leq
 e^{-\sqrt{\frac{\mu}{2}} (k^{\epsilon}-1)}
\]
Applying this to bound \eqref{yexpress} with $y=9\mu T k^{\frac{1}{2}+2\epsilon}$ gives
\[
{\mathbb P}  \left( \int_0^{kT}  \tilde{I} (t) dt \geq  9\mu T k^{\frac{1}{2}+2\epsilon} \right)
\leq
9\mu k^{3/2 }T   e^{-\sqrt{\frac{\mu}{2}} (k^{\epsilon}-1)} + e^{-10\mu T k^{3/2}}.
\]
Let $B$ be the event in the above probability (right-hand side above). We can now bound the integral of interest
\begin{align*}
{\mathbb P}\left( \int_0^{kT} I^{(k)}(t) {\mathbb I} [ M_{\geq 3}(t) \geq k^{1/2 + \epsilon}]dt \geq k^{1/2+2\epsilon}\right) &\leq {\mathbb P} (A) + {\mathbb P}(B\cap A^c)  \\
&\leq {\mathbb P} (A) + {\mathbb P} (B)\\
& \leq  a_1 k^2 e^{-b_1 k^{\epsilon}} + \left[ 9\mu k^{3/2 }T   e^{-\sqrt{\frac{\mu}{2}} (k^{\epsilon}-1)} + e^{-10\mu T k^{3/2}} \right].
\end{align*} 
From this inequality, it is clear that the required bound holds for suitably chosen parameters $a_2$ and $b_2$.
\hfill \Halmos	\endproof

\subsubsection{Bounding Queue Sizes}
To apply Lemma \ref{Lem:AMIdentity} and consequently prove Proposition \ref{Thrm:StateSpaceCollapse}, we must analyse $N_{\geq 3}^{(k)}$, the number of jobs queued in $3$rd position or above. 
However thus far in Section \ref{idlesness}, we have only gained bounds on $M^{(k)}_{\geq 3}$ the number of servers with $3$ or more jobs. To relate $M^{(k)}_{\geq 3}$ with $N^{(k)}_{\geq 3}$, we can apply bounds that bound the total queue size in NDS. (Our analysis here is specific  to JSQ. Other bounds might be used here to provide proofs for other dispatch rules.)

In this subsubsection the main result is the following lemma.

\begin{lemma}\label{lem:Nbound}
\[
{\mathbb P}\left( \max_{0\leq t\leq T} \left\{ \hat{N}^{(k)}(t)-\hat{N}^{(k)}(0) \right\}\geq n+3 \right) \leq a_3 e^{-b_3 \sqrt{n}}
\]
where $a_3$ and $b_3$ are positive constants that do not depend on $k$.
\end{lemma}

To prove this we bound the size of each excursion of the queue size process in Lemma \ref{Lemma8} and then bound the length of these excursions in Lemma \ref{Lemma9}. 

\begin{lemma}\label{Lemma8} 
For all times $s$ with $\hat{M}^{(k)}_{\geq 3}(s) \geq k^{-\frac{1}{2}+\epsilon}$ and $I^{(k)}(ks)=0$, let  $\tau$ be the next time that $\hat{M}^{(k)}_{\geq 3} < k^{-\frac{1}{2}+\epsilon}$ holds, specifically,
	$\tau = \min \left\{ t \wedge T : \hat{M}_{\geq 3}^{(k)}(t)< k^{-\frac{1}{2}+\epsilon} \right\}$. Given $s$ and $\tau$, the following bound holds
	\[
	\max_{s\leq t\leq \tau} \left\{ \hat{N}^{(k)}(t) - \hat{N}^{(k)}(s) \right\}\leq P.
	\]
where $P$ is a random variable such that
	\[
	\mathbb{P}\left(P\geq p+1\right) \leq 4 \exp\left\{-\frac{p^2}{8\mu T}\left( 1-\frac{p}{4\mu k T}\right)\right\}.
	\]
	\end{lemma}

\proof{Proof}
We know from \eqref{NExpression} that for $0\leq s \leq t\leq T$,
\begin{align*}
N^{(k)}(kt)- N^{(k)}(ks) =& {\mathcal N}_a \left( k\mu ( k -\alpha )t) -{\mathcal N}_a ( k\mu ( k -\alpha )s  \right) \\
& + {\mathcal N}_d \left( \int_0^t  \mu(k- I^{(k)}(u) )du \right) - {\mathcal N}_d \left( \int_0^s  \mu(k- I^{(k)}(u) )du \right)
\end{align*}
where $\mathcal{N}_a$ and $\mathcal{N}_d$ are unit Poisson processes. To prove the result we bound the excursions of $\hat{N}^{(k)}(t)$ given the bound on idleness Lemma \ref{2ndIdleBound} (which bounds the integral terms above).

Given $s$ is a time where $M^{(k)}_{\geq 3} (ks) \geq k^{\frac{1}{2} + \epsilon}$ and $I^{(k)}(ks)=0$.
(Note, each time $\max_{0\leq t'\leq t}\hat{N}^{(k)}(t')$ increases above $2+k^{1/2 +\epsilon}$, we must have $I^{(k)}(kt)=0$ and $M^{(k)}_{\geq 3} (ks') \geq k^{\frac{1}{2} + \epsilon}$.) Thus from Lemma \ref{2ndIdleBound} with high probability 
\begin{equation}\label{Intbound}
\int_{ks}^{kt} {I}^{(k)}(s)ds \leq k^{\frac{1}{2}+2\epsilon}
\end{equation}
for any $t$ such that $t\leq \tau$ where $\tau$ is the first time after $s$ such that $M_{\geq 3}^{(k)}(\tau) < k^{\frac{1}{2}+\epsilon}$
holds.
Thus defining $\bar{{\mathcal N}}(t)= {\mathcal N}(t)-t$, from \eqref{Intbound}, we have that 
\begin{align*}
&\hat{N}^{(k)}(t)-\hat{N}^{(k)}(s)\\ 
=& 
\frac{1}{k} 
{\mathcal N}_a \left( k\mu( k -\alpha) t \right) 
-
\frac{1}{k} 
{\mathcal N}_a \left( k\mu( k -\alpha) s \right)\\
&- \frac{1}{k} 
{\mathcal N}_d \left(  \int_0^{kt}  \mu(k- I^{(k)}(u) )du \right) 
+ 
{\mathcal N}_d \left( \int_0^{ks}  \mu(k- I^{(k)}(u) )du \right)
\quad\\
\leq &  
\frac{1}{k} 
{\mathcal N}_a \left( k\mu (k -\alpha) t \right) 
-
\frac{1}{k} 
{\mathcal N}_a \left( k\mu (k -\alpha) s \right)
\\
& 
- \frac{1}{k} 
{\mathcal N}_d \left(  \mu k^2 t -  k^{\frac{1}{2}+2\epsilon}- \mu \int_0^{ks}  I^{(k)}(u) du \right) 
+  
\frac{1}{k} 
{\mathcal N}_d \left(  \mu k^2 s - \mu \int_0^{ks}  I^{(k)}(u) du \right)
\\
=& 
\frac{1}{k} 
\bar{{\mathcal N}}_a \left( k\mu (k -\alpha) t \right) 
-
\frac{1}{k} \bar{{\mathcal N}}_a \left( k\mu (k -\alpha) s \right) - \alpha \mu (t-s)
\\
&- \frac{1}{k} 
\bar{{\mathcal N}}_d 
\left( 
	\mu k^2 t - k^{\frac{1}{2}+2\epsilon}
	- 
	\mu 
	\int_0^{ks}  I^{(k)}(u) du  
\right) 
+  
\frac{1}{k} 
\bar{{\mathcal N}}_d \left( \mu k^2 s - \mu \int_0^{ks}  I^{(k)}(u) du   \right) +k^{-\frac{1}{2}+2\epsilon}.
\end{align*}
We can maximize the left- and right-hand side of the above expression to give that
\begin{align*}
&
\max_{t : s\leq t\leq \tau}\left\{ \hat{N}^{(k)}(t) - \hat{N}^{(k)}(s) \right\}\\
\leq 
&\max_{0\leq s \leq t\leq T}
\Bigg\{  
\frac{1}{k} \bar{{\mathcal N}}_a \left( k\mu (k -\alpha) t \right) 
-\frac{1}{k} \bar{{\mathcal N}}_a \left( k\mu (k -\alpha) s \right) 
- \alpha (t-s)\\
&\qquad\qquad
- \frac{1}{k} \bar{{\mathcal N}}_d \left( \mu k^2 t - k^{\frac{1}{2}+2\epsilon}
- \mu \int_0^{ks}  I^{(k)}(u) du  \right) 
+  \frac{1}{k} \bar{{\mathcal N}}_d \left( \mu k^2 s - \mu \int_0^{ks}  I^{(k)}(u) du   \right) +k^{-\frac{1}{2}+2\epsilon} \Bigg\} \\
\leq & \left[  
\max_{0\leq t \leq T}\frac{1}{k} \bar{{\mathcal N}}_a \left( \mu k^2 t \right) 
- \min_{0\leq s \leq T}\frac{1}{k} \bar{{\mathcal N}}_a \left( \mu k^2 s \right) \right] 
+ \left[\max_{0\leq t \leq T} \frac{1}{k} \bar{{\mathcal N}}_d \left( \mu k^2  t \right) 
- \min_{0\leq s \leq T} \frac{1}{k} \bar{{\mathcal N}}_d \left( \mu k^2  s \right)\right] 
+ k^{-\frac{1}{2}+2\epsilon}=:P.
\end{align*}
So, we define $P$ to be the difference between the maximum and minimum of two Poisson processes. 
We can apply Lemma \ref{PoissonBound} to the 4 maximizations and minimization above. For instance,
\[
\mathbb{P} \left( \max_{0\leq t\leq T} \frac{1}{k}\bar{\mathcal{N}}_a(\mu k^2 t) \geq z \right) = \mathbb{P} \left(  \max_{0\leq t\leq \mu k^2 T} ( {\mathcal{N}}_a( t) - t ) \geq z k \right) \leq \exp\left\{-\frac{z^2}{2\mu T}\left(1-\frac{z}{\mu kT}\right)\right\}.
\]
Further we can assume that $k$ is sufficiently large so that $k^{-\frac{1}{2}+2\epsilon}<1$. Thus 
\[
\mathbb{P}\left(P\geq p+1\right) \leq 4 \exp\left\{-\frac{p^2}{8\mu T}\left( 1-\frac{p}{4\mu k T}\right)\right\}
\]
which gives a bound of the required form.

\hfill \Halmos
\endproof

The following lemma shows that every time there is an excursion away from $\hat{N}^{(k)}(0)=3$ then there is a strictly positive probability that $\hat{N}^{(k)}(t)$ will stay above $5/2$ until time $T$.

\begin{lemma}\label{Lemma9}
	Suppose that $\hat{N}^{(k)}(0)=3$ and $I^{(k)}(ks)=0$, let $\tau^{(k)}$ be the next time  such that $M^{(k)}_{\geq 3}(t)< k^{\frac{1}{2}+\epsilon}$ holds. Then there exists a strictly positive constant $q$ such for all $k$
	\[
	\mathbb{P} ( \tau^{(k)} > T) \geq q.
	\]
	\end{lemma}
\proof{Proof}

First observe
\begin{align}
\hat{N}^{(k)}(t) =& \hat{N}^{(k)}(0) + \frac{1}{k} {\mathcal N}_a \left( k(\mu k -\alpha) t \right) 
- \frac{1}{k} {\mathcal N}_d \left(  \int_0^{kt}  \mu(k- I^{(k)}(u) )du \right) \notag \\
\geq & 3+ \frac{1}{k} {\mathcal N}_a \left( k\mu( k -\alpha) t \right) 
- \frac{1}{k} {\mathcal N}_d \left(   \mu k^2 t \right) \label{InNHat} \\
=& 3+ \frac{1}{k} \bar{{\mathcal N}}_a \left( k \mu( k -\alpha) t \right) 
- \frac{1}{k} \bar{{\mathcal N}}_d \left(   \mu k^2 t \right) - \alpha \mu t.\notag
\end{align}
By the Martingale Functional Central Limit Theorem \ref{MgCLT} the last term converges to a 
Brownian motion with drift. That is
\[
\left(\frac{1}{k} \bar{{\mathcal N}}_a \left( k \mu ( k -\alpha) t \right) 
- \frac{1}{k} \bar{{\mathcal N}}_d \left(   \mu k^2 t \right) - \alpha \mu t : t\geq 0 \right)\Rightarrow \left(\sqrt{2\mu} B(t) -\alpha \mu  t : t\geq 0\right).
\]
Further, we know that 
\[
\bP \left( 3 + \sqrt{2\mu} B(t) -\alpha \mu t > \frac{5}{2}, \forall t \leq T \right) \geq 2q
\]
for some positive constant $q$. The above set is an open set with respect to the Skorohod topology on $[0,T]$. Thus by the Portmanteau Theorem applied to open sets (see \cite{billingsley2013convergence}): for all sufficiently large $k$
\begin{equation}\label{ProbHitBound}
\mathbb{P}\left(3+\frac{1}{k} \bar{{\mathcal N}}_a \left( k \mu (k -\alpha) t \right) 
- \frac{1}{k} \bar{{\mathcal N}}_d \left(   \mu k^2 t \right) - \alpha \mu t > 2\frac{1}{2} , \forall t\geq T \right) \geq q.
\end{equation}
Combining \eqref{InNHat} and \eqref{ProbHitBound} gives the required bound
\[
\mathbb{P}\left(\tau^{(k)} > T\right) = \mathbb{P}\Big( \hat{N}^{(k)}(t) > \frac{5}{2}, \forall t\leq T \Big)\geq q.  
\]
\hfill\Halmos
\endproof
We use the last lemma to obtain a bound on the maximum queue size under the JSQ policy.

\proof{Proof of Lemma \ref{lem:Nbound}} 

	Let $\sigma_0$ be the first time that $\hat{N}^{(k)}(t)=\lceil \hat{N}^{(k)}(0) \rceil + 3$.
	Let $\hat{\tau}_i$ be the first time after $\sigma_i$ that $\hat{N}^{(k)}(t)=1.5$
	and let $\sigma_i$ be the first time after $\hat{\tau}_{i-1}$ that $\hat{N}(t) = 3$. Applying Lemma \ref{lem:VBound} note that  $M_{\geq 3}^{(k)}(t)<k^{\frac{1}{2}+\epsilon}$ must have occurred. So given Lemma \ref{Lemma9}, $\hat{\tau}_0-\sigma_0 > \tau^{(k)}$.
	
	By Lemma \ref{Lemma9}, 
	\begin{equation*}
		\bP (\tau_i - \sigma_i \geq T) \geq q>0
	\end{equation*}
	and by Lemma \ref{Lemma9}
	\begin{equation*}
  \bP \left( \max_{\sigma_i \leq t \leq \tau_i} \hat{N}^{(k)}(t) \geq p+3 \right) 
  \leq 
	\bP \left( P \geq p +3 \right)
\end{equation*}
where $P$ is as specified in Lemma \ref{Lemma9}.

So
\begin{equation*}
  \max_{0\leq t \leq T} \hat{N}^{(k)}(t) 
  \leq
  \sum_{i=1}^{G(q)} P_i + \lceil \hat{N}^{(k)}(0) \rceil + 3
\end{equation*}
where $G(q)$ is geometrically distributed with parameter $q$ and $P_i$ is an independent random variable with distribution $P$. So applying a union bound
\begin{align*}
  \bP \left( \sum_{i=1}^{G(q)} P_i \geq n \right) 
  &
  \leq 
  \bP \left( G(q) \geq \sqrt{n} \right)
  +
  \bP \left( \sum_{i=1}^{\sqrt{n}} P_i \geq  n \right)
 	\\ &
   \leq 
    \bP \left( G(q) \geq \sqrt{n} \right)
  +
  \sqrt{n}\bP \left( P \geq  \sqrt{n} \right) \leq a_3 e^{-b_3 \sqrt{n}}
\end{align*}
where the final inequality holds for appropriate positive constants $a_3$ and $b_3$.
\hfill \Halmos

\endproof

\subsubsection{Proof of Proposition \ref{Thrm:StateSpaceCollapse}}

With Lemma \ref{Lem:AMIdentity},  Lemma \ref{lem:VBound}, and Lemma \ref{lem:Nbound},  we can now prove Proposition \ref{Thrm:StateSpaceCollapse}.
In essence Lemma \ref{lem:VBound} and Lemma \ref{lem:Nbound}, combine to say that we cannot have both $\hat{M}_1^{(k)}$ and $N^{(k)}_{\geq 3}$ large at the same time. This bounds the error in identity Lemma \ref{Lem:AMIdentity} and thus shows that $(2-\hat{N}^{(k)}(t))_+$ and  $\hat{M}^{(k)}_1(t)$ must be close.

\proof{Proof of Proposition \ref{Thrm:StateSpaceCollapse}}
Suppose that the following three events all hold:
\begin{subequations}
	\begin{align}
A 
&= 
\Big\{ I^{(k)}(t) \leq k^{\frac{1}{2}+\epsilon},\, \forall t\leq kT \Big\}\label{eq:Acond}
\\
B 
&= 
\Big\{ 
	(2\Ik (t) + \Mk_1(t) )   
	\leq 
	k^{\frac{1}{2}}, 
	\forall t \leq kT\text{ such that }M^{(k)}_{\geq 3}(t) \geq k^{\frac{1}{2}+\epsilon} \Big\}\label{eq:Bcond}
\\
C
&= 
\Big\{\max_{0\leq s\leq T}  \hat{N}^{(k)}(s) \leq k^{\epsilon}\Big\}
\end{align} 
\end{subequations}
Applying a union bound, we know by Proposition \ref{lem:IdleBound},  Lemma \ref{lem:VBound}, and Lemma \ref{lem:Nbound} that 
\[
\bP( A \cap B \cap C ) \geq 1- a \exp\{-b k^{\frac{\epsilon}{2}}\}
\]
where $a$ and $b$ are constants depending on $\mu$, $\alpha$ and $T$.

From Lemma \ref{Lem:AMIdentity}, we have the identity:
\begin{equation}\label{eq:NMIndentity}
\big(2-\hat{N}^{(k)}(t)\big)_+ - \hat{M}^{(k)}_1(t) = \left(\hat{M}_1^{(k)}(t) + 2\frac{I^{(k)}(kt)}{k} - \frac{N^{(k)}_{\geq 3}(kt)}{k}\right)_+\!\! - \hat{M}_1^{(k)}(t).
\end{equation}
We analyze this equality under two separate conditions: 1. $N^{(k)}_{\geq 3}(kt) \leq k^{1/2+2\epsilon}$ and 2. $N^{(k)}_{\geq 3}(kt) \geq k^{1/2+2\epsilon}$.

First, we assume the condition $N^{(k)}_{\geq 3}(kt) \leq k^{1/2+2\epsilon}$ is met. We note the following elementary bound: for $x>0$,
\begin{equation*}
  |(x+y)_+ - x| = | -(x+y)_- +y | \leq (x + y)_- + |y| \leq 2|y|. 
\end{equation*}
Applying this bound to Identity \eqref{eq:NMIndentity} gives that
\[
\big| \big(2-\hat{N}^{(k)}(t)\big)_+ - \hat{M}^{(k)}_1(t) \big|
\leq 
2\Bigg|2\frac{I^{(k)}(kt)}{k} - \frac{N^{(k)}_{\geq 3}(kt)}{k}\Bigg| 
\leq 
4 k^{-\frac{1}{2}+\epsilon}
\] 
in the final inequality, above, we apply our condition, $N^{(k)}_{\geq 3} (kt) \leq k^{1/2+2\epsilon}$,  and our assumption that Event $A$ holds. 
So for the first case the required bound holds.

Second, assume that the condition ${N}^{(k)}_{\geq 3} (kt) \geq k^{1/2+2\epsilon}$ is met. We know that $1 + \max_{0\leq s\leq T}  \hat{N}^{(k)}(s)$ gives the maximum queue size achieved at any server. Thus,
\begin{equation*}
   {M}_{\geq 3}(kt) \times \max_{0\leq s\leq T}  \hat{N}^{(k)}(s)
   \geq 
    {N}_{\geq 3}^{(k)}(kt).
\end{equation*}
Thus, since Event C holds, we have that 
\begin{equation*}
  {M}_{\geq 3}(kt)
  \geq 
  \frac{
  		{N}^{(k)}_{\geq 3}(kt)
  		}{
  		\max_{s\leq T} \hat{N}^{(k)}(s)
  		}
  \geq  
  \frac{
  		k^{1/2+2\epsilon}
  		}{
  		k^{\epsilon}
  		}
  	=
  	k^{1/2+\epsilon}.
\end{equation*}
By this bound and since Event $B$ holds, we have 
\begin{equation*}
  (2\Ik (kt) + \Mk_1(kt) )   
	\leq 
	k^{\frac{1}{2}}.
\end{equation*}
Now applying this to Identity \eqref{eq:NMIndentity} gives
\begin{equation*}\label{eq:NMIndentity}
\Big|\big(2-\hat{N}^{(k)}(t)\big)_+ - \hat{M}^{(k)}_1(t) \Big| 
	\leq 
	\left(
		\frac{{M}_1^{(k)}(kt)}{k} 
		+ 
		2\frac{I^{(k)}(kt)}{k}
	\right)_+
		\!\! 
		+ 
		\frac{{M}_1^{(k)}(kt)}{k} \leq 2 k^{-1/2}.
\end{equation*}
as required.\hfill\halmos
\endproof

\subsection{Timescale Separation}

The main aim of this subsection is to prove Proposition \ref{TimescaleIdleness} as stated below. This addresses the convergence of term \eqref{NTermsb} as discussed in the outline of proof of Theorem \ref{JSQ:Thrm} in Section \ref{sec:outline}.

From Proposition \ref{Thrm:StateSpaceCollapse} we now see that $\big(2-\hat{N}^{(k)}(t)\big)_+$ provides a good estimate for $\hat{M}^{(k)}_1(t)$, the proportion of queues of length $1$. We present a slightly refined version of this result in Section \ref{subsec:Tight}. This then dictates the transitions of the idleness process which we analyse and bound in Section \ref{BoundingI}. After this we are in a position to prove Proposition \ref{TimescaleIdleness}.  

Our proof is slightly complicated by the necessity to provide a localizing sequence for our weak convergence result. Throughout this subsection we let $\tau^{(k)}(\eta)$ be the first time before time $T$ such that $\hat{N}^{(k)}(t)\leq 1+\eta$ holds. That is
\begin{equation*}
  \tau^{(k)}(\eta) = \min\{ t \leq T : \hat{N}^{(k)}(t)\leq 1+\eta  \} \wedge T
\end{equation*}

We analyze the number of idle servers, $\hat{I}^{(k)}$, over excursions from zero.
Let $T^{(k)}_m$ be the sequence of times, less than $k^2 \tau^{(k)}(\eta)$ where $\hat{I}^{(k)}$ hits zero. Specifically, $T_0^{(k)}=0$ and, for $m\in\mathbb{N}$,
\[
T^{(k)}_m 
= 
\inf\big\{ 
t\geq X_{m-1}^{(k)} 
: 
\hat{I}^{(k)}(t)=0 
\big\}
\wedge
\big\{ k^2 \tau^{(k)}(\eta) \big\}.
\]
In the above, expression $X_{m-1}^{(k)} $ is the next time after time $T_{m-1}^{(k)} $  where $\hat{I}^{(k)}(t)>0$ holds. 

As discussed above, we wish to prove the following Proposition. 

\begin{proposition} \label{TimescaleIdleness}
\[
	\sup_{0\leq t\leq \tau^{(k)}(\eta)} 
		\left| 
			\frac{1}{k^2} 
				\int_{0}^{k^2t}	\hat{I}^{(k)}(s)ds 
			- 
				\int_{0}^{t}	
					\frac{(2-\hat{N}^{(k)}(s))_+}{(\hat{N}^{(k)}(s)-1 )}
				ds 
			\right| 
	\xrightarrow[k\rightarrow \infty]{} 
	0
\]
where the above convergence is convergence in probability.
\end{proposition}

\subsubsection{Tightness bound}\label{subsec:Tight}
We wish to show that $\hat{I}^{(k)}$  behaves approximately as an M/M/1 queue. 
 We know that the $\hat{I}^{(k)}$ observes transition rates
 \begin{align*}
 &\hat{I}^{(k)}\mapsto \hat{I}^{(k)}+1 \qquad\text{at rate}\qquad \mu \hat{M}_1^{(k)}(t),\\
 &\hat{I}^{(k)}\mapsto (\hat{I}^{(k)}-1)_+ \quad\text{at rate}\qquad \mu \left( 1 -\frac{\alpha}{k} \right).
 \end{align*}
By Proposition \ref{Thrm:StateSpaceCollapse}, $\hat{M}_1^{(k)}(t)$ is well approximated by $\big(2-\hat{N}^{(k)}(t)\big)_+$. In this subsubsection, to show that these terms are approximately constant over excursions of $\hat{I}^{(k)}$. Thus, in the next subsubsection, we discuss how $\hat{I}^{(k)}$ behaves approximately as an M/M/1 queue.

We prove the following proposition, which is the main result in this subsubsection. It shows that $\hat{N}^{(k)}$ cannot oscillate greatly between excursions of $\hat{I}^{(k)}$.

\begin{proposition}\label{TightProp}
\[
	\mathbb{P} 
		\left( 
			\sup_{m\in\mathbb{Z}_+} 
			\sup_{T^{(k)}_m < t \leq T^{(k)}_{m+1}  } \left| 
					\hat{M}_1^{(k)}
						\left(\frac{t}{k^2} \right)
					-
					\left( 
							2
							-
							\hat{N}^{(k)}\left(\frac{T^{(k)}_m}{k^2} \right) 
					 \right)_+
				 \right| 
			\geq 
			2 k^{-1+\epsilon}   
		\right) 
	\leq 
	\gamma e^{-\beta k^{\epsilon/2}}.
\]
\end{proposition}
To prove this proposition, we require two lemmas. The first, Lemma \ref{lem:NhatBound}, shows that $\hat{N}$ does not vary greatly on the fast, $O(k^2)$, time scale -- which is the timescale at which the idle server process evolves. The second, Lemma \ref{lem:NhatBound3} shows that there are no long excursion of the Idle server process; that is, of order greater that $O(k^{\epsilon/2})$.

In the following two lemmas and the proof of Proposition \ref{TightProp}, we let 
\[
\delta^{(k)} = k^{-1+\epsilon} \qquad \text{and}\qquad S^{(k)} = k^{\epsilon/2}.
\]

\begin{lemma}\label{lem:NhatBound} For positive constants $\gamma_1$ and $\beta_1$,
\begin{equation*}\label{NhatBound}
	\mathbb{P}
	\left(
		\sup_{m\in\mathbb{Z}_+} 
		\sup_{T_m^{(k)}\leq t \leq T_m^{(k)}+S^{(k)}} 
		\left|
			\hat{N}^{(k)}\left(\frac{t}{k^2}\right) - \hat{N}^{(k)} \left(\frac{T_m^{(k)}}{k^2}\right)
		\right|
		\geq 
		{\delta^{(k)}}  
	\right)
	\leq 
 \gamma_1 e^{-\beta_1 k^{\epsilon/2}}.
\end{equation*}	
\end{lemma}

\proof{Proof.}
Let $ 0= J_0^{(k)}< J_1^{(k)} < J_2^{(k)}  <... $ be the jump times of the process $(\hat{N}^{(k)}(t/k^2) : 0\leq t \leq k^2 \tau^{(k)}(\eta))$. Since each time $T_m^{(k)}$ is initiated by a jump time, it suffices to prove the result for jump times $J^{(k)}_m$. 

We prove the lemma in three steps. First, we prove a similar bound to \eqref{NhatBound} but for a single jump. Second, we prove that the number of jumps is of order $O(k^2)$. Third, we combine the two bounds using a union bound to give the result.

First, let $J^{(k)}=J_m^{(k)}$ be a jump of our process. Similar to expression \eqref{NExpression}, there exist independent unit rate Poisson Processes $\mathcal{N}_a$ and $\mathcal{N}_d$ (representing arrival and departure processes) such that, for $t\geq 0$
\begin{align}\label{NPois}
\hat{N}^{(k)}\left(\frac{t+J^{(k)}}{k^2}\right)- \hat{N}^{(k)}\left( \frac{J^{(k)}}{k^2} \right) =&  \frac{1}{k}{\mathcal N}_a \left( \left(1  -\frac{\alpha}{k} \right)\mu t\right)   - \frac{1}{k}{\mathcal N}_d \left( \mu t - \frac{\mu }{k^3} \int_0^\frac{t}{k^2}  \hat{I}^{(k)}(s)ds \right).
\end{align}
Thus, for $0\leq t\leq S^{(k)}$
\begin{equation}\label{OscB}
\sup_{0\leq t \leq S^{(k)}}\left|\hat{N}^{(k)}\left(\frac{t+J^{(k)}}{k^2}\right)- \hat{N}^{(k)}\left( \frac{J^{(k)}}{k^2} \right)\right| \leq  \sup_{0\leq t \leq S^{(k)}} \left[ \frac{1}{k} \left[\mathcal{N}_a (\mu t ) + \mathcal{N}_d (\mu t ) \right] \right] \leq \frac{1}{k} \left[\mathcal{N}_a (\mu 
S^{(k)} ) + \mathcal{N}_d (\mu S^{(k)} ) \right].
\end{equation}
Since $k\delta^{(k)} \geq e^2 2\mu S^{(k)}$ for all suitably large $k$, we can apply Lemma \ref{PoissonBound} to give that
\begin{equation}\label{Po2}
\mathbb{P} \left(\frac{1}{k} \left[\mathcal{N}_a (\mu 
S^{(k)} ) + \mathcal{N}_d (\mu S^{(k)} ) \right] \geq \delta^{(k)}  \right) \leq e^{- 2 \mu S^{(k)}-\delta^{(k)} k }.
\end{equation}
Inequalities \eqref{Po2} and \eqref{OscB} give the bound
\begin{equation}\label{NhatBound2}
\mathbb{P}\left(\sup_{J^{(k)}\leq t \leq J^{(k)}+S^{(k)}} \left|\hat{N}^{(k)}\left(\frac{t}{k^2}\right) - \hat{N}^{(k)} \left(\frac{J^{(k)}}{k^2}\right)\right| \geq {\delta^{(k)}}  \right)\leq e^{-2\mu S^{(k)} -{\delta^{(k)}} k}.
\end{equation}	

Second, we count the number of jumps that can occur in the time interval $[0,k^2 T]$. Again, since $(\hat{N}^{(k)}(t/k^2) : t\geq 0)$ has arrivals and departures that occur as a Poisson process of rate less than $2\mu$ (c.f. \eqref{NPois}), we have that from Lemma \ref{PoissonBound} that
\[
\mathbb{P}( |m : J^{(k)}_m\leq k^2T | \geq 2e^2\mu k^2 T  ) \leq \mathbb{P}\left( Po(2\mu k^2 T) \geq 2e^2\mu k^2 T \right) \leq e^{- (1+e^2) 2\mu k^2 T}.
\]

Third, we condition on whether the event $\{ |m : J^{(k)}_m\leq k^2T | \geq 2e^2\mu k^2 T\}$ occurs or not. When we condition on the event holding we apply  the above inequality, and when the event does not hold we apply a union bound to the finite number of transitions holding:
\begin{align*}
&\mathbb{P}\left(\sup_{m : J^{(k)}_m\leq k^2T} \sup_{0\leq t \leq S^{(k)}} \left|\hat{N}^{(k)}\left(\frac{t+J_m^{(k)}}{k^2}\right) - \hat{N}^{(k)} \left(\frac{J_m^{(k)}}{k^2}\right)\right| \geq {\delta^{(k)}}  \right)\\
\leq & e^{- (1+e^2) 2\mu k^2 T} + 2e^2 \mu k^2 T \sup_{m \leq 2e^2 \mu k^2 T}  	\mathbb{P}\left(\sup_{0\leq t \leq S^{(k)}} \left|\hat{N}^{(k)}\left(\frac{t+J_m^{(k)}}{k^2}\right) - \hat{N}^{(k)} \left(\frac{J_m^{(k)}}{k^2}\right)\right| \geq {\delta^{(k)}}  \right)\\
\leq & e^{- (1+e^2) 2\mu k^2 T} + 2e^2 \mu k^2 T e^{-2\mu S^{(k)} -{\delta^{(k)}} k}\\
\leq & \gamma_1 e^{-\beta_1 k^{\epsilon/2}}
\end{align*}
where the last inequality holds for appropriate constants $\gamma_1$, $\beta_1$.
This gives the required inequality.
\hfill \halmos\endproof 

Next, in the following lemma, we show that there are no long excursions of $\hat{I}^{(k)}$.

\begin{lemma} \label{lem:NhatBound3}
	\[
	\mathbb{P} 
		\left( 
			\exists m 
			\text{ s.t. } 
				T^{(k)}_{m+1} 
				- 
				T^{(k)}_{m}
				\geq 
				S^{(k)} 
		\right) 
	\leq  
	 \gamma_2 e^{-\beta_2 k^{\epsilon/2}}
	\]
\end{lemma}

\proof{Proof}
Our State Space Collapse Proposition, Proposition \ref{Thrm:StateSpaceCollapse}, states that
\[
\bP\Big( 
	\sup_{0\leq t\leq T} 
		\Big|
			\big(2-\hat{N}^{(k)}(t)\big)_+ 
			- 
			\hat{M}^{(k)}_1(t) 
		\Big| 
	\leq 
	k^{-1/2+\epsilon}
	\Big) 
\geq 
1-  
a \exp \left\{-bk^{\frac{\epsilon}{2}}\right\}.
\]
On the above event, the following sequence of bounds hold for times $t\leq \tau(\eta)$:
\[
\hat{M}_1^{(k)} (t)
\leq 
(2-\hat{N}^{(k)})_+
+
k^{-\frac{1}{2}+\epsilon}
\leq 
1-\eta
+
k^{-\frac{1}{2}+\epsilon}
\leq 
1- \frac{\eta}{2},
\]
where in the final inequality we assume that $k$ is sufficiently large  that $k^{-\frac{1}{2}+\epsilon} \leq \frac{\eta}{2}$. Thus we have that that
\[
\mathbb{P} 
	\left( 
		\hat{M}^{(k)}_1 (t) 
		\leq 
		1-\frac{\eta}{2},
		\quad \forall t \leq \tau^{(k)}(\eta) 
	\right) 
\geq 
1
-
a 
	\exp 
		\left\{
			-bk^{\frac{\epsilon}{2}}
		\right\}.
\]

Thus since 
 \begin{align*}
 &\hat{I}^{(k)}\mapsto \hat{I}^{(k)}+1 \qquad\text{at rate}\qquad \mu \hat{M}_1^{(k)},\\
 &\hat{I}^{(k)}\mapsto \hat{I}^{(k)}-1 \qquad\text{at rate}\qquad \mu \left( 1 -\frac{\alpha}{k}\right), \quad\text{ when }\quad \hat{I}^{(k)}>0,
 \end{align*}
the above inequality show that, with probability greater that $1-a\exp\left\{-bk^{\frac{\epsilon}{2}}\right\}$, we can bound all excursions of $\hat{I}^{(k)} (t)$ above by an $M/M/1$ queue with arrival rate $\mu (1 - \eta /2)$ and departure rate $\mu\left( 1- \frac{\alpha}{k}\right)$. By Lemma \ref{ExcusionLemma}, $\tilde{T}$, the length of an excursion for this $M/M/1$ queue can be bounded by
	\begin{equation}\label{TtildeBound}
	{\mathbb P} (\tilde{T} \geq S^{(k)} ) 
	\leq 
	\sqrt{\tfrac{1- \frac{\alpha}{k}}{ (1 - \eta /2)}}
	e^{
		-(\sqrt{\mu- \frac{\alpha \mu }{k}}- \sqrt{\mu (1 - \eta /2)})^2 
		S^{(k)}	
		}
	 \leq  \tfrac{1}{ \sqrt{(1 - \eta /2)}}e^{-\mu \frac{\eta^2}{16} S^{(k)}}
	\end{equation}
In the second inequality, above, we apply the bound that $1- \sqrt{1-x} \geq x/2$ for $x\geq 0$ as follows:
\[
(\sqrt{\mu- {\alpha \mu }/{k}}- \sqrt{\mu (1 - \eta /2)})^2 
\geq 
(\sqrt{\mu}- \sqrt{\mu (1 - \eta /2)})^2 = \mu ( 1 - \sqrt{1-{\eta}/2} )^{2}\geq \mu \frac{\eta^2}{16}.
\]

Each excursion of $\hat{I}^{(k)}$ is initiated by a departure and the number of departures is bounded above by a Poisson random variable of parameter $\mu k^2 T$. We know that the probability that this random viable exceeds $e^2 \mu k^2 T$ can be bounded as in Lemma \ref{PoissonBound} that is
\[
\mathbb{P} 
	\left( 
		|m : T^{(k)}_m\leq \tau(\eta) | 
		\geq e^2\mu k^2 T
	\right)
\leq
e^{- (1+e^2) \mu k^2 T} 
\]
and the thus applying \eqref{TtildeBound} and a union bound when $e^2 \mu k^2 T$ excursions occur gives the bound
	\begin{align*}
	&
	\mathbb{P} 
		\left( 
			\exists m 
			\text{ such that } 
			T^{(k)}_{m+1} - T^{(k)}_{m} \geq S^{(k)}  
		\right) 
	\\
	&
	\leq 
	\mathbb{P} 
	\left( 
	|m : T^{(k)}_m\leq \tau(\eta) | 
	\geq e^2\mu k^2 T
	\right)
	+
	e^2\mu k^2 T
		{\mathbb P} (\tilde{T} \geq S^{(k)} ) 
	\\
	&
	\leq
	e^{- (1+e^2) \mu k^2 T}  	
	+
	\tfrac{e^2 \mu k^2 T}{ \sqrt{(1 - \eta /2)}}
		e^{-\mu \frac{\eta^2}{4} S^{(k)}} 
	\end{align*}
	which is a bound of the required form.
\hfill \halmos
\endproof

We now combine Lemmas \ref{lem:NhatBound} and \ref{lem:NhatBound3} to prove Proposition \ref{TightProp}.

\proof{Proof of Proposition \ref{TightProp}.}
Applying a union bound to Lemma \ref{lem:NhatBound} and Lemma \ref{lem:NhatBound3} gives
	\begin{align*}
	&
	\mathbb{P}
		\left(
			\sup_{m\in\mathbb{Z}_+} 
			\sup_{T^{(k)}_{m}\leq t \leq T^{(k)}_{m+1}} 
			\left|
				\hat{N}^{(k)}\left(\frac{t+T^{(k)}_{m}}{k^2}\right) 
				- 
				\hat{N}^{(k)} \left(\frac{T^{(k)}_{m}}{k^2}\right)
			\right|
		 \geq {\delta^{(k)}}  
		 \right)\\
	\leq &
	\mathbb{P}
	\left(
	\sup_{m\in\mathbb{Z}_+} 
	\sup_{T_m^{(k)}\leq t \leq T_m^{(k)}+S^{(k)}} 
	\left|
	\hat{N}^{(k)}\left(\frac{t}{k^2}\right) - \hat{N}^{(k)} \left(\frac{J_m^{(k)}}{k^2}\right)
	\right|
	\geq 
	{\delta^{(k)}}  
	\right)
	+
		\mathbb{P} 
		\left( 
		\exists m 
		\text{ s.t. } 
		T^{(k)}_{m+1} 
		- 
		T^{(k)}_{m}
		\geq 
		S^{(k)} 
		\right) 
	\\
	\leq & \alpha_1 e^{-\beta_1 k^{\epsilon/2}} + \alpha_2 e^{-\beta_2 k^{\epsilon/2}}.
	\end{align*}	
By Lipschitz continuity of  the function $\hat{N}\mapsto (2-\hat{N})_+$, the same bound holds with $(2-\hat{N})_+$ in place of $\hat{N}$, namely,
	\begin{align*}
	&\mathbb{P}
		\left(
			\sup_{m \in \mathbb{Z}_+} 
			\sup_{T^{(k)}_{m}\leq t \leq T^{(k)}_{m+1}}
				\left|
					\left(2- \hat{N}^{(k)}\left(\frac{t+T^{(k)}_{m}}{k^2}\right) \right)_+ 
					- 
					\left(2-  \hat{N}^{(k)} \left(\frac{T^{(k)}_{m}}{k^2}\right) \right)_+
				\right| 
		\geq {\delta^{(k)}} 
		\right)\\
	\leq & \alpha_1 e^{-\beta_1 k^{\epsilon/2}} + \alpha_2 e^{-\beta_2 k^{\epsilon/2}} .
	\end{align*}	
	Finally applying a union bound along with state space collapse, Proposition \ref{Thrm:StateSpaceCollapse}, we have that 
		\begin{align*}
		&\mathbb{P}
			\left(
				\sup_{m \in \mathbb{Z}_+} 
				\sup_{T^{(k)}_{m}\leq t \leq T^{(k)}_{m+1}} 				
					\left|
						\hat{M}^{(k)}_1\left(\frac{t+T^{(k)}_{m}}{k^2}\right)  
						- 
						\left(
							2
							-  
							\hat{N}^{(k)} \left(\frac{T^{(k)}_{m}}{k^2}\right) 
						\right)_+
					\right| 
				\geq {\delta^{(k)}} + k^{-1/2+\epsilon} 
			\right)\\
		\leq & \alpha_1 e^{-\beta_1 k^{\epsilon/2}} + \alpha_2 e^{-\beta_2 k^{\epsilon/2}} + a e^{-b k^{\epsilon/2}} \\
		\leq &  \gamma e^{-\beta k^{\epsilon/2}}
		\end{align*}	
for appropriate positive constants $\gamma $ and $\beta$. \hfill\halmos
\endproof 

We can now see that the idleness process, $\hat{I}^{(k)} $, behaves approximately as an M/M/1 queue over each excursion. Thus in what follows we will mostly focus on behaviour of $\hat{I}^{(k)}$ on time intervals $[T_{m-1}^{(k)},T_{m}^{(k)}]$. Considering the statement of Proposition \ref{TimescaleIdleness}, the following Lemma shows that there is a not a great deal of error introduced by effectively rounding down to the nearest excursion time $T_m$.

\begin{lemma}\label{LemTm}
	\begin{equation*}
		\sup_{m} \sup_{T^{(k)}_m \leq t \leq T^{(k)}_{m+1} }
		\left|
		\frac{1}{k^2} 
			\int_{T^{(k)}_m}^{k^2t} \hat{I}^{(k)} (s) ds
			-
			\int_{T^{(k)}_m/k^2}^{t}	
					\frac{(2-\hat{N}^{(k)}(s))_+}{(\hat{N}^{(k)}(s)-1 )}
				ds 
		\right|
			\xrightarrow[k\rightarrow\infty]{} 0
\end{equation*}
where the above convergence is convergence in probability.
\end{lemma}
\proof{Proof}
By Lemma \ref{lem:NhatBound3}
	\[
	\mathbb{P} 
		\left( 
			\exists m 
			\text{ s.t. } 
				T^{(k)}_{m+1} 
				- 
				T^{(k)}_{m}
				\geq 
				k^{\epsilon/2}
		\right) 
	\leq  
	 \gamma_2 e^{-\beta_2 k^{\epsilon/2}}
	\]	
and by Lemma \ref{lem:IdleBound}i)	For $\epsilon>0$,
	\[
	\bP \left( \sup_{0\leq t\leq kT} \Ik(t) \geq k^{\frac{1}{2}+\epsilon} \right) \leq a k^2 e^{-b k^{\epsilon}} 
	\]
Further, by assumption
\begin{equation*}
	\frac{(2-\hat{N}^{(k)}(s))_+}{(\hat{N}^{(k)}(s)-1 )}\leq \frac{2}{\eta}
\end{equation*}
for $s\leq \tau(\eta)$. 

The three bounds above show that with high probability
\begin{align*}
\sup_{m} \sup_{T^{(k)}_m \leq t \leq T^{(k)}_{m+1} }
	\left|
		\frac{1}{k^2} 
			\int_{T^{(k)}_m}^{k^2t} \hat{I}^{(k)} (s) ds
			-
			\int_{T^{(k)}_m/k^2}^{t}	
					\frac{(2-\hat{N}^{(k)}(s))_+}{(\hat{N}^{(k)}(s)-1 )}
				ds 
	\right|
\leq 
	{k^{-3/2+3\epsilon/2}} 
+ 
	\frac{2}{\eta}
	{k^{-2+\epsilon/2}}
\end{align*}
as required.
	\hfill\Halmos
\endproof

\subsubsection{Bounding Processes}\label{BoundingI}
 Proposition \ref{TightProp} says that $\hat{M}_1^{(k)}$, which dictates the arrival of idle servers, can be approximated by $(2-\hat{N}(T_m^{(k)}/k^2))_+ \pm 2 k^{-1/2+\epsilon}$  where here $T_{m}^{(k)}$ is the value of the last excursion of $\hat{I}^{(k)}$.  
So we now have that $\hat{I}^{(k)}(t)$, for times $t\in [T_m^{(k)},T_{m+1}^{(k)})$ observes transition rates
 \begin{align*}
 &\hat{I}^{(k)}\mapsto \hat{I}^{(k)}+1 \qquad\text{at rate}\qquad \mu (2-\hat{N}(T_m^{(k)}/k^2))_+ +o(1),\\
 &\hat{I}^{(k)}\mapsto (\hat{I}^{(k)}-1)_+ \quad\text{at rate}\qquad \mu \left( 1 -\frac{\alpha}{k} \right).
 \end{align*}

Thus, as we will now prove, $(\hat{I}^{(k)}(t) : t \in [T_m^{(k)},T_{m+1}^{(k)}])$, can be bounded by processes ${I}_+^{(k)}$ and ${I}_-^{(k)}$. Here, for $t\in [T_m^{(k)},T_{m+1}^{(k)}]$, we define  ${I}_+^{(k)}(t)$ and ${I}_-^{(k)}(t)$ to be the continuous time Markov chain with ${I}_\pm^{(k)}(T_m^{(k)})=0$ and with non-zero transition rates
 \begin{align*}
 &{I}^{(k)}_{\pm}\mapsto {I}^{(k)}_{\pm}+1 
	 \qquad\text{at rate}
	\qquad \lambda^{(k)}_{m,\pm},\\
 &{I}^{(k)}_{\pm}\mapsto ({I}^{(k)}_{\pm}-1)_+ 
	 \quad\text{at rate}\qquad 
	 \mu^{(k)}_{m,\pm}.
 \end{align*}
For the $m$th idleness excursion, that is for times $T^{(k)}_m\leq t < T^{(k)}_{m+1}$, we can define the rates
\begin{align}\label{Lambdaplus}
\lambda^{(k)}_{m,+}& := \mu (2-\hat{N}^{(k)}(T^{(k)}_m))_++ 2\mu k^{-1/2+\epsilon}, \qquad \mu^{(k)}_{m,+} := \mu - \frac{\alpha}{k},\\
\lambda^{(k)}_{m,-} & := \mu (2-\hat{N}^{(k)}(T^{(k)}_m))_+ - 2\mu k^{-1/2+\epsilon}, \qquad \mu^{(k)}_{m,-} := \mu .
\end{align}
Note that both ${I}_+^{(k)}$ and ${I}_-^{(k)}$ are M/M/1 queues on time intervals $[T_m^{(k)},T_{m+1}^{(k)}]$.
	
The following lemma shows that $I^{(k)}_+$ and $I^{(k)}_-$ bound $\hat{I}^{(k)}$  above and below.
\begin{lemma}\label{IdleCoupling}
The processes $I^{(k)}_+$ and $I^{(k)}_-$ may be chosen so that
	\[
	\mathbb{P} \left( I^{(k)}_-(t) \leq \hat{I}^{(k)}(t) \leq I_+^{(k)}(t), \forall t\leq \tau^{(k)}(\eta) \right) \geq 1 -\gamma e^{-\beta k^{\epsilon/2}}
	\]
	 where $\gamma$ and $\beta$ are positive constants (as given in Proposition \ref{TightProp}).
\end{lemma}
\proof{Proof}
For each excursion $T^{(k)}_m\leq t < T^{(k)}_{m+1}$, we can write $\hat{I}$ in the form 
\[
\hat{I}^{(k)}(t) = {\mathcal N}_{m,d}\left( \mu \int_{T^{(k)}_m}^t \hat{M}^{(k)}_1(u) du\right)
-
\mathcal{N}_{m,a}\left( \int_{T^{(k)}_m}^t  \mu \Big( 1 -\frac{\alpha}{k} \Big)  ds \right)
\]
where ${\mathcal N}_{m,d}$, $\mathcal{N}_{m,a}$ are independent Poisson processes. By Proposition \ref{TightProp} with probability greater that $1 -\gamma e^{-\beta k^{\epsilon/2}}$, we have that
\[
 \int_{T^{(k)}_m}^t  \lambda_-^{(k)}(s)ds \leq \mu \int_{T^{(k)}_m}^t \hat{M}^{(k)}_1(u) du \leq  \int_{T^{(k)}_m}^t  \lambda_+^{(k)}(s)ds
\]
Thus 
\begin{align*}
I^{(k)}_-(t) = & {\mathcal N}_{m,d}\left( \mu \int_{T^{(k)}_m}^t  \lambda_-^{(k)}(s)ds \right)
 -
 \mathcal{N}_{m,a}\left(  \int_{T^{(k)}_m}^t  \lambda_+^{(k)}(s)ds \right)
\\
\leq & 
\hat{I}^{(k)}(t) 
\\
\leq &
 {\mathcal N}_{m,d}\left( \mu \int_{T^{(k)}_m}^t  \lambda_+^{(k)}(s)ds \right)
-
\mathcal{N}_{m,a}\left(  \int_{T^{(k)}_m}^t  \mu_+^{(k)}(s)ds \right)
=I^{(k)}_+(t)
\end{align*}
as required. \hfill \halmos

\endproof

We now work to prove a result similar to Proposition \ref{TimescaleIdleness} for the processes $I_{+}^{(k)}$ and $I_{-}^{(k)}$.

\begin{lemma}\label{IntNI}
\begin{equation*}
\sup_{0 \leq t \leq \tau^{(k)}(\eta)}
	\left|
	\int_{0}^{t}
	\bigg[
		\frac{
			\lambda_{\pm}(k^2 s)
			}{
			\mu_{\pm}(k^2 s)-\lambda_{\pm}(k^2 s)}	 
		- \
		\frac{
			\big( 2- \hat{N}^{(k)}(s)\big)_{+} 
			}{ 
			\hat{N}^{(k)}(s) - 1
			}	
	\bigg]
	ds 
	\right| 
\xrightarrow[k\rightarrow\infty]{} 0
\end{equation*}
where here convergence is convergence in probability.
\end{lemma}
\proof{Proof}
Recalling the definitions of $\lambda_{\pm}(s)$ we have for instance, that,
\begin{align}
\frac{
	\lambda^{(k)}_{-}(k^2s)
	}{
	\mu_{-}(k^2 s)-\lambda^{(k)}_{-}(k^2s)
	}  
:= 
\frac{
	(2-\hat{N}^{(k)}(T^{(k)}_m))_+ - 2 k^{-1/2+\epsilon}
	}{
	\hat{N}^{(k)}(T^{(k)}_m)-1+2 k^{-1/2+\epsilon}.
}
\end{align}	
Now since $n\mapsto \frac{(2-n)_+}{n-1}$ is Lipschitz continuous on $n\geq 1+\eta/2$ with Lipschitz constant $2/\eta$. We can apply 
by Lemma \ref{lem:NhatBound} (with $S_k$ chosen so that Lemma \ref{lem:NhatBound3} holds):
\begin{align*}
&
\bigg|
		\frac{
			\lambda_{-}(k^2 s)
		}{
		\mu_{-}(k^2 s)-\lambda_{-}(k^2 s)}	 
	- \
	\frac{
		\big( 2- \hat{N}^{(k)}(s)\big)_{+} 
	}{ 
	\hat{N}^{(k)}(s) - 1
}	
\bigg|
\\
\leq &
\frac{2}{\eta}	\sup_{m\in\mathbb{Z}_+} 
	\sup_{T_m^{(k)}\leq t \leq T_{m+1}^{(k)}} 
	\left|
	\hat{N}^{(k)}\left(\frac{t}{k^2}\right) - \hat{N}^{(k)} \left(\frac{T_m^{(k)}}{k^2}\right)
	\right| + \frac{4}{\eta} k^{-1/2 + \epsilon}
\\	
\leq & 
	\frac{2}{\eta} \delta^{(k)}+ \frac{4}{\eta} k^{-1/2+\epsilon}.
	\end{align*}
	Applying this bound to the integral stated in the lemma gives the result. (The same bound for $\lambda_+$ and $\mu_+$ holds my a more-or-less identical argument.)
	\hfill\Halmos
\endproof

\begin{lemma}\label{IntIPlus}
Almost surely,
	\begin{align}\label{I+Integral}
	\limsup_{k\rightarrow\infty}
	\sup_{m \in \mathbb{Z}_{+}} 
			\frac{1}{k^2}  
			\int_0^{T_m^{(k)}}	
			\left[ 
				I_+^{(k)}(s) 
				-
				\frac{
					\lambda^{(k)}_+(s)
					}{
					\mu_+^{(k)}(s)-\lambda_+^{(k)}(s)
					}
			\right]
			ds 
=
&
0
\\
	\liminf_{k\rightarrow\infty}
	\inf_{m \in \mathbb{Z}_{+}} 
			\frac{1}{k^2}  
			\int_0^{T_m^{(k)}}	
			\left[ 
				I_-^{(k)}(s) 
				-
				\frac{
					\lambda^{(k)}_-(s)
					}{
					\mu_-^{(k)}(s)-\lambda_-^{(k)}(s)
					}
			\right]
			ds 
=
&
0
	\end{align}
\end{lemma}
\proof{Proof}
We prove the result for $I_+^{(k)}$. 
The result holds for $I_-^{(k)}$ by a more-or-less identical argument.
Let $T^+_{m+1}$ be the time after $T_m^{(k)}$ that $I^{(k)}_+(s)$ next empties. 
(Here we assume that the process $I^{(k)}_+(s)$ is allowed to continue past time $T_{m+1}^{(k)}$ to $T_{m+1}^{+}$ time with transition rates \eqref{Lambdaplus}.)
 Further, let 
\[
Z_m = 
	\int_{T_m}^{T_{m+1}^+} 
		I_+^{(k)}(s) 
		- 
		\frac{
			\lambda^{(k)}_+(s)
			}{
			\mu_+^{(k)}(s)-\lambda_+^{(k)}(s)
			} 
	ds.
\]
We can now expand the integral in \eqref{I+Integral} as follows
\begin{align}
&\frac{1}{k^2}  
	\int_0^{T_m^{(k)}}	
	\left[ 
		I_+^{(k)}(s) 
		-
		\frac{
			\lambda^{(k)}_+(s)
		}{
			\mu_+^{(k)}(s)-\lambda_+^{(k)}(s)
		}
	\right]
	ds \notag
	\\
= 
&
\frac{1}{k^2}
	\sum_{m'\leq m}
		\int_{T^{(k)}_{m'-1}}^{T^{(k)}_{m'}}
			I_+^{(k)} (s)
			-
			\frac{
				\lambda^{(k)}_+(s)
				}{
				\mu_+^{(k)}(s)-\lambda_+^{(k)}(s)
				}
		ds \notag
\\
=
&
\frac{1}{k^2} 
	\sum_{m'\leq m}
		Z_{m'}
-
\frac{1}{k^2}
	\sum_{m'\leq m}
	\int_{T_{m'}^{(k)}}^{T_{m'}^+}
	\bigg[
		I_+^{(k)}(s)
		-
		\frac{
			\lambda^{(k)}_+(s)
		}{
			\mu_+^{(k)}(s)-\lambda_+^{(k)}(s)
		}
		\bigg]	
	ds 
\\
\leq 
& 
\frac{1}{k^2} 
	\sum_{m'\leq m}
		Z_{m'}
+
\frac{1}{k^2}
	\sum_{m'\leq m}
	\int_{T_{m'}^{(k)}}^{T_{m'}^+}
		\frac{
			\lambda^{(k)}_+(s)
		}{
			\mu_+^{(k)}(s)-\lambda_+^{(k)}(s)
		}
	ds 	\label{ZExpand}
\end{align}
We now separately bound the summation over $Z$ and integral in \eqref{ZExpand}.

Firstly, for the summation over $Z$. 
The term $Z_m$ is the difference between the area under the excursion of an $M/M/1$ queue and its mean, see Lemma \ref{ExcusionLemma4}.
Thus, by Lemma \ref{ExcusionLemma4}, 
\[
\mathbb{E} \left[ Z_m \Big| \mathcal{F}_{T^{(k)}_m} \right] =0
\]
for all $m\in\mathbb{Z}_+$. 
In other words, the sum over $Z_m$ is a Martingale difference sequence. 
By Lemma \ref{ExcusionLemma}, $Z_m$ is of finite variance. So, by Lemma \ref{MgSLLN},
\[
\sup_{m} \left| \frac{1}{k^2}\sum_{m'\leq m} Z_{m'}\right| \xrightarrow[k\rightarrow\infty]{ a.s.}0
\]

Now consider the second integral in \eqref{ZExpand}. Since we chose times $s\leq k^2 \tau^{(k)}(\eta)$, we have that 
\[
0\leq \frac{
	\lambda^{(k)}_+(s)
}{
	\mu_+^{(k)}(s)-\lambda_+^{(k)}(s)
}
\leq 
\frac{6\mu}{\eta}.
\]
Further we have that
\begin{align*}
	\frac{1}{k^2}
	\sum_{m'\leq m} 
		T^+_{m'} - T^{(k)}_{m'} 
&\leq 
	\frac{1}{k^2}
	\sum_{m'\leq m} 
		T^+_{m'} - T^{-}_{m'} 
\\
&=
	\frac{1}{k^2}
	\sum_{m'\leq m} 
		T^+_{m'} - \mathbb{E} T^+_{m'} 	
		- 
		\frac{1}{k^2}
	\sum_{m'\leq m} 
		T^{-}_{m'} -\mathbb{E} T^{-}_{m'} 
		+
				\frac{1}{k^2}
	\sum_{m'\leq m} 
		\mathbb{E} T^{+}_{m'} -\mathbb{E} T^{-}_{m'} \\
&
\xrightarrow[k\rightarrow\infty]{ a.s.}0 
\end{align*}
The first two terms above converge to zero, by Lemma \ref{MgSLLN}, since again they are a martingale difference sequence. The expectations in the final term can be calculated explicitly, in particular,
\begin{equation*}
		\mathbb{E} T^{+}_m =
		\frac{1}{1-\lambda_{m,+}^{(k)}/\mu_{m,+}^{(k)}} - \frac{1}{\lambda_{m,+}^{(k)}},
		\qquad 
		\mathbb{E} T^{-}_m =
		\frac{1}{1-\lambda_{m,-}^{(k)}/\mu_{m,-}^{(k)}} - \frac{1}{\lambda_{m,-}^{(k)}}
\end{equation*}
and thus $\mathbb{E} T^{+}_m - \mathbb{E} T^{+}_m $ convergences to zero since $\lambda_{m,+}^{(k)}$ and $\lambda_{m,-}^{(k)}$ (and $\mu_{m,+}^{(k)}$ and $\mu_{m,-}^{(k)}$) converge to the same value. 
Thus the 2nd integral in \eqref{ZExpand} also convergences to zero:
\begin{equation}
0\leq \frac{1}{k^2}
	\sum_{m'\leq m}
	\int_{T_{m'}^{(k)}}^{T_{m'}^+}
		\frac{
			\lambda^{(k)}_+(s)
		}{
			\mu_+^{(k)}(s)-\lambda_+^{(k)}(s)
		}
	ds 
\leq 
\frac{6\mu}{\eta}	\frac{1}{k^2}
	\sum_{m'\leq m} 
\left(	
	T^+_{m'} - T^{(k)}_{m'} 
\right) 
\xrightarrow[k\rightarrow\infty]{ a.s.}0 .
\end{equation}
Thus since each term in \eqref{ZExpand} converges to zero we have the required result.
\hfill \halmos
 
\endproof

\subsubsection{Proof of Proposition \ref{TimescaleIdleness}}\label{sec:Prop3Proof}
We are now in a position to prove Proposition \ref{TimescaleIdleness}. 
With Lemma \ref{IdleCoupling}, we have found a suitable processes $I_-^{(k)}$ and $I_+^{(k)}$ to bound to the idleness processes $\hat{I}^{(k)}$. Further, in Lemmas \ref{IntNI} and \ref{IntIPlus}, we see that $I_-^{(k)}$ and $I_+^{(k)}$ obey a similar limit result to that required by Proposition  \ref{TimescaleIdleness}. We now put these together to prove the proposition.

\proof{Proof of Proposition \ref{TimescaleIdleness}.}
By Lemma \ref{IdleCoupling}, with probability greater than $1-\gamma e^{-\beta k^{\epsilon/2}}$ we have that
\begin{subequations}\label{Expanded0}
\begin{align}
&
\frac{1}{k^2}
	\int_{0}^{T^{(k)}_m} 
	\left[ 
		I_-^{(k)} (s) 
		-
		\frac{\lambda_-(s)
		}{
		\mu_-(s) - \lambda_-(s)
		} 
	\right]
		ds
+
	\int_{0}^{T^{(k)}_m/k^2}
	\Bigg[ 
		\frac{\lambda_-(k^2s)
		}{
		\mu_-(k^2s) - \lambda_-(k^2s) 
		}		
		-
		\frac{\big( 2- \hat{N}^{(k)}(s)\big)_+}{\hat{N}^{(k)}(s)-1} 
	\Bigg]
	ds 
\label{Expanded1}
\\
\leq 
&
	\frac{1}{k^2} 
		\int_{0}^{T_{m}^{(k)}} 
			\hat{I}^{(k)}(s) 
		ds 
		-
		\int_{0}^{ T_{m}^{(k)}/k^2} 
			\frac{\big( 2- \hat{N}^{(k)}(s)\big)_+}{\hat{N}^{(k)}(s)-1} 
		ds 
\notag
\\
\leq
&
\frac{1}{k^2}
	\int_{0}^{T^{(k)}_m} 
		\left[ 
			I_+^{(k)} (s) 
			-
			\frac{\lambda_+(s)
				}{
				\mu_+(s) - \lambda_+(s)
			} 
		\right]		
		ds
+
	\int_{0}^{T^{(k)}_m/k^2}
	\Bigg[ 
		\frac{
			\lambda_+(k^2 s)
			}{
			\mu_+(k^2 s) - \lambda_+(k^2 s)
			}		
	-
	\frac{\big( 2- \hat{N}^{(k)}(s)\big)_+}{\hat{N}^{(k)}(s)-1} 
\Bigg]
	ds.
	\label{Expanded3} 
\end{align}
\end{subequations}
The first integral in \eqref{Expanded1} and \eqref{Expanded3} converges in probability to zero by Lemma \ref{IntIPlus}. The second integral in \eqref{Expanded1} and \eqref{Expanded3} converges in probability by Lemma \ref{IntNI}. Thus 
\begin{equation*}
	\sup_{m} 
	\left|
	\frac{1}{k^2} 
	\int_{0}^{T^{(k)}_m} \hat{I}^{(k)} (s) ds
		-
		\int_{0}^{T^{(k)}_m/k^2}	
			\frac{(2-\hat{N}^{(k)}(s))_+}{(\hat{N}^{(k)}(s)-1 )}
				ds 
	\right|		\xrightarrow[k\rightarrow\infty]{} 0
\end{equation*}
where the above convergence is convergence in probability. Combining this with Lemma \ref{LemTm} gives Proposition \ref{TimescaleIdleness}.
\hfill\Halmos 
\endproof

\subsection{Properties of Stochastic Differential Equations for $\hat{N}$}
Thus far we have focused on the weak convergence of the stochastic processes $\hat{N}^{(k)}$. 
In this subsection we focus on the properties of the limit stochastic differential equation:
\begin{equation*}
	d \hat{N} (t) = \mu  \Big[ \frac{(2-\hat{N})_+}{\hat{N}-1}-\alpha\Big] dt + \sqrt{2\mu} \;dB (t), \qquad t\geq 0.
\end{equation*}
In Lemma \ref{CTS_SDE} we prove continuity of solutions of the above SDE with respect to data given by the driving process $(B(t): t \geq 0)$.
In Lemma \ref{DiffZero} we prove path-wise uniqueness of solutions.
In Proposition \ref{JSQ:EqPro} we characterize the stationary distribution of this SDE.  

\begin{lemma}\label{CTS_SDE}
	Consider the integral form
	\begin{equation}\label{IntForm}
	n(t) = n(0) + \sqrt{2\mu} b(t) + \mu \int_0^t \left[  \frac{(2-n(s))_+}{(n(s)-1)} - \alpha \right] dt.
	\end{equation}
	Let $b(t)$ be a continuous function and $n(t)$ a solution to the integral expression, above, with $n(t)$ strictly greater than 1. 
%
If we let $\tau^{(k)}(\eta)$ be the first time that $n^{(k)}(s) \leq 1+ \eta $ and $b^{(k)}$ converges to $b$ uniformly on compacts in $[0,\tau^{(k)}(\eta)]$, i.e. for all $\eta>0$ and $T>0$
	\[
	\sup_{0\leq t \leq \tau^{(k)}(\eta)\wedge T} |b^{(k)}(t) - b(t)| \xrightarrow[k\rightarrow\infty]{} 0
	\]
	then $n^{(k)}$ approaches $n$ uniformly on compacts, that is,
\begin{equation*}
		\sup_{0\leq t \leq T} 
|n^{(k)}(t) - n(t)| \xrightarrow[k\rightarrow\infty]{} 0
\end{equation*}

\end{lemma}

\proof{Proof}
Suppose $b(t)$ is continuous and such that $n(t)>0$. 
Let $\eta$ be the minimum value of $n(t)-1$ for $t\leq T$.
Note that on this domain the function $n\mapsto (2-n)_+/(n-1)$ is Lipschitz continuous with Lipschitz constant $\eta^{-1}$. 
Let $\tau^{(k)}=\tau^{(k)}(\eta/2)$ be the first time that $n^{(k)}\leq 1+ \eta/2$ occurs.
By assumption for each $\epsilon>0$, there exists a $K$ such the for all $k\geq K$ we have that 
\[
|n^{(k)}(0) - n(0) | 
+
\sqrt{2\mu}\sup_{t\leq T} \left| b^{(k)}(t) - b(t) \right| \leq \epsilon .
\]
Further we choose $\epsilon>0$ such that 
\begin{equation} \label{epsiloncond}
\epsilon < \frac{\eta}{2} e^{-\eta^{-1}T}
\end{equation}

With the triangle-inequality, we can bound the process for times $t\leq \tau^{(k)}\wedge T$ as follows
\[
| n^{(k)}(t) - n(t)| 
\leq 
| n^{(k)}(0) - n(0) | 
+ \sqrt{2\mu} |b^{(k)}(t)-b(t)| 
+  \mu \int_0^t 
\left|
\frac{(2-n^{(k)}(s))_+}{(n^{(k)}(s)-1)} - \frac{(2-n(s))_+}{(n(t)-1)}
\right| ds
\]
Thus applying the Lipschitz condition and maximizing over $t$ we have that
\begin{align*}
&\sup_{t\leq \tau^{(k)} \wedge T} 
| n^{(k)}(t) - n(t)| \\
&\leq 
| n^{(k)}(0) - n(0) | 
+ \sqrt{2\mu} \sup_{t\leq \tau^{(k)} \wedge T} |b^{(k)}(t)-b(t)| 
+  \mu \int_0^{\tau^{(k)} \wedge T}
\eta^{-1} \sup_{u\leq s} 
| n^{(k)}(u) - n(u)| ds .
\end{align*}
Applying Gronwall's Lemma, we have that
\begin{align}
\sup_{t\leq \tau^{(k)} \wedge T} 
| n^{(k)}(t) - n(t)|  
&\leq   
\left[ | n^{(k)}(0) - n(0) | 
+ \sqrt{2\mu} \sup_{t\leq T\wedge \tau^{(k)}} |b^{(k)}(t)-b(t)| \right] e^{\eta^{-1} T} \\
&\leq \epsilon e^{\eta^{-1} T} \label{epsilonBound} \\
&< \frac{\eta}{2} \label{HitBound}
\end{align}
where the second inequality above holds by the definition of $\epsilon$ and the final inequality holds by condition \eqref{epsiloncond}. By Inequality \eqref{HitBound} we see that for our choice of $\epsilon$ that
\[
(n^{(k)}(t) -1) > (n(t) -1)  -\frac{\eta}{2} \geq \frac{\eta}{2} 
\]
for all times $t\leq \tau^{(k)}\wedge T$. Thus $\tau^{(k)}>T$, and therefore, Inequality \eqref{epsilonBound} now states 
\begin{equation}
\sup_{t\leq T} | n^{(k)}(t) - n(t)|
\leq   \epsilon e^{\eta^{-1} T}
\end{equation}
In other words, we see that as required
\[
\lim_{k\rightarrow\infty}\sup_{t\leq T} | n^{(k)}(t) - n(t)| = 0.
\]
\hfill \Halmos
\endproof

The following proposition summarizes the properties of SDE \eqref{JSQ:Diffusion}.
\begin{lemma}\label{DiffZero}
	A solution to the SDE 
	\begin{equation}\label{SDE1}
	d \hat{N} (t) = \mu  \Big[ \frac{(2-\hat{N})_+}{\hat{N}-1}-\alpha\Big] dt + \sqrt{2\mu} \;dB (t), \qquad t\geq 0
	\end{equation}
	exists, and is path-wise unique. Further the explosion time is almost surely not finite.
\end{lemma}
\proof{Proof}
We can prove the result with a localization argument. In particular, for $\epsilon \in (0,\frac{1}{2})$, we consider the stopping times 
$\tau_\epsilon= \inf\{t : \hat{N}(t) = 1+\epsilon \}$ and define $\tau_1 = \lim_{\epsilon\rightarrow 0} \tau_{\epsilon}$ to be the explosion time of the SDE \eqref{SDE1}. The processes $\hat{N}(t)$, for $t\leq \tau_\epsilon$ satisfies an SDE with Lipschitz coefficients. Such SDEs are known to be path-wise unique, for instance, see \cite{rogers2000diffusions}, Theorem V.11.2. So solutions to the SDE, \eqref{SDE1}, are path-wise unique until time $\tau_\epsilon$.
So path-wise uniqueness holds on each time interval $[0,\tau_\epsilon)$.

We will now show that, almost surely, $\tau_1=\infty $, and thus, path-wise uniqueness holds for all time.
Notice that without loss of generality, we may assume that the initial condition satisfies $1<\hat{N}(0)<2$ (because if the limit $\lim_{\epsilon\rightarrow 0} \tau_\epsilon $ were finite then there exists a stopping time such that $\hat{N}(\tau_\epsilon ) < 2$ and invoking the strong Markov property we can start from this time instead -- Note, local Lipschitz continuity of coefficient ensures the process is Strong Markov up to its time of explotion.) Further, let $\tau_2 = \inf\{t : \hat{N}(t) = 2 \}$. Notice that, by a similar argument, if $\mathbb P ( \tau_1 < \infty ) >0$ then $\mathbb P ( \tau_1 < \tau_2 ) >0$. (If explosion occurs then, every time there must be a positive probability that explosion occurs before $\hat{N}$ hits to $2$.)

For times $t < \tau_1 \wedge \tau_2$ the process obeys the SDE
\begin{equation}
d \hat{N} (t) = \mu  \Big[ \frac{1}{\hat{N}-1}-(\alpha+1)\Big] dt + \sqrt{2\mu} \;dB (t), \qquad t\geq 0.
\end{equation}
Notice that aside from the drift of $-(\alpha+1)$ the above processes is similar to a Bessel process. (In particular, if we remove the drift term $-(\alpha+1)$ then $(\hat{N}-1)/\sqrt{2\mu} $ would be is a Bessel process of Dimension 2.) By Girsanov's Theorem we can add and remove this drift, with an exponential change of measure. Thus on each compact time interval, a zero probability event with drift is a zero probability event without a drift. It is well known that a $2$ dimensional Bessel process is non-explosive, since two dimensional Brownian motion is neighbourhood recurrent at zero, see \cite{rogers2000diffusions} Section VI.35 for a proof. Thus $\mathbb P ( \tau_1 < \tau_2 ) = 0$. Thus the result now follows from the existence and uniqueness of SDEs with Lipschitz coefficients. 
	(Further, an alternative proof is provided by the Watanabe-Yamada path-wise uniqueness theorem; see \cite{rogers2000diffusions} Theorem V.40.1.)
	\hfill \Halmos
\endproof 

\begin{proposition}\label{JSQ:EqPro}
	The SDE 
	\begin{equation}
	d \hat{N} (t) = \mu  \Big[ \frac{(2-\hat{N})_+}{\hat{N}-1}-\alpha\Big] dt + \sqrt{2\mu} \;dB (t), \qquad t\geq 0
	\end{equation}
has a stationary distribution with probability density function 
	\begin{equation}\label{JSQ:Equilibrium2}
	\frac{d\pi}{dn} \propto 
	\begin{cases}
	\exp \left\{ -\alpha (n-2) \right\}, & n \geq 2, \\
	(n-1)\exp \left\{ -(\alpha+1)(n-2) \right\}, & 1\leq n\leq 2.
	\end{cases}
	\end{equation}
and this distribution has expected value:
	\begin{align} \label{ENHat}
	\bE [\hat{N}] 
=	1+ \left[ \frac{1}{\alpha (1+\alpha)^2} + \frac{e^{1+\alpha}}{(1+\alpha)^2} \right]^{-1}\left[  \frac{\alpha^2+4\alpha+1}{\alpha^2 (1+\alpha)^3} + \frac{2 e^{\alpha+1}}{(1+\alpha)^3} \right]
	\end{align}
\end{proposition}
\proof{Proof}
Applying the substitution $X= \frac{1}{\sqrt{\mu}} \hat{N}$, we have that 
\begin{align}
d X =& \sqrt{2} d B_t + \sqrt{\mu} \left[ \frac{(2-\sqrt{\mu}X)_+}{(\sqrt{\mu}X-1)} - \alpha \right] dt \notag\\
=& \sqrt{2} d B_t - V'(X) dt \label{DiffForm}
\end{align}
where $V(x)$ is a function with derivative 
\[
V'(x) =\sqrt{\mu} \left[ \alpha - \frac{(2-\sqrt{\mu}x)_+}{(\sqrt{\mu}x-1)}   \right]
\]
The diffusion \eqref{DiffForm} is a Langevin Diffusion. It is well known that such diffusions have an invariant measure given by
\[
\tilde{\pi}(x)= \exp\{-V(x) \}.
\]
This can be verified directly as $\pi$ is a stationary solution to the Fokker-Plank equations. (A more exacting analysis giving exponential convergence of Langevin diffusions to stationarity is given in \cite{roberts1996}.)

What remains are somewhat tedious calculations resulting in the required expression \eqref{JSQ:Equilibrium2} and \eqref{ENHat}. For \eqref{JSQ:Equilibrium2}, integrating the expression for $V'(x)$ we get that 
\begin{align*}
V(x) 
&= \sqrt{\mu} \int_{\frac{2}{\sqrt{\mu}}}^{x} \alpha - \frac{(\frac{2}{\sqrt{\mu}}-x)_+}{(x-\frac{1}{\sqrt{\mu}})}dx\\
&= \sqrt{\mu} \alpha \left(x-\frac{2}{\sqrt{\mu}}\right) + \sqrt{\mu}\mathbb{I}\left[1<\sqrt{\mu}X<2\right] \int_{x}^{\frac{2}{\sqrt{\mu}}} \frac{\frac{2}{\sqrt{\mu}}-x}{x-\frac{1}{\sqrt{\mu}}}dx\\
&=
\sqrt{\mu} \alpha \left(x-\frac{2}{\sqrt{\mu}}\right) + \sqrt{\mu}\mathbb{I}\left[1<\sqrt{\mu}X<2\right] \int_{x}^{\frac{2}{\sqrt{\mu}}} \left[ \frac{\frac{1}{\sqrt{\mu}}}{x-\frac{1}{\sqrt{\mu}}} - 1 \right] dx\\
&=\sqrt{\mu} \alpha \left(x-\frac{2}{\sqrt{\mu}}\right) + \sqrt{\mu}\mathbb{I}\left[1<\sqrt{\mu}X<2\right]
\left[ 
\frac{1}{\sqrt{\mu}} \log \left(\frac{1}{\sqrt{\mu}}\right) - \frac{1}{\sqrt{\mu}} \log \left(x-\frac{1}{\sqrt{\mu}}\right) + \left(x-\frac{1}{\sqrt{\mu}}\right)
\right]
\end{align*}
Therefore we arrive at an invariant distribution $\tilde{\pi}$ with density
\[
\frac{d \tilde{\pi}}{d x}= \exp\{-V(x) \} = 
\begin{cases}
\exp \left\{-\sqrt{\mu}\alpha \left( x- \frac{2}{\sqrt{\mu}} \right)\right\} & \text{if } x \geq \frac{2}{\sqrt{\mu}} \\
\left( \sqrt{\mu} x -1 \right)\exp \left\{ -\sqrt{\mu}(\alpha +1) \left( x -\frac{2}{\sqrt{\mu}}\right)  \right\} & \text{if } x < \frac{2}{\mu}.
\end{cases}
\]

So substituting back $X= \frac{1}{\sqrt{\mu}} \hat{N}$, we see the stationary distribution for $\hat{N}$, $\pi(n)$, is given by.
\[
\frac{d \pi }{d n} =
\begin{cases}
C\exp\{-\alpha (n-2)\}, & \text{if } n \geq 2, \\
C (n-1)\exp \{- (\alpha +1) (n-2) \}, & \text{if } 1\leq n \leq 2.
\end{cases}
\]
where 
\[
C= \int_0^2 (n-1)\exp \{- (\alpha +1) (n-2) \} dn + \int_2^\infty \exp\{-\alpha (n-2)\} dn
\] 
is a normalizing constant that we will now calculate.

To calculate the normalizing constant $C$ and the expectation \eqref{ENHat}. We make use of the following standard (Gamma function) identity
\[
\int_{0}^{\infty} n^s e^{-\theta n}dn= \frac{s!}{\theta^{s+1}}.
\] 

First, we calculate the two terms in the expression for $C$, above. Namely
\begin{align*}
&\int_1^2 (n-1) e^{- (\alpha+1)(n-2) } dn \\
= &\int_{1}^{\infty} (n-1) e^{-(\alpha+1)(n-1)} \cdot e^{-(\alpha+1)}dn- \int_2^{\infty} (n-2) e^{-(\alpha+1)(n-2)}dn - \int_2^\infty e^{- (\alpha +1) (n-2)}dn \\
=& e^{\alpha+1} \frac{1}{(\alpha+1)^2} - \frac{1}{(\alpha+1)^2} - \frac{1}{(\alpha+1)}
\end{align*}
and
\begin{align*}
\int_2^\infty e^{-\alpha(n-2)} dn = \frac{1}{\alpha}. 
\end{align*}
Thus
\begin{equation}\label{Ccalc}
	C= \frac{e^{(\alpha+1)}-1}{(\alpha+1)^2} + \frac{1}{\alpha} - \frac{1}{\alpha+1}.
\end{equation}

Second, we calculate the expected value of $\hat{N}$. 
\begin{align*}
C\mathbb{E} \hat{N} & = \int_1^2 n(n-1) e^{-(\alpha+1)(n-2)} dn + \int_2^\infty n e^{-\alpha(n-2)}dn \\
&=\int_1^\infty \underbrace{n(n-1)}_{n(n-1)=(n-1)^2 + (n-1)} e^{-(\alpha+1)(n-2)} dn \\
&\quad - \int_2^\infty \underbrace{n(n-1)}_{n(n-1)=(n-2)^2 + 3(n-2)+2} e^{-(\alpha+1)(n-2)} dn \\
&\quad\quad + \int_2^\infty \underbrace{n}_{n=(n-2)+2} e^{-\alpha(n-2)}dn\\
&= e^{(\alpha+1)} \left[ \int_0^\infty n^2 e^{-(\alpha+1)n}dn + \int_0^\infty n e^{-(\alpha+1) n}dn \right]\\
&\quad - \left[ \int_0^\infty n^2 e^{-(\alpha+1)n } dn + 3 \int_0^\infty n e^{-(\alpha+1)n}dn + 2 \int_0^\infty e^{-(\alpha+1)n}dn\right]\\
& \quad\quad + \left[ \int_0^\infty ne^{-\alpha n}dn + 2\int_0^\infty e^{-\alpha n}dn \right] \\
& = e^{(\alpha+1)} \left[\frac{2}{(\alpha+1)^3} + \frac{1}{(\alpha+1)^2}\right] \\
 &\quad - \left[ \frac{2}{(\alpha+1)^3} + \frac{3}{(\alpha+1)^2} + \frac{2}{(\alpha+1)}\right]\\
 &\quad\quad + \left[ \frac{1}{\alpha^2} + \frac{2}{\alpha} \right].
\end{align*}
This and the value of $C$ calculated above gives the stationary expectation
\begin{equation*}
\mathbb{E} \hat{N} =	\frac{
		e^{\alpha+1} 
		\Big[
		\frac{2}{(\alpha+1)^3}
		+
		\frac{1}{(\alpha+1)^2}
		\Big]
		-
		\Big[
		\frac{2}{(\alpha+1)^3}
		+
		\frac{3}{(\alpha+1)^2}
		+
		\frac{2}{(\alpha+1)}
		\Big]		
		+	
		\Big[
		\frac{1}{\alpha^2}
		+
		\frac{2}{\alpha}
		\Big]
	}
	{
		\frac{e^{(\alpha+1)}-1}{(\alpha+1)^2}
		+
		\frac{1}{\alpha(1+\alpha)}
	}.
\end{equation*}
The formula \eqref{ENHat} is a somewhat simplified form of this expression.
\hfill \Halmos 
\endproof 

\begin{lemma}\label{StationaryProp}
The Idle-Queue-First NDS approximation given by
\begin{equation*}
	d\hat{N}_{IQF} (t)
= 
\mu \left[
		\frac{1}{\hat{N}_{IQF}-1}
		-
		\alpha 
	\right]  dt
+
\sqrt{2\mu}dB(t),\qquad t\geq 0
\end{equation*}
has a stationary expected queue size of
\begin{equation*}
	\mathbb{E} [ \hat{N}_{IQF} ] 
	=
	1+ \frac{2}{\alpha}.
\end{equation*}
\end{lemma}
\proof{Proof}
Note that $\hat{N}_{IQF} $ is a Bessel process with drift $-\alpha \mu t$. Following a more-or-less identical steps from \eqref{DiffForm} through to \eqref{Ccalc}. We see that this process has stationary distribution
$({\pi}_{IQF}(n): n \geq 1)$ with density:
\begin{equation*}
	\frac{d\pi_{IQF}}{dn}
	= 	
\alpha^{2}(n-1)\exp \{-\alpha(n-1) \}.
\end{equation*}
From this it is straightforward calculation to show that
\[
\mathbb{E} [ \hat{N}_{IQF} ]  = 1+ \int_1^\infty \alpha^2 (n-1)^2 e^{-\alpha (n-1)} dn = 1 + \frac{2}{\alpha}.
\]
\hfill\Halmos 
\endproof

Notice from Lemma \ref{StationaryProp} and Central Queue expected queue size \eqref{CQED} that 
\begin{equation}
 \sup_{\alpha} \frac{\mathbb{E}\big[{\hat{N}_{IQF}}\big]}{\mathbb{E}{\big[\hat{N}_{CQ}\big]}} 
=
\sup_{\alpha}
\left\{
\frac{
1 + \frac{2}{\alpha}
}{
1 + \frac{1}{\alpha}
}
\right\}
=
2
\end{equation}
This obtains the expression \eqref{IQF/CQ} from Section \ref{Comparison}. 
Further from \eqref{ENHat} we obtain that 
\begin{align*}
\frac{
\mathbb{E} \hat{N}_{JSQ}
}{
\mathbb{E} \hat{N}_{CQ}
}
&=
\frac{\alpha}{1+\alpha}
\left[ 
	1+ \left[ \frac{1}{\alpha (1+\alpha)^2} + \frac{e^{1+\alpha}}{(1+\alpha)^2} \right]^{-1}\left[  \frac{\alpha^2+4\alpha+1}{\alpha^2 (1+\alpha)^3} + \frac{2 e^{\alpha+1}}{(1+\alpha)^3} \right]
\right]
\\
&=
1- \frac{1}{1+\alpha}
+
\frac{1}{(1+\alpha)^2}
\left[ 
\frac{
\alpha^{2} +4\alpha +1
+2\alpha^2 e^{1+\alpha}
}{
1+\alpha e^{1+\alpha}
}
\right].
\end{align*}
We see that the above expression is continuous for  $\alpha\in (0,\infty)$ and that
\begin{equation*}
\lim_{\alpha\rightarrow \infty}
\frac{
\mathbb{E} \hat{N}_{JSQ}
}{
\mathbb{E} \hat{N}_{CQ}
}	
= 1 
\qquad \text{and} \qquad 
\lim_{\alpha\rightarrow 0}
\frac{
\mathbb{E} \hat{N}_{JSQ}
}{
\mathbb{E} \hat{N}_{CQ}
}	= 1 .
\end{equation*}
Thus this ratio between JSQ and CQ is bounded. 
Due to the combination of exponentials and polynomials there does not appear to be closed form solution for the value of the ratio between JSQ and CQ. However placing in the numerical solver we see that
 \begin{equation*}
\sup_{\alpha >0}
\frac{
\mathbb{E} \hat{N}_{JSQ}
}{
\mathbb{E} \hat{N}_{CQ}
}
\approx 1.13547 	< 1.14
\end{equation*}
and the maximum is achieved at $\alpha =0.209082$ (here we give the first six significant digits). This is the stated expression \eqref{JSQ/CQ} and further is plotted in Figure \ref{fig:NDS_JSQ_IQF}.

\subsection{Proof of Theorem}\label{sec:main_theorem_proof}

What now follows is the proof of Theorem \ref{JSQ:Thrm}. 

\proof{Proof of Theorem \ref{JSQ:Thrm}}
Recall from \eqref{NTerms}, we may write $\hat{N}^{(k)}(t)$ in the form
\begin{subequations}\label{NTerms2}
	\begin{align}
	\hat{N}^{(k)}(t)- \hat{N}^{(k)}(0)  =& {\mathcal M}_a\left( \mu t - \frac{\alpha}{k} \mu t \right) - {\mathcal M}_d\left( \mu t - \frac{\mu}{k^2} \int_0^{kt} I^{(k)}(s)ds  \right) \label{NTerms2a} \\
	& + \frac{\mu}{k^2}\int_0^{k^2t} \hat{I}^{(k)}(s)ds - \mu \int_0^{t}  \frac{(2-\hat{N}^{(k)}(s))_+}{(\hat{N}^{(k)}_1(s)-1 ) } ds \label{NTerms2b}\\
	&- \alpha \mu t + \mu \int_0^{t}  \frac{(2-\hat{N}^{(k)}(s))_+}{(\hat{N}^{(k)}_1(s)-1 ) } ds.\label{NTerms2c}
	\end{align}
\end{subequations}
First consider the right-hand side of \eqref{NTerms2a}. By Corollary \ref{Cor:MtoB} (and the Martingale FCLT) we have that 
\[
\left( {\mathcal M}^{(k)}_a\left( \mu t - \frac{\alpha}{k}\mu t \right) - {\mathcal M}^{(k)}_d\left( \mu t - \frac{\mu}{k^2} \int_0^{kt} I^{(k)}(s)ds  \right)  : t\geq 0\right)\Rightarrow ( \sqrt{2\mu} B(t) : t \geq 0)
\]

Next consider \eqref{NTerms2b}. By Proposition 
\ref{TimescaleIdleness}, for all $\eta >0$
\[
	\sup_{0\leq t\leq \tau^{(k)}(\eta)} 
		\left| 
			\frac{1}{k^2} 
				\int_{0}^{k^2t}	\hat{I}^{(k)}(s)ds 
			- 
				\int_{0}^{t}	
					\frac{(2-\hat{N}^{(k)}(s))_+}{(\hat{N}^{(k)}(s)-1 )}
				ds 
			\right| 
	\xrightarrow[k\rightarrow \infty]{} 
	0.
\]
So defining 
\[
\epsilon^{(k)}(t \wedge \tau^{(k)}(\eta))
=
\frac{1}{k^2} 
		\int_{0}^{k^2(t\wedge \tau^{(k)}(\eta))}	\hat{I}^{(k)}(s)ds 
	- 
		\int_{0}^{(t\wedge \tau^{(k)}(\eta))}	
			\frac{(2-\hat{N}^{(k)}(s))_+}{(\hat{N}^{(k)}(s)-1 )}
				ds ,
\]
we have that 
\[
(\epsilon^{(k)}(t \wedge \tau^{(k)}(\eta)) : t\geq 0, \eta^{-1}\in {\mathbb N}) \Rightarrow 0.
\]
By the Skorohod representation theorem, there exists a probability space were this convergence is almost sure (See \cite{billingsley2013convergence}). That is, almost surely 
\begin{align*}
&\left( \left( {\mathcal M}^{(k)}_a\left( \mu t - \frac{\alpha}{k}\mu t \right) - {\mathcal M}^{(k)}_d\left( \mu t - \frac{\mu}{k^2} \int_0^{kt}\!\! I^{(k)}(s)ds  \right) ,\; \epsilon^{(k)}(t \wedge \tau^{(k)}(\eta))   \right)  : 0\leq t \leq T, \eta^{-1}\in {\mathbb N}\right)\\
&\xrightarrow[k\rightarrow \infty]{u.o.c.} ( (\sqrt{2\mu} B(t),0) : 0\leq t \leq T, \eta^{-1}\in {\mathbb N})	
\end{align*}

Therefore, Lemma \ref{CTS_SDE} applies to give that, almost surely,
\[
(\hat{N}^{(k)}(t): t\geq 0) \xrightarrow[k\rightarrow\infty]{u.o.c.} (\hat{N}(t): t\geq 0) 
\] 
where $(\hat{N}(t): t\geq 0)$ satisfies the S.D.E. 
	\begin{equation}
	d \hat{N} (t) = \mu  \Big[ \frac{(2-\hat{N})_+}{\hat{N}-1}-\alpha\Big] dt + \sqrt{2\mu} \;dB (t), \qquad t\geq 0
	\end{equation}
Here we note that by Proposition \ref{DiffZero} almost surely $\hat{N}(t)>1$ for all $t$, so the conditions of the limit process in Lemma \ref{CTS_SDE} are satisfied. This gives the required weak convergence result. 

The remaining properties on the stationary distribution of $\hat{N}(t)$ are given in Proposition \ref{JSQ:EqPro}. This completes the proof of Theorem \ref{JSQ:Thrm}. \hfill \Halmos
\endproof
\newpage 

\section{Auxiliary Lemmas}\label{sec:auxiliary_lemmas}
This section contains a few somewhat elementary lemmas that will be used repeatedly. 
This Section can be skipped on first reading and referred back to when required.

\begin{lemma}\label{PoissonBound}
	For $Po(\mu)$, a Poisson random variable with parameter $\mu$, it holds that
	\[
	{\mathbb P} \left( Po(\mu) \geq x  \right) \leq e^{-x-\mu },\qquad	\forall x\geq \mu e^2. 
	\]
	and for $\mathcal{N}(t)$ a unit Poisson process
	\begin{align*}		
	&{\mathbb P} \left( \sup_{0\leq t \leq T} (\mathcal{N}(t) -t)\geq z \right) \leq \exp\big\{-\frac{z^2}{2T}\big(1-\frac{z}{T}\big)\big\},\\
		&		{\mathbb P} \left( \inf_{0\leq t \leq T} (\mathcal{N}(t) -t) \leq -z \right) \leq 2\exp\big\{-\frac{z^2}{2T}\big(1-\frac{z}{T}\big)\big\}.
	\end{align*}
	\end{lemma}
	
\proof{Proof}
	A  Chernoff bound gives
	\[
	{\mathbb{P}}\left(Po(\mu) \geq x \right) \leq e^{-\theta x }\mathbb{E}e^{\theta Po(\mu)} 
	= \exp \left\{ -\theta x + \mu \left(e^\theta -1\right)  \right\}.
	\]
	Minimizing over $\theta \geq 0$, the above is minimized at $\theta = \log (x/\mu)$. Therefore 
	\[
	\mathbb{P} \left(Po(\mu) \geq x  \right) \leq \exp \left\{ - \mu \left[ \left( 1- \frac{x}{\mu}\right)+ \frac{x}{\mu}\log\frac{x}{\mu}\right] \right\}.
	\]
	By assumption $\log \frac{x}{\mu}\geq 2$, substituting this into the logarithms in the square brackets above gives
	\[
	\mathbb{P}\left( Po(\mu) \geq x \right) \leq e^{-\mu - x}
	\]
This gives the first bound. The second follows similarly. By Doob's sub-martingale inequality
\[
{\mathbb P} \left( \sup_{0\leq t \leq T} (\mathcal{N}(t) -t) \geq z \right) \leq e^{-\theta (T+z) }\mathbb{E}e^{\theta Po(T)} = \exp \left\{ - T \left[ -\frac{z}{T} + \left(1+ \frac{z}{T} \right)\log\left(1+ \frac{z}{T} \right)\right] \right\}
\]
where we take $\theta =  \log ((z+T)/\mu)$ as above. Applying the bound $\log (1+y) \geq y-{y^2}/{2}$, gives that
\[
- T \left[ -\frac{z}{T} + \left(1+ \frac{z}{T} \right)\log\left(1+ \frac{z}{T} \right)\right] \leq - T \left[ -\frac{z}{T} + \left(1+ \frac{z}{T} \right)\frac{z}{T}\left(1- \frac{z}{2T} \right)\right] = -\frac{z^2}{2T}\big(1-\frac{z}{T}\big).
\]
The same bound holds for the event $\left\{\inf_{0\leq t \leq T} (\mathcal{N}(t) -t) \leq -z\right\}$ by an identical argument.
\hfill \halmos
\endproof


The following is a basic equality on hitting probabilities of Markov chains.
\begin{lemma}\label{HittingLemma}
	Consider a continuous time birth death process taking values $n$ in $\{0,...,x\}$. Let $f(n)$ be the rate of transitions from $n$ to $n+1$, for $n \neq x$, and $g(n)$ be the rate of transitions from $n$ to $n-1$ for $n\neq 0$. All other transition rates from $n$ are zero. The probability of hitting $x$ before $0$ starting from initial state $1$ is given by the expression
	\[
	\frac{1}{\sum_{n=1}^x \prod^{n-1}_{m=1} \frac{g(m)}{f(m)}}. 
	\]
	(here we apply the convention that $\prod_{m=1}^0 \frac{f(m)}{g(m)}=1$  )
\end{lemma}
\proof{Proof}
The $h_n^x$ be the probability of hitting $x$ before zero when starting from state $n$. 
It is clear that $h_n^x$ satisfies the expression
\[
(f(n) + g(n))  h_n^x = f(n) h_{n+1}^x + g(n) h_{n-1}^x 
\]
for $n=1,...x-1$ and $h^x_0=0$ and $h^x_x=1$. Let $u_n = h_n-h_{n-1}$ then the above expression simplifies to give
\[
0=f(n) u_{n+1} - g(n) u_n
\]
for $n=1,...,x-1$, which implies 
\[
u_n = u_1 \prod_{m\leq n-1} \frac{g(m)}{f(m)}
\]
for $n=1,..,x$. Now since 
\[
u_1^x =h_1^x - h_0^x = h_1^x \qquad \text{and} \qquad 1 = h^x_x - h^x_0 = \sum_{n=1}^x u_n ,
\]
we have that
\[
1=h_1^x   \sum_{n=1}^x \prod_{m\leq n-1} \frac{g(x)}{f(x)} 
\]
as required.\hfill \Halmos
\endproof

\begin{lemma}\label{ExcusionLemma}
	Suppose that $Q$ is an M/M/1 queue with arrival rate $\alpha$ and service rate $\beta$, with $\beta > \alpha$. Let $\tilde{T}$ be the length of the renewal cycle from the queue being of length $1$ to the queue being empty again. Then,
	\begin{equation}\label{ExcusionLemma1}
	{\mathbb P} (\tilde{T} \geq t ) \leq e^{\theta + \phi(\theta)t} \qquad \text{and}\qquad {\mathbb P} \left( \int_0^{\tilde{T}} Q(s) ds \geq x\right) \leq e^{-\theta_0 (\sqrt{cx}-1)}
	\end{equation}
	where  $\phi(\theta)=\alpha(e^{\theta}-1) +\beta(e^{-\theta}-1) $, $c>0$
	and  $\theta_0$ is such that $\phi(\theta_0)+ c\theta_0 \leq 0.$
	In particular, $\phi(\theta)$ is minimized by $\theta = \frac{1}{2}\log \frac{\beta}{\alpha}$ and thus
	\begin{equation}\label{ExcusionLemma2}
	{\mathbb P} (\tilde{T} \geq t ) \leq \sqrt{\tfrac{\beta}{\alpha}}e^{-(\sqrt{\beta}- \sqrt{\alpha}) t} .
	\end{equation}
\end{lemma}
\proof{Proof.}
We have that $Q(0)=1$ and that
\[
Q(t) -1 = {\mathcal N}_{\alpha}(t) - {\mathcal N}_{\beta}(t)
\] 
for $t\leq \tilde{T}$ where ${\mathcal N}_\alpha$ and ${\mathcal N}_\beta$ are independent Poisson processes of rates $\alpha$ and $\beta$, respectively.  Note that
\[
\E \left[ e^{\theta ({\mathcal N}_{\alpha}(t) - {\mathcal N}_{\beta}(t))} \right]= e^{t \phi(\theta)}  \quad\text{where}\quad \phi(\theta)=\alpha(e^{\theta}-1) +\beta(e^{-\theta}-1)  .
\]
(Note that, since $\alpha <\beta$, $\phi'(\theta)<0$. So $\phi(\theta)<0$ for some $\theta >0$.)	
Thus 
\[
X(t) = \exp \left\{\theta Q(t\wedge \tilde{T}  )  - (t\wedge \tilde{T})\phi(\theta) \right\}
\]
is a positive martingale which then gives
\[
{\mathbb E} \left[ e^{-\tilde{T}\phi(\theta)} \right] = \bE\left[ X({\infty}) \right]  \leq  {\mathbb E} \left[ X(0) \right] = e^{\theta}.
\]
Applying a Chernoff bound, we get the inequality
\[
{\mathbb P} (\tilde{T} > t) \leq  {\mathbb E} \left[e^{-\phi(\theta) \tilde{T}} \right] e^{\phi(\theta)t} \leq e^{\theta + \phi(\theta)t}.
\]
This gives the required bound on the left-hand side of \eqref{ExcusionLemma1}.
We now optimize $\phi(\theta)$ to gain inequality \eqref{ExcusionLemma2}. In particular, 
\[
0=\phi'(\theta) = \alpha e^{\theta} - \beta e^{-\theta}  \quad \text{which implies} \quad  \theta^* = \frac{1}{2} \log \frac{\beta}{\alpha}.
\]
Substituting back into $\phi$ gives
\[
\phi(\theta^*)= \alpha \left(\sqrt{\frac{\beta}{\alpha}} -1 \right) + \beta \left(\sqrt{\frac{\alpha}{\beta}} -1 \right)=- \left(\sqrt{\beta} -\sqrt{\alpha} \right)^2
\] 
which then gives the required bound
	\[
	{\mathbb P} (\tilde{T} \geq t ) \leq \sqrt{\tfrac{\beta}{\alpha}}e^{-(\sqrt{\beta}- \sqrt{\alpha})^2t} .
	\]

The bound  on the right-hand side of \eqref{ExcusionLemma1} follows by a similar argument, as follows. Define 
\[
Q^c(s)  := Q(s) +cs= 1+{\mathcal N}_{\alpha}(s) - {\mathcal N}_{\beta}(s) + cs \quad \text{and}\quad 
X^c (t) := \exp \left\{\theta Q^c(t)  -\phi(\theta)t-ct\right\}
\]
for some constant $c$ such that $\alpha + c \leq \beta$. 
(Here we have essentially added an additional upward drift to our queueing process.) 
Again,
\[
\E \left[ e^{\theta ({\mathcal N}_{\alpha}(s) - {\mathcal N}_{\beta}(s) + cs) }  \right]= e^{s (\phi(\theta)+c\theta)}  \quad\text{where}\quad \phi(\theta)=\alpha(e^{\theta}-1) +\beta(e^{-\theta}-1).  
\]
Let $\tilde{T}_z$ be the first time that $Q^c(t)\geq z$ holds, for $z\geq 1$. For any $\theta$ with $\phi(\theta)+c\theta \geq 0$, $X^c({t\wedge \tilde{T}_z})$ is a positive martingale bounded above by $e^{\theta z}$. By the Optional Stopping Theorem
\[
e^\theta = {\mathbb E} [ X^c(0)] = {\mathbb E} [X^c({\tilde{T}_z})] = {\mathbb E} [ e^{ \theta z - \tilde{T}_z (\phi(\theta) + c\theta)}  {\mathbb I} [ \tilde{T}_z < \infty] ].
\]
Let $\theta^*$ be the largest solution to the equation $\phi(\theta^*) + c\theta^*=0$. We have that 
\[
{\mathbb P} (\tilde{T}_z<\infty) \leq e^{-\theta^* (z-1)} 
\]
and thus for any $\theta_0\leq \theta^*$
\[
{\mathbb P} (\tilde{T}_z<\infty) \leq e^{-\theta_0 (z-1)}.
\]
Note that since $\phi(\theta) +c \theta$ is a convex function, zero at $\theta=0$ with $\phi'(0)+ c <0$. Thus $\theta_0\leq \theta*$ iff $\phi(\theta_0)+ c\theta_0 \leq 0$.

Now if $\tilde{T}_z =\infty$ then $Q(t)$ lies below the line $(z-cs : s\geq 0)$. Thus, the area under the $Q(t)$ before it hits zero is less that the area of the triangle in the positive orthant with hypotenuse $(z-cs : s\geq 0)$. That is
\[
\int_0^{\tilde{T}} Q(t) dt < \frac{z^2}{c} .
\] 
Therefore
\[
{\mathbb P} \Big( \int_0^{\tilde{T}}   Q(t) dt \geq  \frac{z^2}{c}\Big) \leq \mathbb{P}(\tilde{T}_z <\infty) \leq e^{-\theta_0 (z-1)}
\]
which after substituting $x= z^2/c$ gives the required bound \eqref{ExcusionLemma1}.

\hfill \Halmos\endproof

\begin{lemma}\label{ExcusionLemma4}
	Suppose that $Q$ is an M/M/1 queue with arrival rate $\alpha$ and service rate $\beta$, with $\beta > \alpha$ let $\tilde{T}$ be the length of the excursion of the $M/M/1$ queue started from zero then 
	\[
	\mathbb{E} \left[ \int_0^{\tilde{T}}\!\! \Big( Q(s) - \frac{\alpha}{\alpha-\beta} \Big) ds \right] = 0.
	\]
\end{lemma}
\proof{Proof}
The stationary distribution of $Q$ (or indeed any irreducible Markov chain) can be expressed as 
\[
\mathbb{P}_\pi (Q=q) 
	=
	\frac{
		1
		}{
		\mathbb{E} [\tilde{T}]  
		}
	\mathbb{E}  
	\left[  
		\int_0^{\tilde{T}} 
			\mathbb{I} [ Q(t) = q]
		dt
	\right] .
\]
where here $\mathbb{P}_\pi$ is the stationary distribution of $Q$. In other words, since the stationary distribution of $Q$ is geometric, we have that
\begin{equation*}
	\mathbb{E} \left[
		\int_0^{\tilde{T}} 
			\Big(
				1- \frac{\alpha}{\beta}
			\Big) 
			\Big(
				\frac{\alpha}{\beta}
			\Big)^{q}
		dt
		\right]  
	=
	\mathbb{E}  
	\left[  
		\int_0^{\tilde{T}} 
			\mathbb{I} [ Q(t) = q]
		dt
	\right] .
\end{equation*}
Multiplying by $q$ and summing over $q\in \mathbb{Z}_+$ gives the result. \hfill \Halmos 
\endproof

The following is a Functional Strong Law of Large Numbers for $L^2$ martingales and is an extension of \cite{williamsprobability} Theorem 12.13a) and Section 12.4.

\begin{lemma}[A Martingale Functional Strong Law of Large Numbers]\label{MgSLLN} Suppose that 
	\[
	M_n:= \sum_{i=1}^{n} Z_i
	\]
	defines a Martingale and is such that, for some constant $C$,
	\[
	\mathbb{E}\left[ Z_k^2 | \mathcal{F}_{k-1} \right] \leq C
	\]
	then, for all $T>0$,
	\[
	\lim_{n\rightarrow \infty}\sup_{0\leq t \leq T } \left| \frac{M_{\lfloor tn \rfloor }}{n} \right| =0.
	\]	
\end{lemma}
\proof{Proof}
By assumption, the increasing process
\[
<\! M\! >_n := \sum_{k \leq n} \mathbb{E} \left[ Z_k^2 | \mathcal{F}_{k-1} \right]  
\]
is such that $<\! M\! >_n\leq nC$. By \cite[Theorem 12.13a and 12.4]{williamsprobability}  
\begin{align*}
&\lim_{n\rightarrow \infty} \frac{M_n}{<\! M\! >_n} =0 \qquad \text{on the event}\qquad \left\{<\! M\! >_n=\infty\right\}, \\
&\lim_{n\rightarrow \infty} M_n \quad \text{exists and is finite on the event }\quad \left\{<\! M\! >_n<\infty\right\} .
\end{align*} 
Since $<\! M\! >_n\leq nC$, in both instances we have that
\[
\lim_{n\rightarrow \infty}\frac{M_n}{n} =0.
\]

Expanding this statement, 
\begin{equation*}
\forall \epsilon >0,\quad \exists \tilde{N} \quad\text{s.t.} \quad \forall n> \tilde{N}, \qquad \left| \frac{M_n}{n}  \right| \leq \epsilon.
\end{equation*}
Notice, this also bounds our process for all suitably large times:
\begin{equation*}
\forall \epsilon >0,\quad \exists \tilde{N} \quad\text{s.t.} \quad \forall t n> \tilde{N} ,\qquad \left| \frac{M_{\lfloor tn \rfloor}}{n}  \right| \leq \epsilon t.  
\end{equation*}
This just leaves a finite number of terms, that is $\tilde{N}$ terms, to deal with. We apply the crude bound, if $tn \leq \tilde{N}$
\begin{equation*}
\left| \frac{M_{\lfloor tn \rfloor}}{n}  \right| \leq \frac{1}{n} \max_{k=1,...,\tilde{N}}  |M_k|
\end{equation*}
Notice, we can make the right-hand side of this expression small by taking $n$ much larger that $\tilde{N}$. Thus we can collect together these two cases
\begin{equation*}
\left| \frac{M_{\lfloor tn \rfloor}}{n}  \right| \leq  
\begin{cases}
\epsilon t, & \text{if } tn > \tilde{N}\\
\frac{1}{n}\max_{k=1,...,\tilde{N}}  |M_k|, &\text{if } tn \leq \tilde{N}
\end{cases}
\end{equation*}
By bounding $t\leq T$ the first case can be made arbitrarily small and, then, by choosing $n$ suitably large compared to $\tilde{N}$ the second case can be made arbitrarily small. This then proves that 
\[
\lim_{n\rightarrow \infty}\sup_{0\leq t \leq T } \left| \frac{M_{\lfloor tn \rfloor }}{n} \right| =0
\]	
as required.\hfill \halmos
\endproof

The following result is the central limit theorem counterpart of the above strong law of large numbers result. For it's proof, we refer the interested reader to the excellent survey of \cite{whitt2007}.

\begin{theorem}[Martingale Functional Central Limit Theorem]\label{MgCLT}
If $M^{(k)}=(M^{(k)}(t) : t\geq 0)$, $k\geq 0$ is sequence of continuous time Martingales whose quadratic variation $[M^{(k)}]$ converge in the Skorohod topology to the identify function:
\[
[M^{(k)}] \Rightarrow \mathbf{e}\qquad\text{as} \quad k\rightarrow \infty
\]
then 
\[
M^{(k)}\Rightarrow B
\]
where $B$ is a standard Brownian motion.
\end{theorem}
For a proof of this result see \cite[Theorem 2.1 and Lemma 2.1]{whitt2007}. The proof of the Martingale Functional Central Limit consists of verifying tightness of the Martingale sequence and then verifying L\'evy's characterization of Brownian motion for any limit of this relatively compact sequence.

\subsection{Comparison of stationary distributions}
In Section \ref{Comparison}, we state that the following stochastic bounds hold
\[
	\pi_{CQ} \leq_{st} \pi_{JSQ} =_{st} \pi_{I1F} \leq_{st} \pi_{IQF},
\]
and we state that this follows since the drift terms of each stochastic process, corresponding to $CQ$, $JSQ$, $I1F$ and $IQF$ dominate each other. This result follows from the following lemma. 

\begin{lemma}\label{St_Dom}
Consider positive recurrent stochastic differential equations 
\begin{align*}
	dX_1(t) &= - d_1(X_1(t)) dt + \sqrt{2} dB(t) \\
	dX_2(t) &= - d_2(X_2(t)) dt + \sqrt{2} dB(t) 	
\end{align*}
with stationary distributions $\pi_1(x)$ and $\pi_2(x)$ with support on $(0,\infty)$. 
If the drift of these processes are continuous on $(0,\infty)$ and satisfy $d_1(x)\geq d_2(x)$ then
\[
\pi_1 \leq_{st} \pi_2\ .
\]
\end{lemma}

(Within the proof, we provide comments as to how the result extends to a reflected Brownian motion with drift, as is given for the Central Queue NDS diffusion.)
\proof{Proof}
The diffusions $X_1$ and $X_2$ are Langevin Diffusions. Specifically, their generators and forward-equation are of the form
\[
\mathcal L  = \frac{d^2}{d x^2} - V'(x) \frac{d}{dx}
\]
for some differentiable function $V(x)$. 
Setting $V'_1(x) = d_1(x)$ and $V'_2(x)=d_2(x)$ 
This can be verified directly with It\^{o}'s formula. The invariant measure for such diffusions has density given by $m(x) = e^{-V(x)}$. This holds  since $\mathcal L^\dagger m(x) = 0 $ which can be verified substitution into the adjoint of the generator above, or by substitution into the Kolmogorov forward-equation. (For reflected Brownian motion with drift the same form of invariant measure holds. So our results apply here also.)

We now directly calculate and compare the stationary distributions for these processes. First note that
\[
-V_1(y) +V_1(x)=\int_x^y -d_1(u) {\rm d}u \leq \int_x^y -d_2(u) {\rm d}u = -V_2(y) + V_2(x)\ .
\]
Taking exponentials and integrating once more gives
\begin{equation}\label{Vineq}
\int_x^\infty e^{-V_1(y) + V_1(x)} {\rm d}y \leq \int_x^\infty e^{-V_2(y) + V_2(x)} {\rm d}y\ .
\end{equation}

Now note that
\[
\mathbb P_{\pi_1} ( X_1 \geq x) 
= 
\frac{
\int_x^\infty e^{-V_1(y)} dy 
}{
\int_0^\infty e^{-V_1(y)} dy
}.
\]
Taking logarithms and differentiating with respect to $x$ gives
\begin{align*}
	\frac{d}{dx} \log \mathbb P_{\pi_1} ( X_1 \geq x)  
	&
	= \int_x^\infty e^{-V_1(y) + V_1(x)} {\rm d}y
	\\
	&
	\leq  
	 \int_x^\infty e^{-V_2(y) + V_2(x)} {\rm d}y
	 =
	 \frac{d}{dx} \log \mathbb P_{\pi_2} ( X_2 \geq x)  \ .
\end{align*}
In the inequality above, we apply \eqref{Vineq}.  Thus after integrating and taking exponentials, we have that
\[
\mathbb P_{\pi_1} ( X_1 \geq x)   \leq \mathbb P_{\pi_2} ( X_2 \geq x)  
\]
Thus are required we have that $\pi_1 \leq_{st} \pi_2$.
\hfill \Halmos
\endproof

\end{APPENDICES}

%
%






\end{document}